\documentclass[aap,seceqn,MSNbibl,dvips]{arximspdf}
\usepackage{mathrsfs}

\doi{10.1214/09-AAP649}
\volume{20}
\issue{3}
\pubyear{2010}
\firstpage{1126}
\lastpage{1176}

\makeatletter

\newtheorem{theorem}{Theorem}[section]
\newtheorem{prop}{Proposition}[section]
\newtheorem{col}{Corollary}[section]
\newtheorem{lem}{Lemma}[section]

\makeatother

\begin{document}
\begin{frontmatter}

\title{Network stability under max--min fair bandwidth~sharing}
\runtitle{Network stability under max--min fair bandwidth sharing}

\begin{aug}
\author[A]{\fnms{Maury} \snm{Bramson}\corref{}\thanksref{t1}\ead[label=e1]{bramson@math.umn.edu}}
\runauthor{M. Bramson}
\affiliation{University of Minnesota}
\address[A]{School of Mathematics\\
University of Minnesota\\
Twin Cities Campus\\
Institute of Technology\\
127 Vincent Hall\\
206 Church Street S.E.\\
Minneapolis, Minnesota 55455\\
USA\\
\printead{e1}} 
\end{aug}

\thankstext{t1}{Supported in part by NSF Grants DMS-02-26245 and CCF-0729537.}

\received{\smonth{12} \syear{2008}}
\revised{\smonth{7} \syear{2009}}

%
\begin{abstract}
There has recently been considerable interest in the stability of
different fair bandwidth sharing policies for models that arise in the
context of Internet congestion control. Here, we consider a connection
level model, introduced by Massouli\'{e} and Roberts
[\textit{Telecommunication Systems} \textbf{15} (2000) 185--201], that
represents the randomly varying number of flows present in a network.
The weighted $\alpha$-fair and weighted max--min fair bandwidth sharing
policies are among important policies that have been studied for this
model. Stability results are known in both cases when the interarrival
times and service times are exponentially distributed. Partial results
for general service times are known for weighted $\alpha$-fair
policies; no such results are known for weighted max--min fair policies.
Here, we show that weighted max--min fair policies are stable for
subcritical networks with general interarrival and service
distributions, provided the latter have $2+\delta_{1}$ moments for some
$\delta_{1}>0$. Our argument employs an appropriate Lyapunov function
for the weighted max--min fair policy.
\end{abstract}

%

\end{frontmatter}

\section{Introduction}\label{intro}

We consider a connection level model for Internet congestion control
that was first studied by Massouli\'{e} and Roberts \cite{r9}. This
stochastic model represents the randomly varying number of flows in a
network for which bandwidth is dynamically shared among flows that
correspond to the transfer of documents along specified routes.
Standard bandwidth sharing policies are the weighted $\alpha$-fair,
$\alpha\in(0, \infty)$, and the weighted max--min fair policies. An
important example of the former is the proportionally fair policy,
which corresponds to $\alpha= 1$. The weighted max--min fair policy
corresponds to $\alpha= \infty$. These policies allocate service
uniformly to documents along a given route, and allocate service
amongst different routes in a ``fair'' manner. A question of
considerable interest is when such policies are stable.

De Veciana, Lee and Konstantopoulos \cite{r5} studied the stability of
weighted max--min fair and proportionally fair policies; Bonald and
Massouli\'{e} \cite{r2} studied the stability of weighted $\alpha$-fair
policies. Both papers assumed exponentially distributed interarrival
and service times for documents. The first condition is equivalent to
Poisson arrivals, and the second condition corresponds to exponentially
distributed document sizes with documents processed at a constant rate.
Both papers constructed Lyapunov functions which imply the stability of
such models when the models are subcritical, that is, the underlying
Markov process is positive Harris recurrent when the average load at
each link is less than its capacity.

Relatively little is currently known regarding the stability of
subcritical networks with general interarrival and service times.
Massouli\'{e} \cite{r8} showed stability for the proportionally fair
policy for exponentially distributed interarrival times and general
service times that are of phase type. A suitable Lyapunov function was
employed to show stability.

The stability problem for bandwidth sharing policies is in certain
aspects similar to the analogous problem for multiclass queueing
networks. A significant complication that arises in the context of
bandwidth sharing policies is the requirement of simultaneous service
of documents at all links along a route. This can reduce the efficiency
of service, and complicates analysis when the interarrival and service
times are not exponentially distributed.

When the interarrival and service times are exponentially distributed,
finer results are possible. In Kang et al. \cite{r7}, a diffusion
approximation is established under weighted proportionally fair
policies. There and in Gromoll and Williams \cite{r6}, summaries and a
more detailed bibliography are provided for different bandwidth sharing
policies, for both exponentially distributed and more general
interarrival and service times.

Here, we investigate the behavior of weighted max--min fair policies
for subcritical networks whose interarrival and service times have
general distributions. We show that such networks are stable, provided
that the service distributions have $2+\delta_{1}$ moments for some
$\delta_{1}>0$. No conclusion is reached when fewer moments exist. As
in previous papers on stability, we construct a suitable Lyapunov
function. Because of the more general framework here, the Markov
process underlying the model will now have a general state space, and
will require the machinery associated with positive Harris recurrence.

We next give a more detailed description of the model we consider,
after which we state our main results. We then provide some basic
motivation behind their proof together with a summary of the remainder
of the paper.

\subsection*{Description of the model}

In the model we consider, \textit{documents} are assumed to arrive at
one of a finite number of \textit{routes} $r \in\mathcal{R}$ according
to independent renewal processes, with interarrival times denoted by
$\xi_{r}(1), \xi_{r}(2), \ldots.$ Here, $\xi_{r}(1)$ are the initial
residual interarrival times, and are considered part of the initial
state. The remaining variables $\xi_{r}(2), \xi_{r}(3), \ldots$ are
assumed to be i.i.d. with mean $1/\nu_{r}$, $\nu_{r}>0$, for each $r$,
with the sequences being independent of one another; $\xi_{r}$ will
denote a random variable with the corresponding distribution. The
service times of documents are assumed to be independent of one another
and of the interarrival times, and have distribution functions
$H_{r}(\cdot)$ with means $m_{r}<\infty$. The initial state will
include the residual service times of documents initially in the
network.

On each route $r$, there are a finite number of links $l$, where
service is allocated to the documents on the route. For the models
considered in \cite{r2,r5,r8} and \cite{r9}, documents on a route $r$
receive service simultaneously at all links $l$ on the route, with all
such documents being allocated the same rate of service $\lambda_{r}$
at all such links at a given time. Associated with such a network is an
\textit{incidence matrix} $A=(A_{l,r})$, $l \in\mathcal{L}$, $r
\in\mathcal{R}$, with $A_{l,r} = 1$ if link $l$ lies on route $r$ and
$A_{l,r} = 0$ otherwise. When $A_{l,r} = A_{l,r'} = 1$, with $r \ne
r'$, the routes $r$ and $r'$ share a common link.

Setting $z_{r}$ equal to the number of documents on route $r$, $\Lambda
_{r} = \lambda_{r}z_{r}$ denotes the rate of service allocated to the
totality of all documents on the route. Each link $l$ is assumed to
have a given \textit {bandwidth capacity} $c_{l}>0$. A \textit{feasible
policy} requires that this capacity not be exceeded, namely
%
%
\begin{equation}
\label{eq1.12.1}
\sum_{r\in\mathcal{R}}A_{l,r} \Lambda_{r} \le c_{l}\qquad
\mbox{for all } l \in\mathcal{L}.
\end{equation}
Denoting by $\Lambda= (\Lambda_{r})$ and $c=(c_{l})$ the corresponding
column vectors, this is equivalent to $A\Lambda\le c$, with the
inequality being interpreted coordinatewise.

None of the results in this paper relies on the restriction that either
$A_{l,r} = 1$ or $A_{l,r}=0$. Here, we continue to assume that
(\ref{eq1.12.1}) is satisfied, for given $A$, but with the weaker
assumption $A_{l,r}\ge 0$. Under this new setup, each link may be
interpreted as belonging to every route. A given link $l$ now allocates
the same rate of service $\Lambda_{r}$ to each route $r$, which
utilizes this service at rate $A_{l,r}$. For $A_{l,r}\in[0,1]$,
$A_{l,r}$ may be interpreted as the proportion of this potential
service that is actually utilized at link $l$ by route $r$.

The \textit{traffic intensity} $\rho_{r} = \nu_{r}m_{r}$ measures the
average rate over time at which work enters a route $r$. We say a
network is \textit{subcritical} if
%
%
\begin{equation}
\label{eq5.3.8}
\sum_{r\in\mathcal{R}}A_{l,r}\rho_{r}<c_{l} \qquad\mbox{for
all }l\in \mathcal{L},
\end{equation}
or, in matrix form, $A\rho<c$, where $\rho=(\rho_{r})$ is the
corresponding column vector. This corresponds to the definition of
subcriticality that is employed in the context of queueing networks,
where the load at each station (here, load at each link) is strictly
less than its capacity. Condition (\ref{eq5.3.8}) is needed for
stability. It is assumed in, for example, \cite{r2,r5} and \cite{r8}.

The $\alpha$-fair and max--min fair policies are examples of feasible
policies for which the allocation of service to documents at a given
time is determined by the vector $z=(z_{r})$; the weighted
$\alpha$-fair and max--min fair policies are defined analogously, but
with a weight $w_{r}>0$ assigned to route $r$. We do not define
$\alpha$-fair here, or, in particular, proportionally fair, referring
the reader to the previous references. \textit{Weighted max--min fair}
(WMMF) is defined as a feasible policy that, at each time, allocates
service so that
%
%
\begin{equation}
\label{eq1.15.1}
\min_{r \in\mathcal{R}'} \{\lambda_{r}/w_{r}\}
\mbox{ is maximized},
\end{equation}
among nonempty routes $\mathcal{R}'$. That is, the minimum amount of
weighted service each document receives is maximized, on $r$ with
$z_{r}>0$, subject to the constraint (\ref{eq1.12.1}).

As defined above, a WMMF policy always exists, although it need not be
unique, since there may be some flexibility in allocating service among
those routes where documents are receiving more than the minimal amount
of service. Since our results apply to all such policies, we will not
bother to select a ``best'' member that, for instance, maximizes
service on the routes that are already receiving more than the minimum
service. Such a ``best'' policy can be obtained by solving a hierarchy
of optimization problems, as mentioned above display (2) in~\cite{r5}.
[By employing the convexity that is inherent in the constraint
(\ref{eq1.12.1}), it is routine to verify the existence of such
policies.]

Since the vector $z$ of documents changes as time evolves, so will the
allocation of service. From this point on, we reserve the notation
$\lambda_{r}(t)$ and $\Lambda_{r}(t)$ for the allocation of service for
a WMMF policy at time $t$. We find it useful to also introduce
%
%
\begin{equation}
\label{eq1.17.1}
\lambda^{w}(t) = \min_{r \in\mathcal{R}'}
\{\lambda_{r}(t)/w_{r}\}
\end{equation}
with (\ref{eq1.15.1}) in mind. Between arrivals and departures of
documents, $\lambda(\cdot) = (\lambda_{r}(\cdot))$ and
$\lambda^{w}(\cdot)$ will be constant; we specify that they be right
continuous with left limits.

The state of the network evolves over time as documents arrive in the
network, are served, and then depart. For networks with exponentially
distributed interarrival and service times and an assigned policy,
$z=(z_{r})$ suffices to describe its state. As with queueing networks,
one needs to specify the residual interarrival and service times in
general. With this in mind, we employ the notation $z_{r}(B_{r})$ to
denote the number of documents on route $r$ that have residual service
times in $B_{r}\subseteq\mathbb{R}^{+}$, and $u_{r}$ to denote the
residual interarrival time for $r$, with $z(B)=(z_{r}(B_{r}))$,
$B=(B_{r})$, and $u=(u_{r})$ denoting the corresponding vectors.
Setting
%
%
\begin{equation}
\label{eq1.18.1}
x = (z(\cdot), u),
\end{equation}
the state $x$ contains this information. We will employ $X(t),
Z(t,\cdot)$ and $U(t)$ for the corresponding random states at time $t$.
The natural metric space $S$ that corresponds to the states $x$ is no
longer discrete. We will describe $S$ in more detail in Section
\ref{sec2}.

One can specify a Markov process $X(\cdot)$ on $S$ that corresponds to
the network with the assigned WMMF policy. The process $X(\cdot)$ is
constructed in the same manner as is its analog for a queueing network.
More detail is again given in Section \ref{sec2}. We note here that
since $S$ is not discrete, the notion of positive recurrence needs to
be replaced by that of positive Harris recurrence. When $X(\cdot)$ is
positive Harris recurrent, we will say that the network is
\textit{stable}.

In order to demonstrate positive Harris recurrence for $X(\cdot)$, we
will define, in Section \ref{sec3}, an appropriate nonnegative
function, or \textit{norm}, $\|x\|$, for $x \in S$. It is defined in
terms of the norms $|x|_{L}, |x|_{R}$ and $|x|_{A}$, by
%
%
\begin{equation}
\label{eq1.20.1}
\|x\| = |x|_{L} + |x|_{R} + |x|_{A}.
\end{equation}
Without going into detail here, we note that $|x|_{L}$ and $|x|_{R}$
are defined from $z(\cdot)$, where $|x|_{L}$, in essence, measures
residual service times smaller than $N$, for a given large $N$,
$|x|_{R}$ measures residual service times greater than $N$, and
$|x|_{A}$ is a function of the largest residual interarrival time.
(When a distribution function $H_{r}$ has a thin enough tail, we
actually replace $N$ by a smaller value $N_{H_{r}}$ that depends on
$H_{r}$.) As one should expect, as either the total number of documents
$\sum_{r}z_{r} \to\infty$ or $|u| \to\infty$, then $\|x\| \to \infty$.

\subsection*{Main results} We now state our two main results.
\begin{theorem}
\label{thm5.7.1}
Suppose that a subcritical network with a weighted
max--min fair policy has interarrival times with finite means and
service times with $2+\delta_{1}$ moments, $\delta_{1}>0$. For the norm
in (\ref{eq1.20.1}), and appropriate $N,L$ and $\varepsilon_{1}> 0$,
%
%
\begin{equation}
\label{eq5.3.10}
E_{x}[\|X(N^{3})\|] \le(\|x\| \vee L) - \varepsilon_{1}N^{2}\qquad
\mbox{for all } x \in S.
\end{equation}
\end{theorem}

Inequality (\ref{eq5.3.10}) states that, for large $\|x\|$, $X(\cdot )$
has an average negative drift over $[0,N^{3}]$ that is at least of
order $1/N$. This rate will be a consequence of the application of $N$
in the construction of the norm $|x|_{L}$ that appears in
(\ref{eq1.20.1}).

The reader will recognize (\ref{eq5.3.10}) as a version of Foster's
criterion. It will imply the positive Harris recurrence of $X(\cdot)$,
provided that the states in $S$ communicate with one another in an
appropriate sense. Petite sets are typically employed for this purpose;
they will be defined in Section \ref{sec2}. A petite set $A$ has the
property that each measurable set $B$ is ``equally accessible'' from
all points in $A$ with respect to a given measure.
\begin{theorem}
\label{thm1.24.1}
Suppose that a subcritical network with a weighted
max--min fair policy has interarrival times with finite means and
service times with $2+\delta_{1}$ moments, $\delta_{1}>0$. Also,
suppose that $A_{L}=\{x\dvtx\|x\|\le L\}$ is petite for each $L>0$, for
the norm in (\ref{eq1.20.1}). Then, $X(\cdot)$ is positive Harris
recurrent.
\end{theorem}

Theorem \ref{thm1.24.1} will follow from Theorem \ref{thm5.7.1} by
standard reasoning. More detail is given in Section \ref{sec2}.

A standard criterion that ensures the above sets $A_{L}$ are petite is
given by the following two conditions on the interarrival times. The
first condition is that the distribution of $\xi_{r}(2)$ is unbounded
for all $r$, that is,
%
%
\begin{equation}
\label{eq1.25.1}
P\bigl(\xi_{r}(2) \ge s\bigr) > 0 \qquad\mbox{for all } s.
\end{equation}
The second condition is that, for some $l_{r}\in\mathbb{Z}^{+}$, the
$(l_{r}-1)$-fold convolution of $\xi_{r}(2)$ and Lebesque measure are
not mutually singular. That is, for some nonnegative $q_{r}(\cdot)$
with $\int^{\infty }_{0}q_{r}(s)\,ds>0$,
%
%
\begin{equation}
\label{eq1.25.2}
P\bigl(\xi_{r}(2) + \cdots+ \xi_{r}(l_{r}) \in[c,d]\bigr)
\ge\int^{d}_{c} q_{k}(s)\,ds
\end{equation}
for all $c<d$. When the interarrival times are exponentially
distributed, both (\ref{eq1.25.1}) and (\ref{eq1.25.2}) are immediate.
More detail is given in Section \ref{sec2}.

We therefore have the following corollary of Theorem \ref{thm1.24.1}.
\begin{col}
\label{col1.26.1}
Suppose that a subcritical network with a weighted
max--min fair policy has interarrival times with finite means that
satisfy (\ref{eq1.25.1}) and (\ref{eq1.25.2}), and service times with
$2+\delta_{1}$ moments, $\delta_{1}>0$. Then, $X(\cdot)$ is positive
Harris recurrent.
\end{col}

\subsection*{Outline of the paper and main ideas}

In Section \ref{sec2}, we will provide a brief background of Markov
processes and will summarize the construction of the space $S$ and
Markov process $X(\cdot )$ described above. We will also provide
background that will be employed to derive Theorem \ref{thm1.24.1} from
Theorem \ref{thm5.7.1} and to obtain Corollary~\ref{col1.26.1}. The
machinery for this is standard in the context of queueing networks; we
explain there the needed modifications.

The remainder of the paper is devoted to the demonstration of Theorem
\ref{thm5.7.1}. (One minor result, Proposition \ref{prop3.13.1}, is
needed for Theorem \ref{thm1.24.1}.) In Section \ref{sec3}, we will
specify the norms \mbox{$|\cdot|_{L}$}, \mbox{$|\cdot|_{R}$} and
\mbox{$|\cdot|_{A}$} that define \mbox{$\|\cdot\|$} in
(\ref{eq1.20.1}). Employing bounds on these three norms that will be
derived in Sections \ref{sec4}, \ref{sec5} and \ref{sec10}, we obtain
the conclusion (\ref{eq5.3.10}) of Theorem \ref {thm5.7.1}.

For large $\|x\|$, at least one of the norms $|x|_{V}$, with $V$ equal
to $L, R$ or $A$, must also be large. When $|x|_{V}$ is large for given
$V$, it will follow that $E_{x}[|X(N^{3})|_{V}]-|x|_{V}$ is
sufficiently negative so that (\ref{eq5.3.10}) will hold.

The analysis for $|\cdot|_{A}$ is straightforward and is done in
Section \ref{sec4}. The behavior of $E_{x}[|X(N^{3})|_{R}]-|x|_{R}$ is
analyzed in Section \ref {sec5}. The remaining five sections are
devoted to analyzing $E_{x}[|X(N^{3})|_{L}] - |x|_{L}$. In the last two
cases, one needs to reason that, in an appropriate sense, the decrease
in residual service times of existing documents more than compensates
for the increase due to arriving documents, thus producing a net
negative drift.

For such an analysis, it makes sense to decompose the process $X(\cdot
)$ into pro\-cesses $\tilde{X}(\cdot)$ and $X^{A}(\cdot)$, with
\[
X(t) = \tilde{X}(t) + X^{A}(t) \qquad\mbox{for all } t.
\]
The process $\tilde{X}(t)$ is obtained from $X(t)$ by retaining only
those documents, the \textit{original documents}, that were initially
in the network, and $X^{A}(t)$ consists of the remaining documents.
Neither $\tilde {X}(\cdot)$ nor $X^{A}(t)$ is Markov. One defines
$\tilde{Z}(t,B)$ and $Z^{A}(t,B)$ analogously to $Z(t,B)$.

Because of the WMMF policy, all documents that remain on a route $r$,
over the time interval $[0,t]$, receive the same service
$\Delta_{r}(t)$, with $\Delta_{r}(t) = \int^{t}_{0}\lambda_{r}(t')\,dt'$.
Consequently,
%
%
\begin{equation}
\label{eq1.30.1}
\tilde{Z}_{r}(t,B) = z_{r}\bigl(B+\Delta_{r}(t)\bigr) \qquad\mbox{for }
t \ge0, r \in \mathcal{R}, B \subseteq\mathbb{R}^{+}.
\end{equation}
The norms $|\cdot|_{L}$ and $|\cdot|_{R}$ will be defined so that
documents with greater residual service times contribute more heavily
to the norms. On account of (\ref{eq1.30.1}), $|\tilde{X}(t)|_{L}$ and
$|\tilde{X}(t)|_{R}$ will therefore decrease over time; one can also
obtain upper bounds on $|X^{A}(t)|_{L}$ and $|X^{A}(t)|_{R}$. One can
use this to obtain a negative net drift on
$E_{x}[|X(N^{3})|_{L}]-|x|_{L}$ and $E_{x}[|X(N^{3})|_{R}]-|x|_{R}$, as
mentioned earlier.

Only limited use of inequalities arising from (\ref{eq1.30.1}) is
needed in Section \ref{sec5} for $|\cdot|_{R}$. More detailed versions
are needed for $|\cdot|_{L}$, which are presented in the first part of
Section \ref{sec6}.

In Section \ref{sec6}, we also introduce the sets $\mathcal{A}(t)$,
along which we will be able to obtain good pathwise upper bounds on
$|X^{A}(t)|_{L}$. We show in Section \ref{sec6}, by using elementary
large deviation estimates, that the probabilities of the complements
$\mathcal{A}(t)^{c}$ are small enough so that
\[
E_{x} [|X(N^{3})|_{L} - |x|_{L}; \mathcal{A}(N^{3})^{c} ]
\]
is negligible with respect to $E_{x}[|X(N^{3})|_{L}]-|x|_{L}$.

Sections \ref{sec7}--\ref{sec10} analyze the behavior of
$|X(N^{3})|_{L}$ on $\mathcal {A}(N^{3})$. Section \ref{sec7} considers
the contribution to $|X(N^{3})|_{L}$ of residual times $s > N_{H_{r}}$;
$N_{H_{r}}$ was mentioned parenthetically after (\ref {eq1.20.1}) and
satisfies $N_{H_{r}} \le N$. Sections \ref{sec8} and \ref{sec9}
consider the contribution to $|X(N^{3})|_{L}$ of residual times $s \le
N_{H_{r}}$. In Section \ref{sec8}, this is done for $\Delta_{r}(N^{3})
> 1/b^{3}$, for given $r$, with the constant $b$ introduced in
(\ref{eq5.1.4}). Here, service of individual documents is intense
enough to provide straightforward upper bounds for
$|X(N^{3})|_{L}-|x|_{L}$.

Section \ref{sec9} considers the case with $\Delta_{r}(N^{3})
\le1/b^{3}$. This is the only place in the paper where the
subcriticality of the network is employed; estimation for
$|X(N^{3})|_{L}$ must therefore be more precise. The short Section
\ref{sec10} combines the results of Sections \ref{sec6}--\ref {sec9} to
give the desired bounds on $E_{x}[|X(N^{3})|_{L}]-|x|_{L}$.

\subsection*{Notation}
For the reader's convenience, we list here some of the notation in the
paper, part of which has already been employed. We set $\bar{H}_{r}(s)
= 1 - H_{r}(s)$; quantities such as $\bar{H}^{*}_{r}(s)$ and
$\bar{\Phi}^{*}_{r}(s)$, are defined analogously in terms of
$H^{*}_{r}(s)$ and $\Phi^{*}_{r}(s)$, which will be introduced later
on. The term $x$ indicates a state in $S$ and the corresponding term
$X(t)$ indicates a random state at time $t$; $z(\cdot)$ and
$Z(t,\cdot)$, and $u$ and $U(t)$ play analogous roles. We will
abbreviate $\Delta _{r} = \Delta_{r}(N^{3})$ and set $i_{r}(s) = s +
\Delta_{r}$; $i_{r}(s)$ is the initial residual service time of an
original document that has residual service time $s$ at time $N^{3}$.
We employ $C_{1}, C_{2}, \ldots$ and $\varepsilon_{1}, \varepsilon_{2},
\ldots$ for different positive constants that appear in our bounds,
whose precise values are unimportant. The symbols $\mathbb{Z}^{+}$ and
$\mathbb{R}^{+}$ denote the positive integers and positive real
numbers, and $\mathbb{Z}^{+,0} = \mathbb{Z}^{+} \cup\{0\}$; $\lfloor y
\rfloor$ and $\lceil y \rceil$ denote the integer part of
$y\in\mathbb{R}^{+}$ and the smallest integer $n$ with $n \ge y$; and
$c\vee d$ and $c\wedge d$ denote the greater and smaller value of $c,d
\in\mathbb{R}$. The acronyms LHS and RHS will stand for ``left-hand
side'' and ``right-hand side'' when referring to equations or
inequalities. Since the paper is devoted to demonstrating Theorems
\ref{thm5.7.1} and \ref{thm1.24.1}, we will implicitly assume that the
network under consideration has a WMMF policy, except when stated
otherwise, and that the moment conditions on the interarrival and
service times given in Theorem \ref{thm5.7.1} hold. We assume the
network is subcritical only when explicitly stated.

\section{Markov process background}\label{sec2}

In this section, we provide a more detailed description of the
construction of the Markov process $X(\cdot)$ that underlies a WMMF
network. We then show how Theorem \ref{thm1.24.1} and its corollary
follow from Theorem~\ref{thm5.7.1}. Analogs of this material for
queueing networks are given in Bramson \cite{r1}. Because of the
similarity of the two settings, we present a summary here and refer the
reader to \cite{r1} for additional detail.

\subsection*{Construction of the Markov process}

As in (\ref{eq1.18.1}), we define the state space $S$ to be the set of
pairs $x=(z(\cdot),u)$, where $z(\cdot) = (z_{r}(\cdot))$ and
$z_{r}(\cdot )$ is a counting measure that maps $B
\subseteq\mathbb{R}^{+}$ to $\mathbb{Z}^{+,0}$, and $u=(u_{r})$, $r \in
\mathcal{R}$, has positive components. Here, $z(\cdot)$ corresponds to
the residual service times of documents and $u$ to the residual
interarrival times. (One could, as in (4.1) of \cite{r1}, distinguish
documents based on their ``age,'' which is not needed here.)

For the purpose of constructing a metric $d(\cdot, \cdot)$ on $S$, we
assign to each document the pair $(r_{i}, s_{i})$, $i=1,2,\ldots,$
where $r_{i} \in \{1, \ldots, |\mathcal{R}|\}$ denotes its route and
$s_{i}>0$ its residual service time. Documents are ordered so that
$s_{1} \le s_{2} \le\cdots ,$ with the decision for ties being made
based on a given ordering of the routes. When $i$ exceeds the number of
documents belonging to $x$, we assign the value $(r_{i}, s_{i}) =
(0,0)$. For $x, x' \in S$, with the coordinates labeled
correspondingly, we set
%
%
\begin{equation}
\label{eq2.3.1} d(x,x') = \sum^{\infty}_{i=1}
\bigl((|r_{i}-r'_{i}|+|s_{i}-s'_{i}|) \wedge 1 \bigr) +
\sum_{r}|u_{r}-u'_{r}|.
\end{equation}
One can check that $d(\cdot, \cdot)$ is separable and locally compact.
(See page 82 of \cite{r1} for details.) We equip $S$ with the standard
Borel $\sigma$-algebra inherited from $d(\cdot, \cdot)$, which we
denote by $\mathscr{S}$. In Proposition \ref{prop3.13.1}, we will show
$|\cdot|_{L}$, $|\cdot|_{R}$ and $|\cdot|_{A}$ are continuous in
$d(\cdot, \cdot)$.

The Markov process $X(t) = (Z(t,\cdot), U(t))$ underlying the network,
with $Z(t, \cdot)$ and $U(t)$ taking values $z(\cdot)$ and $u$ as
above, is defined to be the right continuous process whose evolution is
determined by the assigned WMMF policy. Documents are allocated service
according to the rates $\lambda_{r}(\cdot)$, which are constant in
between arrivals and departures of documents on routes. Upon an arrival
or departure, rates are re-assigned according to the policy. We note
that this procedure is not policy specific, and also applies to
$\alpha$-fair policies. By modifying the state space descriptor to
contain more information, one could also include more general networks.

Along the lines of page 85 of \cite{r1}, a filtration $(\mathcal
{F}_{t})$, $t\in [0,\infty]$, can be assigned to $X(\cdot)$ so that
$X(\cdot)$ is Borel right and, in particular, is strong Markov. The
processes $X(\cdot)$ fall into the class of piecewise-deterministic
Markov processes, for which the reader is referred to Davis \cite{r4}
for more detail.

\subsection*{Recurrence}

The Markov process $X(\cdot)$ is said to be \textit{Harris recurrent}
if, for some nontrivial $\sigma$-finite measure $\varphi$,
\[
\varphi(B) > 0 \mbox{ implies } P_{x}(\eta_{B} = \infty) = 1
\qquad\mbox{for all } x \in S,
\]
where $\eta_{B} = \int^{\infty}_{0}1\{{X(t)\in B}\}\,dt$. If $X(\cdot )$
is Harris recurrent, it possesses a stationary measure $\pi$ that is
unique up to a constant multiple. When $\pi$ is finite, $X(\cdot)$ is
said to be \textit{positive Harris recurrent}.

A practical condition for determining positive Harris recurrence can be
given by using petite sets. A nonempty set $A \in\mathscr{S}$ is said
to be \textit{petite} if for some fixed probability measure $a$ on $(0,
\infty)$ and some nontrivial measure $\nu$ on $(S, \mathscr{S})$,
\[
\nu(B) \le\int^{\infty}_{0}P^{t}(x,B)a(dt)
\]
for all $x \in A$ and $B \in\mathscr{S}$. Here, $P^{t}(\cdot, \cdot )$,
$t\ge 0$, is the semigroup associated with~$X(\cdot)$. As mentioned in
the \hyperref[intro]{Introduction}, a petite set $A$ has the property
that each set $B$ is ``equally accessible'' from all points $x \in A$
with respect to the measure~$\nu $. Note that any nonempty measurable
subset of a petite set is also petite.

For given $\delta> 0$, set
\[
\tau_{B}(\delta) = \inf\{t \ge\delta\dvtx X(t) \in B\}
\]
and $\tau_{B} = \tau_{B}(0)$. Then, $\tau_{B}(\delta)$ is a stopping
time. Employing petite sets and $\tau_{B}(\delta)$, one has the
following characterization of Harris recurrence and positive Harris
recurrence. (The Markov process and state space need to satisfy minimal
regularity conditions, as on page 86 of \cite{r1}.) The criteria are
from Meyn and Tweedie \cite{r10}; discrete time analogs of the
different parts of the proposition have long been known; see, for
instance, Nummelin \cite{r11} and Orey \cite{r12}.
\begin{theorem}
\label{thm2.8.1}
\textup{(a)} A Markov process $X(\cdot)$ is Harris recurrent if and only if
there exists
a closed petite set $A$ with
%
%
\begin{equation}
\label{eq2.8.2}
P_{x}(\tau_{A} < \infty) = 1 \qquad\mbox{for all } x \in S.
\end{equation}

\textup{(b)} Suppose the Markov process $X(\cdot)$ is Harris recurrent.
Then, $X(\cdot)$ is positive Harris recurrent if and only if there
exists a closed petite set $A$ such that for some $\delta> 0$,
%
%
\begin{equation}
\label{eq2.8.3}
\sup_{x \in A}E_{x}[\tau_{A}(\delta)] < \infty.
\end{equation}
\end{theorem}

One can apply Theorem \ref{thm2.8.1}, together with a stopping time
argument, to show the following version of Foster's criterion. It is
contained in Proposition 4.5 in \cite{r1}.
\begin{prop}
\label{prop2.9.1}
Suppose that $X(\cdot)$ is a Markov process, with
norm $\|\cdot\|$, such that for some $\varepsilon> 0$, $L > 0$ and
$M>0$,
%
%
\begin{equation}
\label{eq2.9.2}
E_{x}[\|X(M)\|] \le(\|x\| \vee L) - \varepsilon
\qquad\mbox{for all } x.
\end{equation}
Then, for $0 < \delta\le M$,
%
%
\begin{equation}
\label{eq2.9.3}
E_{x}[\tau_{A_{L}}(\delta)]
\le\frac{M}{\varepsilon}(\|x\| \vee L) \qquad\mbox{for all } x,
\end{equation}
where $A_{L}=\{x\dvtx\|x\| \le L\}$. In particular, if $A_{L}$ is
closed petite, then $X(\cdot)$ is positive Harris recurrent.
\end{prop}

\subsection*{Theorem \protect\ref{thm1.24.1} and its corollary}

Proposition \ref{prop2.9.1} and Theorem \ref{thm5.7.1} provide the main
tools for demonstrating Theorem \ref{thm1.24.1}. We also require
Proposition \ref{prop3.13.1}, which states that the norm $\|\cdot\|$ in
(\ref {eq1.20.1}) is continuous in the metric $d(\cdot,\cdot)$, and
hence that $A_{L}=\{ x\dvtx \|x\| \le L\}$ is closed for each $L$.
Together, they give a quick proof of the theorem.
\begin{pf*}{Proof of Theorem \protect\ref{thm1.24.1}}
From the conclusion (\ref{eq5.3.10}) in Theorem \ref{thm5.7.1}, we know
that the assumption (\ref{eq2.9.2}) in Proposition \ref{prop2.9.1} is
satisfied for some $L$, with $M=N^{3}$ and $\varepsilon=
\varepsilon_{1}N^{2}$. In Theorem \ref{thm1.24.1}, it is assumed that
$A_{L}$ is petite; by Proposition \ref{prop3.13.1}, we also know it is
closed. So, all of the assumptions in Proposition \ref{prop2.9.1} are
satisfied, and hence $X(\cdot)$ is positive Harris recurrent.
\end{pf*}

Corollary \ref{col1.26.1} follows immediately from Theorem \ref
{thm1.24.1} and the assertion, before the statement of the corollary,
that the sets $A_{L}$ are petite under the assumptions (\ref{eq1.25.1})
and (\ref{eq1.25.2}). A somewhat stronger version of the analogous
assertion for queueing networks is demonstrated in Proposition 4.7 of
\cite{r1}. (The proposition states that the sets $A$ are uniformly
small.) The reasoning is the same in both cases and does not involve
the policy of the network. The argument, in essence, requires that one
wait long enough for the network to have at least a given positive
probability of being empty; this time $t$ does not depend on $x$ for
$\| x\|\le L$. One uses (\ref{eq1.25.1}) for this. By using
(\ref{eq1.25.2}), one can also show that the joint distribution
function of the residual interarrival times has an absolutely
continuous component at this time, whose density is bounded away from
0. It will follow that the set $A_{L}$ is petite with respect to $\nu$,
with $a$ chosen as the point mass at $t$, if $\nu$ is concentrated on
the empty states, where it is a small enough multiple of $|\mathcal
{R}|$-dimensional Lebesque measure restricted to a small cube.

\section[Summary of the proof of Theorem 1.1]{Summary of the proof of Theorem \protect\ref{thm5.7.1}}\label{sec3}

As mentioned in Section \ref{intro}, the norm $\|\cdot\|$ in Theorem
\ref{thm5.7.1} consists of three components, with
%
%
\begin{equation}\label{eq5.3.6}
\|x\| = |x|_{L} + |x|_{R} + |x|_{A}
\end{equation}
for each $x$. After introducing these components, we will state the
bounds associated with each of them that we will need, leaving their
proofs to the remaining sections. We then show how Theorem
\ref{thm5.7.1} follows from these bounds.

\subsection*{Definition of norms}

We first define $|x|_{L}$. This requires a fair amount of notation,
which we will introduce shortly. We begin by expressing $|x|_{L}$ in
terms of this notation; when the notation is then specified, we
motivate it by referring back to~$|x|_{L}$.

We set $|x|_{L} = \sup_{r,s} |x|_{r,s}$ for $r \in\mathcal{R}$ and $s>0$,
where
%
%
\begin{equation}
\label{eq5.3.1} |x|_{r,s} =
\frac{w_{r}(1+as_{N})z^{*}_{r}(s)}{\nu_{r}\Gamma
(\bar{H}^{*}_{r}(s_{N}))} .
\end{equation}
We need to define the terms ${H}^{*}_{r}(\cdot)$, $z^{*}_{r}(\cdot)$,
$\Gamma(\cdot)$, $a$ and $s_{N}$.

Starting with ${H}^{*}_{r}(\cdot)$ and $z^{*}_{r}(\cdot)$, we recall
the distribution functions $H_{r}(\cdot)$ and counting measure
$z_{r}(\cdot)$ from Section \ref{intro}. In (\ref{eq5.3.1}), we will
require their analogs ${H}^{*}_{r}(\cdot) $ and $z^{*}_{r}(\cdot)$ to
have densities with bounded first derivatives and to be ``close'' to
$H_{r}(\cdot)$ and $z_{r}(\cdot)$. For this, we define $H^{*}_{r}
(\cdot)$ and $z^{*}_{r}(\cdot)$ as the convolutions of $H_{r}(\cdot )$
and of $z_{r}(\cdot)$ by an appropriate distribution function
$\Phi(\cdot)$ with density $\phi(\cdot)$. Setting
%
%
\begin{equation}
\label{eq5.1.4}
\phi(s) = \cases{ \frac{2}{3}ebe^{-bs}, &\quad for $s >
1/b$, \vspace*{2pt}\cr
\frac{2}{3}b^{2}s, &\quad for $s \in(0, 1/b]$, \vspace*{2pt}\cr
0, &\quad for $s \le0$,}
\end{equation}
%
for $b \in\mathbb{Z}^{+}$ with $b \ge2$, $\phi(\cdot)$ is the density
of $\Phi(s) = \int^{s}_{-\infty}\phi(s')\,ds'$. We note that
$\Phi(\cdot)$ has mean at most $2/b$ and that $\phi(\cdot)$ satisfies
%
%
\begin{equation}
\label{eq5.1.10}
\phi'(s) \le b^{2} \quad\mbox{and}\quad \phi(s+s')/\phi(s) \ge e^{-bs'}
\end{equation}
for $s,s'>0$. The above properties and the exponential tail of $\phi
(\cdot)$ will be useful later when analyzing $|\cdot|_{L}$ and
$|\cdot|_{R}$ [as in (\ref{eq5.2.1}), (\ref{eq10.5.3}),
(\ref{eq6.88.1}), (\ref{eq6.88.2}) and (\ref{eq46.5.4})].

Convoluting by $\Phi(\cdot)$, we set
%
%
\begin{eqnarray}
\label{eq5.1.6}
H^{*}_{r}(s) &=& (H_{r}*\Phi)(s) = \int^{\infty}_{0} \Phi
(s-s')\,dH(s'),\nonumber\\[-8pt]\\[-8pt]
z^{*}_{r}((0,s]) &=& (z_{r}*\Phi)((0,s]) = \int^{\infty}_{0}
\Phi(s-s')\,
dz_{r}((0,s'])\nonumber
\end{eqnarray}
with $z^{*}_{r}(B)$ being defined analogously for $B
\subseteq\mathbb{R}^{+}$. Differentiating both quantities in
(\ref{eq5.1.6}), we also set
%
%
\begin{eqnarray}
\label{eq5.1.8}
h^{*}_{r}(s) &=& (H_{r}*\Phi)'(s) = \int^{\infty}_{0} \phi(s-s')
\,dH(s'), \nonumber\\[-8pt]\\[-8pt]
%
z^{*}_{r}(s) &=& (z_{r}*\Phi)'((0,s]) = \int^{\infty}_{0} \phi(s-s')
\,dz_{r}((0,s']).\nonumber
\end{eqnarray}
Convolution by $\Phi(\cdot)$, as in (\ref{eq5.1.6}), produces a measure
$z_{r}^{*}(\cdot)$ that approximates $z_{r}(\cdot)$ and possesses a
density.

Since $H_{r}(\cdot)$ is assumed to have a finite $(2+\delta_{1})$th
moment for all $r$, the same is true for $H^{*}_{r}(\cdot)$. This
implies that for appropriate $C_{1}\ge1$,
%
%
\begin{equation}
\label{eq5.2.1}
\bar{H}^{*}_{r}(s) \le\frac{C_{1}}{(1+s)^{2+\delta_{1}}} \qquad\mbox{for all }
s > 0 \mbox{ and } r \in\mathcal{R},
\end{equation}
for $\delta_{1}$ chosen as in Theorem \ref{thm5.7.1}. We assume wlog
that $\delta_{1} \le1$. Since the difference of the means of
$H^{*}_{r}(\cdot)$ and $H_{r}(\cdot)$ is at most $2/b$ for each $r$ and
$H(\cdot)$ is subcritical, $H^{*}(\cdot)$ will also be subcritical for
large enough $b$.

We set
%
%
\begin{equation}
\label{eq5.1.3}
\Gamma(\sigma) = \sigma+ C_{2}a \sigma^{\gamma} \qquad\mbox{for }
\sigma\in[0,1].
\end{equation}
We choose $\gamma\in(0, \delta_{1}/24]$, $C_{2}\ge2C_{1}/\gamma$ and
$a$ small enough so that $a \le(1/C_{2}) \wedge1$ and (\ref{eq5.2.5})
is satisfied. One can think of $\Gamma(\cdot)$ as being almost linear
for values of $\sigma$ that are not too small; the power $\gamma$ needs
to be small in order to be able to bound $|x|_{r,s}$ later on for
$H^{*}_{r}(s_{N})$ small; $\gamma>0$ is needed so that the integral in
(\ref{eq5.2.4}) is finite.

We set $s_{N} = s \wedge(N_{H_{r}}+1)$ for $N \in\mathbb{Z}^{+}$, where
%
%
\begin{equation}
\label{eq90.1.1}
N_{H_{r}} = (\bar{H}^{*}_{r})^{-1}(1/N^{4})\wedge N.
\end{equation}
It follows that
%
%
\begin{equation}
\label{eq90.3.1}
1/N^{4} \le\bar{H}^{*}_{r}(N_{H_{r}}) \le C_{1}/ N^{2+\delta_{1}}.
\end{equation}
If $H^{*}_{r}(\cdot)$ has a relatively fat tail, say $\bar
{H}^{*}_{r}(s) \sim s^{-3}$, (\ref{eq90.1.1}) implies that $N_{H_{r}} =
N$; otherwise, $N_{H_{r}} < N$ and $\bar{H}^{*}_{r}(N_{H_{r}}) =
1/N^{4}$. In either case, it will follow from (\ref{eq90.3.1}) that
$\Gamma(\bar{H}^{*}_{r}(N_{H_{r}}))$ is ``large enough'' for us to
adequately bound $|x|_{r,s}$. We will assume that $N\in \mathbb{Z}^{+}$
is chosen large enough so $N \ge1/a$ and $N_{H_{r}} \ge1$ for all~$r$.

The norm $|\cdot|_{L}$ has been defined with the following motivation.
As the process $X(\cdot)$ evolves, documents arrive at each route, are
served, and eventually depart. In Proposition \ref{prop46.3.2}, we will
show that, under certain assumptions for $X(t)$ on $t \in[0, N^{3}]$,
for large enough $b$,
%
%
\begin{equation}
\label{eq3.6.1}
\lambda^{w}(t) \ge(1 + \varepsilon_{2})/|x|_{L} \qquad\mbox{on } t\in[0, N^{3}],
\end{equation}
for some $\varepsilon_{2}> 0$, because of the subcriticality of
$H^{*}(\cdot)$.
%
%
Reasoning as below (\ref{eq1.30.1}), this will imply that individual
documents receive enough service so that $|X(t)|_{L}$ decreases on
average over $[0, N^{3}]$. More specifically, the increase in the term
$\Gamma(\bar{H}^{*}_{r} (s_{N}))$ in (\ref{eq5.3.1}), after translating
$s_{N}$ according to the service of documents, will compensate for the
arrival of new documents. For documents with residual service $s \le
N_{H_{r}} \le N+1$ at $t=0$, the term $1+as_{N}$ in (\ref{eq5.3.1}),
after translating $s_{N}$ according to the service of documents, will
decrease sufficiently over $[0, N^{3}]$ to produce the term
$-\varepsilon_{1}N^{2}$ in (\ref{eq5.3.10}). For documents with
residual service $s > N_{H_{r}}$, we will instead need to employ the
norm \mbox{$|\cdot |_{R}$}, which we introduce next. (On $(N_{H_{r}},
N_{H_{r}}+1]$, the intervals overlap.)

The norm $|\cdot|_{R}$ in (\ref{eq5.3.6}) is given by
%
%
\begin{equation}
\label{eq5.3.2}
|x|_{R} = M_{1}\sum_{r}\kappa_{N,r}\int^{\infty}_{N_{H_{r}}}
N_{r}(s)z^{*}_{r}(s)\,ds.
\end{equation}
We need to identify the terms $\kappa_{N,r}$, $N_{r}(\cdot)$ and
$M_{1}$. We set
%
%
\begin{equation}
\label{eq3.8.1}
\kappa_{N,r} = 1/\Gamma(\bar{H}^{*}_{r}(N_{H_{r}}))
\end{equation}
and
%
%
\begin{equation}
\label{eq90.1.2}
N_{r}(s) = \cases{
s^{2}/N, &\quad for $s > N$, \cr
s, &\quad for $s \le N$.}
\end{equation}
Later on, we will also employ $\kappa_{N} \stackrel{\mathrm{def}}{=}
\max_{r}\kappa_{N,r}$. For the term $M_{1}$, we will require that
%
%
\begin{equation}
\label{eq3.8.2}
M_{1} \ge8 C_{3}\Bigl(\max_{r,r'} w_{r}/w_{r'}\Bigr),
\end{equation}
where $C_{3}$ is chosen as in (\ref{eq5.4.2}).

Since \mbox{$|\cdot|_{R}$} is given by a weighted sum of the residual service
times of the different documents, it will be easier to work with than
$|\cdot|_{L}$, which is a supremum. For smaller values of $s$, we
required $|\cdot|_{L}$ because\vspace*{1pt} of the nature of the WMMF policy.
Because of the bound on $\bar{H}^{*}_{r}(\cdot)$ in (\ref{eq5.2.1}),
the impact of large residual service times on the evolution of
$X(\cdot)$ will typically be small, and so one can employ the ``more
generous'' definition over $(N_{H_{r}}, \infty)$ given in
(\ref{eq5.3.2}).

As we will see in Section \ref{sec5}, we will require the presence of
the term $N_{r}(s)$ in the integrand in (\ref{eq5.3.2}) to ensure that
the integral decreases sufficiently rapidly from the service of
documents when the integral is large. This will rely on $N'_{r}(s)
\ge1$ on $(N_{H_{r}}, \infty )$. For $s > N$, the denominator $N$ in
$s^{2}/N$ is needed so that the expected increase due to arrivals does
not dominate the term $-\varepsilon_{1}N^{2}$ in (\ref {eq5.3.10}),
which was mentioned in the motivation for the definition of
\mbox{$|\cdot|_{L}$}. This denominator is not needed for $s \in[N_{H_{r}},N)$
because (\ref {eq90.1.1}) will guarantee that the integrand is already
sufficiently small there. The terms $\kappa_{N,r}$ are needed when we
combine the norms $|\cdot|_{L}$ and $|\cdot|_{R}$ in $\|\cdot\|$,
because of the denominator $\Gamma (\cdot)$ in $|\cdot|_{r,s}$.

The norm $|\cdot|_{A}$ in (\ref{eq5.3.6}) is needed for the residual
interarrival times. It is given by
%
%
\begin{equation}
\label{eq5.3.5}
|x|_{A} = \frac{1}{N}\max_{r}\theta(u_{r}),
\end{equation}
where $\theta(y)$, $y>0$, satisfies the following properties. We assume
that $\theta(y)>0$ for all $y$ and that $\theta(\cdot)$ and $\theta
'(\cdot)$ are strictly increasing, with
%
%
\begin{equation}
\label{eq5.8.1}
\theta'(y) \to\infty\qquad\mbox{as } y \to\infty.
\end{equation}
We also assume that
%
%
\begin{equation}
\label{eq5.8.3}
\theta(y) \le y^{2} \qquad\mbox{for all } y,
\end{equation}
and that $\theta(\cdot)$ grows sufficiently slowly so that
%
%
\begin{equation}
\label{eq5.8.2}
E[\theta(\xi_{r})] < \infty\qquad\mbox{for all } r.
\end{equation}
Since $E[\xi_{r}]<\infty$, it is possible to specify such $\theta
(\cdot)$
that also satisfy the previous two displays.

The above properties for $\theta(\cdot)$ will enable us to show that
the expected value of $|X(t)|_{A}$ will decrease over time when
$|X(t)|_{A}$ is large. In particular, because of (\ref{eq5.8.1}) and
(\ref{eq5.8.2}), the decrease in $|\cdot|_{A}$ due to decreasing
residual interarrival times will, on the average, dominate the increase
in $|\cdot|_{A}$ due to new interarrival times that occur when a
document joins a route. The argument for this is given in Section
\ref{sec4} and is fairly quick. We note that when $\xi_{r}$ are all
exponentially distributed, the term $|\cdot|_{A}$ may be omitted in the
definition of \mbox{$\|\cdot\|$}.

The reader attempting to understand the norm \mbox{$\|\cdot\|$} should first
concentrate on $|\cdot|_{L}$, which was chosen to accommodate the WMMF
policy. When the service distributions $H_{r}(\cdot)$ all have compact
support and the interarrival times are exponentially distributed, one
may, in fact, set $\|x\|=|x|_{L}$ for a large enough choice of~$N$.

We note that the norm \mbox{$|\cdot|_{L}$} is not appropriate for weighted
$\alpha$-fair policies. In particular, the supremum and the function
$\Gamma(\cdot)$ in its definition are not appropriate factors in this
context. On the other hand, $|\cdot|_{R}$, with suitable $M_1$, and
$|\cdot|_{A}$ should still be applicable to $\alpha$-fair policies,
provided a suitable replacement of \mbox{$|\cdot|_{L}$} can be found.

In order to apply Proposition \ref{prop2.9.1} in the proof of Theorem
\ref{thm1.24.1} in Section \ref{sec2}, we needed to know that the sets
$A_{L} = \{x\dvtx\|x\|\le L\}$ are closed. For this, it suffices to
show the norm \mbox{$\|\cdot\|$} is continuous in the metric $d(\cdot,\cdot)$
that is given in (\ref{eq2.3.1}).
\begin{prop}
\label{prop3.13.1}
The norm $\|\cdot\|$ in (\ref{eq5.3.6}) is continuous in the metric
$d(\cdot,\cdot)$ given by (\ref{eq2.3.1}).
\end{prop}
\begin{pf}
It suffices to show $|\cdot|_{L}$, $|\cdot|_{R}$ and $|\cdot|_{A}$ are
each continuous in $d(\cdot,\cdot)$. For $|\cdot|_{L}$, note that the
coefficients of $z^{*}_{r}(s)$ in (\ref{eq5.3.1}) are bounded. On the
other hand, if $d(x,x') \le\varepsilon< 1$, then one can show, by using
the first part of (\ref{eq5.1.10}), that
%
%
\begin{equation}
\label{eq13.3.2}
|z^{*}_{r}(s) - z{}^{\prime,*}_{r}(s)| \le b^{2} \varepsilon\qquad
\mbox{for all } s
\mbox{ and } r,
\end{equation}
where $z^{\prime,*}_{r} = (z')^{*}_{r}(s)$. It follows from this and
(\ref{eq5.3.1}) that $|\cdot|_{L}$ is in fact Lipschitz in
$d(\cdot,\cdot)$.

For $|\cdot|_{R}$, one can apply both parts of (\ref{eq5.1.10}) to show
with a bit of work that, if $d(x,x') \le\varepsilon< 1$ and $x$ has no
residual service times greater than $M$, for given $M$, then
%
%
\begin{equation}
\label{eq13.3.3}
\int^{\infty}_{N_{H_{r}}}N_{r}(s) |z^{*}_{r}(s) - z^{\prime,*}_{r}(s)|\,ds
\le(M+1)^{2}b^{2}\varepsilon+ (1-e^{-b\varepsilon})|x|_{R}
\end{equation}
for all $r$. Since the coefficients of
$\int^{\infty}_{N_{H_{r}}}N_{r}(s) z^{*}_{r}(s)\,ds$ in $|x|_{R}$ are
bounded and the RHS of (\ref{eq13.3.3}) goes to 0 as $\varepsilon\to0$,
the continuity of $|\cdot|_{R}$ follows.

Since $\theta'(u_{r})$ is bounded for bounded values of $u_{r}$,
$|\cdot|_{A}$ is also continuous.
\end{pf}

In addition to the norms in (\ref{eq5.3.6}), we will employ the
following norms in showing Theorem \ref{thm5.7.1}:
%
%
%
\begin{equation}
\label{eq5.3.7}
|x| = \sum_{r}z_{r}(\mathbb{R}^{+}) = \sum_{r}z^{*}_{r}(\mathbb{R}^{+})
\end{equation}
and
%
%
\begin{equation}
\label{eq5.3.3}
|x|_{K} = \sum_{r}\kappa_{N,r}z^{*}_{r}((N_{H_{r}},\infty)).
\end{equation}
Although we will not employ them in this section, we also introduce the norms
%
%
\begin{equation}
\label{eq5.3.11}
|x|_{1} = \sum_{r}z^{*}_{r}((0,N_{H_{r}}]),\qquad %
|x|_{2} = \sum_{r}z^{*}_{r}((N_{H_{r}}, \infty))
\end{equation}
and
%
%
\begin{equation}
\label{eq5.3.12}
|x|_{S} = |x|_{L} +
\max_{r}\frac{w_{r}}{\rho_{r}}z^{*}_{r}((N_{H_{r}},\infty)).
\end{equation}
It obviously follows from (\ref{eq5.3.7}) and (\ref{eq5.3.11}) that
$|x| = |x|_{1} + |x|_{2}$. The norm $|\cdot|_{S}$ will be employed in
Proposition \ref{prop46.3.2} to derive the bound given in
(\ref{eq3.6.1}).

\subsection*{Bounds on $|\cdot|_{L}, |\cdot|_{R}$ and $|\cdot|_{A}$}

In order to derive (\ref{eq5.3.10}), we need bounds on $|\cdot|_{L}$,
$|\cdot|_{R}$ and $|\cdot|_{A}$ as the process $X(t)$ evolves from
$t=0$ to
$t=N^{3}$. We first need to specify the term $L$ 
appearing in (\ref{eq5.3.10}). We choose $l_{1}$ large enough so that
%
%
\begin{equation}
\label{eq5.8.6}
\theta'(l_{1}/2) \ge M_{1}N
\end{equation}
and, for all $r$,
%
%
\begin{equation}
\label{eq5.8.7}
E[\theta(\xi_{r}); \xi_{r}>l_{1}/2] \le(1 /|\mathcal{R}|)
P(\xi_{r}>N^{3}).
\end{equation}
%
We set
%
%
\begin{equation}
\label{eq5.8.5}
L_{1} = \frac{1}{N}\theta(l_{1})
\end{equation}
and
%
%
\begin{equation}
\label{eq5.8.8}
L = 6(\kappa^{2}_{N}N^{17} \vee L_{1}).
\end{equation}

For $|\cdot|_{L}$, we employ the bound from Proposition \ref
{prop50.1.1} that, for large enough $N$ and $b$, small enough $a$, and
appropriate $C_{3}$ and $\varepsilon_{3}> 0$,
%
%
\begin{eqnarray}\qquad
\label{eq5.4.2}
&&E_{x}[|X(N^{3})|_{L}] - |x|_{L} \nonumber\\[-8pt]\\[-8pt]
&&\qquad\le C_{3}N^{3} \cdot1\{|x| \le N^{6}\}
+ [C_{3}(|x|_{K}/|x|)N^{3} - \varepsilon_{3}N^{2}]
\cdot1\{|x|>N^{6}\}\nonumber
\end{eqnarray}
for all $x$. The precise value of $\varepsilon_{3}$ is not important;
in Proposition \ref{prop50.1.1}, it is given by $\frac{1}{4}
\min_{r}w_{r}$. We assume wlog that $\varepsilon_{3}\le C_{3}$.

For $|\cdot|_{R}$, we employ the bound from Proposition \ref
{prop10.7.1} that, for given $\varepsilon_{4}> 0$, large enough $N$,
and $M_{2}=\frac{1}{8}
(1\wedge\min_{l}c_{l})(\min_{r,r'}(w_{r}/w_{r'}))M_{1} \ge C_{3}$,
%
%
\begin{eqnarray}
\label{eq5.4.1}
&&
E_{x}[|X(N^{3})|_{R}] - |x|_{R} \nonumber\\
&&\qquad\le\varepsilon_{4}N^{2} - M_{2}(|x|_{K}/|x|) N^{3} \cdot1\{|x|>N^{6}\}
\\
&&\qquad\quad{} - \kappa_{N}N^{4} \cdot1 \{|x|_{R}>\kappa^{2}_{N}N^{17}, |x| \le
N^{6}\}\nonumber
\end{eqnarray}
for all $x$. We will later choose $\varepsilon_{4}$ small with respect
to $\varepsilon_{3}$;
the constant $C_{3}$ is chosen as in (\ref{eq3.8.2}) and (\ref{eq5.4.2}).

For $|\cdot|_{A}$, we will show in Proposition \ref{prop6.1.4} and Proposition
\ref{prop6.2.5} that, for this choice of $\varepsilon_{4}$ and large
enough $N$,
%
%
\begin{equation}
\label{eq5.5.1}
E_{x}[|X(N^{3})|_{A}] - |x|_{A} \le\varepsilon_{4}N^{2} - M_{1}
N^{3}\cdot1\{|x|_{A}>L/6\}
\end{equation}
for all $x$.

\subsection*{Derivation of (\protect\ref{eq5.3.10}) from (\protect
\ref{eq5.4.2}), (\protect\ref{eq5.4.1})
and (\protect\ref{eq5.5.1})}

We now derive (\ref{eq5.3.10}) from these three bounds. Adding the RHS
of (\ref{eq5.4.2}), (\ref{eq5.4.1}) and (\ref{eq5.5.1}), one obtains,
for large enough $N$ and $b$, and small enough $a$,
%
%
\begin{equation}
\label{eq5.5.2}
E_{x}[\|X(N^{3})\|] - \|x\| \le2 C_{3}N^{3}
\end{equation}
for all $x$. 
We next consider the behavior of the LHS of (\ref{eq5.5.2}) for $\|x\|
> L/2$, where $L$ is given by (\ref{eq5.8.8}). This condition implies
that either $|x|_{L}>\kappa^{2}_{N} N^{17}$,
$|x|_{R}>\kappa^{2}_{N}N^{17}$ or $|x|_{A}\ge L/6$.

Suppose first that $|x|_{L}>\kappa^{2}_{N}N^{17}$. We note that if $|x|
\le N^{6}$, then
\[
|x|_{L} \le C_{4}N^{8}
\]
for some constant $C_{4}$. This bound follows from the definition of
$|x|_{L}$ in (\ref{eq5.3.1}), together with the bounds $z^{*}_{r}(s)\le
12b^{2}|x|$ for all $s$, $s_{N} \le N$, and
$\Gamma(\bar{H}^{*}_{r}(s_{N})) \ge C_{5} /N$, for some $C_{5}> 0$
[which follows from (\ref{eq90.3.1}) and $\gamma\le1/4$]. Therefore, if
$|x|_{L} > \kappa^{2}_{N}N^{17}$ and $N$ is large enough so that
$\kappa_{N} \ge1$, one must have $|x|>N^{6}$.

On the other hand, it follows from (\ref{eq5.4.1}) that, on $|x|>N^{6}$,
%
%
\begin{equation}
\label{eq5.5.3}
E_{x}[|X(N^{3})|_{R}] - |x|_{R} \le\varepsilon_{4}N^{2} -
M_{2}(|x|_{K}/|x|) N^{3}.
\end{equation}
Adding the terms corresponding to $|x|>N^{6}$ in (\ref{eq5.4.2}) and
(\ref{eq5.5.1}) to this implies that, 
for $|x|>N^{6}$, and hence for $|x|_{L} > \kappa^{2}_{N} N^{17}$,
%
%
%
\begin{eqnarray}\quad
\label{eq5.5.4}
E_{x}[\|X(N^{3})\|] - \|x\| &\le& (2\varepsilon_{4}- \varepsilon_{3})N^{2} +
(C_{3}- M_{2})
(|x|_{K}/|x|)N^{3}\nonumber\\[-8pt]\\[-8pt]
&\le& -\varepsilon_{1}N^{2},\nonumber
\end{eqnarray}
where the latter inequality follows for $\varepsilon_{4}\le\varepsilon
_{3}/3$ and $\varepsilon_{1} \stackrel{\mathrm{def}}{=}
\varepsilon_{3}/3$, since $M_{2} \ge C_{3}$.

Suppose next that $|x|_{R} > \kappa^{2}_{N}N^{17}$ and $|x| \le N^{6}$.
Adding up the corresponding terms from (\ref{eq5.4.2}), (\ref{eq5.4.1})
and (\ref{eq5.5.1}) implies that the LHS of (\ref{eq5.5.4}) is at most
%
%
\begin{equation}
\label{eq5.11.1}
2\varepsilon_{4}N^{2} + C_{3}N^{3} - \kappa_{N} N^{4} \le- \varepsilon_{1}N^{3}
\end{equation}
for large $N$, which is better than the bound in (\ref{eq5.5.4}).

Suppose finally that $|x|_{A} \ge L/6$. We need to consider only the
case $|x| \le N^{6}$, since $|x| > N^{6}$ is covered by
(\ref{eq5.5.4}). In this case, it follows from (\ref{eq5.4.2}),
(\ref{eq5.4.1}) and (\ref{eq5.5.1}) that the LHS of (\ref{eq5.5.4}) is
at most
%
%
\begin{equation}
\label{eq5.5.5}
(3C_{3}- M_{1})N^{3} \le-\varepsilon_{1}N^{3},
\end{equation}
since $M_{1} \ge4 C_{3}$.

Together, (\ref{eq5.5.4}), (\ref{eq5.11.1}) and (\ref{eq5.5.5}) imply
that, for large enough $N$ and $b$, and small enough $a$,
%
%
\begin{equation}
\label{eq3.22.1}
E_{x}[\|X(N^{3})\|] - \|x\| \le- \varepsilon_{1}N^{2}
\end{equation}
for all $\|x\| > L/2$. Since for large $N$,
\[
L - L/2 \ge2C_{3}N^{3} + \varepsilon_{1}N^{2},
\]
(\ref{eq5.3.10}) follows easily form (\ref{eq5.5.2}) and (\ref{eq3.22.1}).

\section{Upper bounds on $E_{x}[|X(N^{3})|_{A}]$}\label{sec4}

In this section, we will demonstrate the inequality (\ref{eq5.5.1}) for
the upper bounds on $E_{x}[|X(N^{3})|_{A}] - |x|_{A}$. In Proposition
\ref{prop6.1.4}, we obtain the first term on the RHS of (\ref
{eq5.5.1}); this holds for all $x$. We then obtain a better bound in
Proposition \ref{prop6.2.5}, which is valid on $|x| \ge L/6$. Both
parts require just standard techniques.

The first bound employs the following elementary inequality on the
residual interarrival times at time $N^{3}$:
%
%
\begin{equation}
\label{eq6.1.1}
|X(N^{3})|_{A} \le|x|_{A} \vee\frac{1}{N}\max
\{\theta(\xi_{r}(k))\dvtx
r \in\mathcal{R}, k \in[2, A_{r} (N^{3})+1] \}. 
\end{equation}
Here and in later sections, $A_{r}(t)$ denotes the cumulative number of
arrivals at the route $r$ by time $t$; $A(t)$ will denote the
corresponding vector. The inequality $k \le A_{r}(t)+1$ implies that
the interarrival epoch associated with $\xi_{r}(k)$ has already begun
by time $t$. Recall that $\xi_{r}(1)$ is the initial residual time at
route $r$ and $\xi_{r}(2), \xi_{r}(3), \ldots$ are i.i.d. random
variables, and $\theta(\cdot)$ satisfies (\ref
{eq5.3.5})--(\ref{eq5.8.2}).
\begin{prop}
\label{prop6.1.4}
For any $\varepsilon> 0$ and large enough $N$, not depending on $x$,
%
%
\begin{equation}
\label{eq6.1.3}
E_{x}[|X(N^{3})|_{A}] - |x|_{A} \le\varepsilon N^{2}.
\end{equation}
\end{prop}
\begin{pf}
By (\ref{eq5.8.2}), $E[\theta(\xi_{r})] < \infty$ for all $r$. One can
therefore show with some estimation that, for each $r$,
%
%
\begin{equation}
\label{eq6.1.2}
\frac{1}{t}E_{x} \Bigl[\max_{k\in[2,A_{r}(t)+1]} \theta
(\xi_{r}(k)) \Bigr] \to0,
\end{equation}
uniformly in $x$ as $t \to\infty$. For fixed $x$, (\ref{eq6.1.2})
follows immediately from (4.83) of~\cite{r1}; since $A_{r}(t)$
decreases when $\xi_{r}(1)$ increases, this limit is uniform in $x$.

Inequality (\ref{eq6.1.3}) follows immediately from (\ref{eq6.1.1}) and
(\ref{eq6.1.2}), with $t = N^{3}$.
\end{pf}

We proceed to Proposition \ref{prop6.2.5}. For the proposition, it will
be useful to decompose $|X(t)|_{A} - |x|_{A}$ as
%
%
\begin{equation}
\label{eq6.2.6}
|X(t)|_{A} - |x|_{A} = I_{A}(t) - D_{A}(t),
\end{equation}
where $I_{A}(t)$ and $D_{A}(t)$ are the nondecreasing functions
corresponding to the cumulative increase and decrease of
$|X(\cdot)|_{A}$ up to time $t$. That is, \mbox{$I_{A}(0) = D_{A}(0) = 0$},
with $I_{A}(t)$ being the jump process, with
\[
I_{A}(t) - I_{A}(t-) = |X(t)|_{A} - |X(t-)|_{A}
\]
and $D'_{A}(t)$ being the rate of decrease of $|X(t)|_{A}$ at other
times. We note that $D_{A}(t)$ is locally Lipschitz, with $D'_{A}(t)$
defined except at arrivals. In particular, since $U'_{r}(t)=-1$ except
at arrivals,
%
%
\begin{equation}
\label{eq6.2.7}
D'_{A}(t) = \frac{1}{N}\max_{r}\theta'(U_{r}(t))\qquad
\mbox{almost everywhere.}
\end{equation}
We recall the definitions for $l_{1}, L_{1}$ and $M_{1}$ in
(\ref{eq5.8.6})--(\ref{eq5.8.5}) and (\ref{eq3.8.2}).
\begin{prop}
\label{prop6.2.5}
Suppose that $|x|_{A} \ge L/6$. Then, for large enough $N$ not
depending on
$x$,
%
%
\begin{equation}
\label{eq6.2.8}
E_{x}[|X(N^{3})|_{A}] - |x|_{A} \le1 - M_{1}N^{3} \le- M_{1}N^{3}/2.
\end{equation}
\end{prop}
\begin{pf}
We first show that
%
%
\begin{equation}
\label{eq6.2.9}
D_{A}(N^{3}) \ge M_{1}N^{3}.
\end{equation}
Since $|x|_{A} \ge L/6 = \kappa^{2}_{N}N^{17} \vee L_{1}$ and $\theta
(y) \le y^{2}$ for all $y$, one has, for $N \ge2$, that $\max_{r}u_{r}
\ge N^{8} \vee l_{1}$. So, for all $t \in[0, N^{3}]$,
%
%
\begin{equation}
\label{eq6.2.10}
\max_{r}u_{r} - \max_{r}U_{r}(t) \le
N^{3} \le\frac{1}{2} \max_{r}u_{r}.
\end{equation}
Consequently, for all $t \in[0, N^{3}]$,
%
%
\begin{equation}
\label{eq6.2.11}
\max_{r}U_{r}(t) \ge\frac{1}{2}
\max_{r}u_{r} \ge N^{3} \vee\frac{1}{2} l_{1}.
\end{equation}
Moreover, by (\ref{eq5.8.6}) and (\ref{eq6.2.7}), for $\max_{r}U_{r}(t)
\ge \frac{1}{2} l_{1}$, $D'_{A}(t) \ge M_{1}$ almost everywhere.
Together with (\ref{eq6.2.11}), this implies $D'_{A}(t) \ge M_{1}$
almost everywhere on $[0,N^{3}]$, and hence (\ref{eq6.2.9}) holds.

On account of (\ref{eq6.2.9}), in order to show (\ref{eq6.2.8}), it suffices
to show
%
%
\begin{equation}
\label{eq6.2.12}
E_{x}[I_{A}(N^{3})] \le1 
\end{equation}
for large $N$. To obtain (\ref{eq6.2.12}), we first note that, for each
$r$, there cannot be more than one interarrival time occurring over
$(0,N^{3}]$ with value greater than $N^{3}$. Moreover, because of (\ref
{eq6.2.11}), only interarrival times with value at least $N^{3}
\vee(l_{1}/2)$ can contribute to $I_{A}(N^{3})$. The expectation of
$\theta(\xi_{r})$, for $\xi_{r}$ conditioned on being greater than
$N^{3}$ and restricted to being greater than $l_{1}/2$, is
%
%
\begin{equation}
\label{eq6.2.13}
E [\theta(\xi_{r}); \xi_{r} > l_{1}/2 ]/P(\xi_{r} > N^{3}).
\end{equation}
(If $\xi_{r}$ is bounded above by $N^{3}$, set the ratio equal to 0.)
It follows that, for any $x$,
%
%
\begin{equation}
\label{eq6.2.14}
E_{x} [I_{A}(N^{3})] \le\frac{1}{N} \sum_{r}
E[\theta(\xi_{r}); \xi_{r} > l_{1}/2]/P(\xi_{r} > N^{3}).
\end{equation}
By (\ref{eq5.8.7}), the RHS of (\ref{eq6.2.14}) is at most $1/N$,
which implies
(\ref{eq6.2.12}).
\end{pf}

\section{Upper bounds on $E_{x}[|X(N^{3})|_{R}]$}\label{sec5}

In this section, we will demonstrate the following proposition for the
upper bounds on $E_{x}[|X(N^{3})|_{R}] - |x|_{R}$, where $|\cdot|_{R}$
is the norm introduced in (\ref{eq5.3.2}).
\begin{prop}
\label{prop10.7.1}
For given $\varepsilon> 0$, large enough $N$ and all $x$,
%
%
\begin{eqnarray}
\label{eq10.7.2}
E_{x}[|X(N^{3})|_{R}] - |x|_{R} &\le& \varepsilon N^{2} -
M_{2}(|x|_{K}/|x|)N^{3} \cdot
1\{|x|>N^{6}\} \nonumber\\[-8pt]\\[-8pt]
&&{} - \kappa_{N}N^{4} \cdot1\{|x|_{R}>\kappa^{2}_{N}N^{17}, |x|\le
N^{6}\},\nonumber
\end{eqnarray}
%
where $M_{2}$ is specified before (\ref{eq5.4.1}).
\end{prop}

The bound (\ref{eq10.7.2}) implies (\ref{eq5.4.1}), which was employed
in Section \ref{sec3}, together with bounds on $E_{x}[|X(N^{3})|_{L}]$
and $E_{x} [|X(N^{3})|_{A}]$, to obtain (\ref{eq5.3.10}) of Theorem
\ref {thm5.7.1}. The bound on $E_{x}[|X(N^{3})|_{A}]$ was derived
relatively quickly, whereas the bound on $E_{x}[|X(N^{3})|_{L}]$ will
require substantial estimation and will be derived in Sections
\ref{sec6}--\ref{sec10}. The bound on $E_{x}[|X(N^{3})|_{R}]$ that is
given here will require a moderate amount of work.

In order to show Proposition \ref{prop10.7.1}, it will be useful to rewrite
$|X(t)|_{R} - |x|_{R}$ as
%
%
\begin{equation}
\label{eq10.7.3}
|X(t)|_{R} - |x|_{R} = I_{R}(t) - D_{R}(t),
\end{equation}
where $I_{A}(t)$ and $D_{A}(t)$ are the nondecreasing functions
corresponding to the cumulative increase and decrease of
$|X(\cdot)|_{R}$ up to time $t$. A similar decomposition was used in
Section \ref{sec4} for $|X(t)|_{A}$. Here, $I_{R}(0) = D_{R}(0) = 0$,
with $I_{R}(t)$ being the jump process with
\[
I_{R}(t) - I_{R}(t-) = |X(t)|_{R} - |X(t-)|_{R}.
\]
%
One can check that $D_{R}(\cdot)$ is continuous except when a document
departs from a route. Its derivative is defined almost everywhere,
being defined except at the arrival or departure of a document. Since
$D_{R}(\cdot)$ is nondecreasing,
\[
D_{R}(t_{2}) - D_{R}(t_{1}) \ge\int^{t_{2}}_{t_{1}}D'_{R}(t)\,dt
\qquad\mbox{for } t_{1} \le t_{2}.
\]

It is easy to obtain a suitable upper bound on $E_{x}[I_{R}(N^{3})]$; a suitable
lower bound on $E_{x}[D_{R}(N^{3})]$ requires more effort. We therefore first
demonstrate Proposition \ref{prop10.7.4}, which analyzes $E_{x}[I_{R}(N^{3})]$.

As in Section \ref{sec4}, $A_{r}(t)$ denotes the cumulative number of
arrivals at route $r$ by time $t$. It follows from elementary renewal
theory that, for appropriate $C_{6}$ and $t \ge1$,
%
%
\begin{equation}
\label{eq10.1.2}
E_{x}[A_{r}(t)] \le C_{6}t \qquad\mbox{for each } r
\end{equation}
(see, e.g., \cite{r3}, page 136). Since large residual interarrival
times can only delay arrivals, the bound is uniform in $x$.
\begin{prop}
\label{prop10.7.4}
For given $\varepsilon>0$ and large enough $N$,
%
%
\begin{equation}
\label{eq10.7.5}
E_{x}[I_{R}(N^{3})] \le\varepsilon N^{2} \qquad\mbox{for all } x.
\end{equation}
%
\end{prop}
\begin{pf}
It follows from (\ref{eq5.3.2}) that the expected increase in
$I_{R}(\cdot)$, due to a document that arrives at route $r$, is
\[
M_{1}\kappa_{N,r}\int^{\infty}_{N_{H_{r}}}N_{r}(s)h^{*}_{r}(s)\,ds.
\]
Since the number of arriving documents by time $N^{3}$ and their initial
service times are independent, it follows that
%
%
\begin{equation}
\label{eq10.7.6}
E_{x}[I_{R}(N^{3})] = \biggl( M_{1}\kappa_{N,r}\int^{\infty}_{N_{H_{r}}}
N_{r}(s)h^{*}_{r}(s)\,ds \biggr)E_{x}[A_{r}(N^{3})].
\end{equation}

In order to bound the first term on the RHS of (\ref{eq10.7.6}), we
decompose the integral there into $\int^{N}_{N_{H_{r}}} +
\int^{\infty}_{N}$. When $N \ge N_{H_{r}}$, one has, by
(\ref{eq90.1.1}) and (\ref{eq3.8.1}),
%
%
\begin{eqnarray}
\label{eq10.7.7}
\kappa_{N,r}\int^{N}_{N_{H_{r}}}N_{r}(s)h^{*}_{r}(s)\,ds &=& \frac
{1}{\Gamma
(1/N^{4})}\int^{N}_{N_{H_{r}}}s h^{*}_{r}(s)\,ds \nonumber\\[-8pt]\\[-8pt]
&\le& \frac{N}{\Gamma(1/N^{4})}\bar{H}^{*}_{r}(N_{H_{r}}) \le(N^{3}
\Gamma
(1/N^{4}) )^{-1}.\nonumber
\end{eqnarray}
This is, for large enough $N$, at most $1/N^{2}$, because of the small power
$\gamma$ in the definition of $\Gamma(\cdot)$. Also,
%
%
\begin{eqnarray}
\label{eq10.7.8}
\kappa_{N,r}\int^{\infty}_{N}N_{r}(s)h^{*}_{r}(s)\,ds
&\le& \frac
{1}{N\Gamma
(1/N^{4})}\int^{\infty}_{N}s^{2}h^{*}_{r}(s)\,ds \nonumber\\
&\le& \frac{1}{N^{1+\delta_{1}/2}\Gamma(1/N^{4})}\int^{\infty}_{N}
s^{2+\delta_{1}/2}h_{r}^{*}(s)\,ds
\\
&\le&\frac{C_{7}}{N^{1+\delta
_{1}/2}\Gamma
(1/N^{4})}\nonumber
\end{eqnarray}
for appropriate $C_{7}$, with the last inequality holding because of
(\ref{eq5.2.1}). Since $\gamma\le\delta_{1}/24$, this is, for large
$N$, at most $1/N^{1+\delta_{1}/4}$. Together, the bounds for the two
integrals imply that, for large enough $N$,
%
%
\begin{equation}
\label{eq10.7.9}
M_{1}\kappa_{N,r}\int^{\infty}_{N_{H_{r}}}N_{r}(s) h^{*}_{r}(s)\,ds
\le
2/N^{1+\delta_{1}/4}.
\end{equation}

Application of (\ref{eq10.7.9}) and (\ref{eq10.1.2}) to (\ref
{eq10.7.6}), with $t=N^{3}$ in (\ref{eq10.1.2}), implies
(\ref{eq10.7.5}).
\end{pf}

We now derive a lower bound on $E_{x}[D_{R}(N^{3})]$. As in (\ref
{eq10.7.2}), we need to consider two separate cases, depending on
whether $|x|>N^{6}$ or both $|x|_{R} > \kappa^{2}_{N}N^{17}$ and
$|x|\le N^{6}$ hold. In both cases, we will employ the following lemma.
Recall that $M_2 = \frac{1}{8}C_{8}M_1$, with
$C_{8}=(1\wedge\min_{l}c_{l})(\min_{r,r'}(w_{r}/w_{r'}))$.
\begin{lem}
\label{lem10.8.1}
\textup{(a)} For all $t$,
%
%
\begin{equation}
\label{eq10.8.2}
D_{R}(t) \ge M_{1}\bigl(|x|_{K} - |X(t)|_{K}\bigr).
\end{equation}

\textup{(b)} For almost all $t$,
%
%
\begin{eqnarray}
\label{eq10.8.3}
D'_{R}(t) &\ge& \frac{C_{8}M_{1}}{|X
(t)|}\sum_{r} \kappa_{N,r}\int^{\infty}_{N_{H_{r}}} \biggl(\frac{s}{N}
\vee1 \biggr)Z^{*}_{r}(t,s)\,ds \nonumber\\[-8pt]\\[-8pt]
&\ge& 8M_{2}|X(t)|_{K}/
|X(t)|.\nonumber
\end{eqnarray}
\end{lem}
\begin{pf}
We first show (a). Recall that $\tilde{X}(\cdot)$ is the stochastic
process constructed from $X(\cdot)$ in Section \ref{intro}, where
service of documents is pathwise identical to $X(\cdot)$, but where the
arrival of documents is suppressed. One can check that, for all $t$ and
$\omega$,
%
%
\begin{equation}
\label{eq10.10.1}
|\tilde{X}(t)|_{K} \le|X(t)|_{K}
\end{equation}
and
%
%
\begin{equation}
\label{eq10.10.2}
D_{R}(t) \ge|x|_{R} - |\tilde{X}(t)|_{R}.
\end{equation}
Inequality (\ref{eq10.10.1}) follows immediately from $\tilde
{Z}^{*}(t,B)\le Z^{*}(t,B)$ for all $B \subseteq\mathbb{R}^{+}$. For
(\ref {eq10.10.2}), note that the LHS gives the cumulative decrease of
$|X(\cdot)|_{R}$ over $[0,t]$ due to the service of all documents,
whereas the RHS gives the decrease due to service of only the original
documents while ignoring the decrease due to service of new documents.

On account of (\ref{eq10.10.1}) and (\ref{eq10.10.2}), to show (\ref
{eq10.8.2})
it suffices to show
%
%
\begin{equation}
\label{eq10.10.3}
|x|_{R} - |\tilde{X}(t)|_{R} \ge M_{1}\bigl(|x|_{K} - |\tilde{X}(t)|_{K}\bigr).
\end{equation}
Substituting in the definition of $|\cdot|_{R}$ given by (\ref
{eq5.3.2}) and integrating by parts on the LHS of (\ref{eq10.10.3})
gives
%
%
\begin{eqnarray}
\label{eq10.10.4}
&& M_{1}\sum_{r}\kappa_{N,r}N_{r}(N_{H_{r}})\bigl(z^{*}_{r}((N_{H_{r}},
\infty)) -
\tilde{Z}^{*}_{r}(t, (N_{H_{r}}, \infty))\bigr) \nonumber\\[-8pt]\\[-8pt]
&&\qquad{} + M_{1}\sum_{r}\kappa_{N,r}\int^{\infty}_{N_{H_{r}}}N'_{r}(s)
\bigl(z^{*}_{r}((s, \infty)) - \tilde{Z}^{*}_{r}(t, (s,
\infty))\bigr)\,ds.\nonumber
\end{eqnarray}
It follows from (\ref{eq90.1.2}) and $N_{H_{r}}\ge1$ that
$N_{r}(N_{H_{r}}) \ge 1$ and
that $N'_{r}(s) \ge1$ for all $s$. Consequently,
(\ref{eq10.10.4}) is at least
\begin{eqnarray*}
&& M_{1}\sum_{r}\kappa_{N,r}\bigl(z^{*}_{r}((N_{H_{r}}, \infty)) -
\tilde{Z}^{*}_{r}
(t, (N_{H_{r}}, \infty))\bigr) \\
&&\qquad= M_{1}\bigl(|x|_{K} - |\tilde{X}(t)|_{K}\bigr),
\end{eqnarray*}
which implies (\ref{eq10.10.3}).

For (b), we first note that because of the weighted max--min fair
protocol and~(\ref{eq1.12.1}), the
rate at which each document is served is at least
%
%
\begin{equation}
\label{eq10.8.4}
\Bigl(\min_{l}c_{l}\Bigr)\Bigl(\min_{r,r'}(w_{r}/w'_{r})\Bigr)\big/|X(t)|.
\end{equation}
Moreover, the rate of decrease of $|X(t)|_{R}$ per unit service of each document
on route $r$ is at least
%
%
\begin{eqnarray}
\label{eq10.8.5}\quad
M_{1}\kappa_{N,r}\int^{\infty
}_{N_{H_{r}}}N'_{r}(s)Z^{*}_{r}(t,s)\,ds
&\ge&
M_{1}\kappa_{N,r}\int^{\infty}_{N_{H_{r}}} \biggl(\frac{s}{N}\vee1 \biggr)
Z^{*}_{r}(t,s)\,ds \nonumber\\[-8pt]\\[-8pt]
&\ge& M_{1}\kappa_{N,r}Z^{*}_{r}(t, (N_{H_{r}}, \infty)).\nonumber
\end{eqnarray}
Summing (\ref{eq10.8.5}) over $r$ and multiplying by (\ref{eq10.8.4}) gives
each of the bounds in (\ref{eq10.8.3}).
\end{pf}

We first derive a lower bound on $E_{x}[D_{R}(N^{3})]$ in the case
where $|x| >
N^{6}$.
%
\begin{prop}
\label{prop10.7.10}
For large enough $N$ and all $|x| > N^{6}$,
%
%
\begin{equation}
\label{eq10.7.11}
E_{x}[D_{R}(N^{3})] \ge M_{2}(|x|_{K}/|x|)N^{3}.
\end{equation}
%
\end{prop}
\begin{pf}
We restrict our attention to the set
\[
B_{1} = \{\omega\dvtx|X(t)| \le|x| + N^{6} \mbox{ for all } t \in
[0,N^{3}]\}.
\]
By applying Markov's inequality to inequality (\ref{eq10.1.2}) with
$t=N^{3}$, one has that, for large $N$,
%
%
\begin{equation}
\label{eq10.2.2}
P_{x} \biggl(\sum_{r}A_{r}(N^{3})>N^{6} \biggr) \le\frac{C_{6}}{N^{3}}|\mathcal{R}|
\le\frac{1}{2}
\end{equation}
for all $x$. Consequently,
%
%
\begin{equation}
\label{eq10.2.1}
P(B_{1}) \ge1/2.
\end{equation}
This bound does not depend on $|x|$.

We now consider two cases, depending on whether the set
\[
B_{2} = \bigl\{\omega\dvtx|X(t)|_{K} > \tfrac{1}{2}|x|_{K} \mbox{ for all } t
\in
[0,N^{3}] \bigr\}
\]
occurs. Since $|x|>N^{6}$, it follows from the second half of (\ref{eq10.8.3})
that, for all $t \in[0, N^{3}]$,
\[
D'_{R}(t) \ge2M_{2}|x|_{K}/|x|
\]
on $B_{1} \cap B_{2}$. Consequently, on $B_{1} \cap B_{2}$,
%
%
\begin{equation}
\label{eq10.2.6}
D_{R}(N^{3}) \ge2M_{2} (|x|_{K}/|x| ) N^{3}.
\end{equation}

On the other hand, on $B_{1} \cap B_{2}^{c}$,
%
%
\begin{equation}
\label{eq10.7.12}
|x|_{K} - |X(\tau)|_{K} \ge\tfrac{1}{2}|x|_{K}
\end{equation}
for some (random) $\tau\in[0, N^{3}]$. By (\ref{eq10.8.2}),
\[
D_{R}(t) \ge M_{1}\bigl(|x|_{K} - |X(t)|_{K}\bigr)
\]
for all $t$. Together with (\ref{eq10.7.12}), this implies that
%
%
\begin{equation}
\label{eq10.3.4}
D_{R}(N^{3}) \ge D_{R}(\tau) \ge\tfrac{1}{2}M_{1}|x|_{K} \ge2M_{2}
(|x|_{K}/|x| )N^{6},
\end{equation}
where $|x|>N^{6}$ was used in the last inequality.

Together, (\ref{eq10.2.6}) and (\ref{eq10.3.4}) imply that, on $B_{1}$,
\[
D_{R}(N^{3}) \ge2M_{2} (|x|_{K}/|x| ) N^{3}.
\]
Inequality (\ref{eq10.7.11}) follows from this and (\ref{eq10.2.1}).
\end{pf}

We now 
derive a lower bound on $E_{x}[D_{R}(N^{3})]$ in the case where
$|x|_{R} >
\kappa^{2}_{N}N^{17}$ and $|x| \le N^{6}$ both hold. We note that,
starting from
(\ref{eq10.5.1}), the argument relies on the discreteness of
documents. If one
wishes to employ a fluid limit based argument rather than the discrete setting
employed in this paper, different reasoning will be required at this
point; it is not
obvious how one would proceed.
\begin{prop}
\label{prop10.7.14}
For large enough $N$,
%
%
\begin{equation}
\label{eq10.7.15}
E_{x}[D_{R}(N^{3})] \ge\kappa_{N}N^{4}
\end{equation}
for all $|x|_{R} > \kappa^{2}_{N}N^{17}$ and $|x| \le N^{6}$.
\end{prop}
\begin{pf}
As in the proof of Proposition \ref{prop10.7.10}, we restrict attention
to the set $B_{1}$ defined there. The bound $P(B_{1}) \ge1/2$ in
(\ref{eq10.2.1}) continues to hold here. In our present setting, since
$|x| \le N^{6}$, $\omega \in B_{1}$ implies that
\[
|X(t)| \le2N^{6} \qquad\mbox{for all } t \in[0, N^{3}].
\]
We also consider two cases, depending on whether
\[
B_{3} = \bigl\{\omega\dvtx|X(t)|_{R} > \tfrac{1}{2}\kappa^{2}_{N}N^{17}
\mbox{ for all } t \in[0, N^{3}] \bigr\}
\]
occurs.

The case $B^{c}_{3}$ is almost immediate. 
It follows from (\ref{eq10.7.3}) that, for large enough $N$ and for
some $\tau \in(0, N^{3}]$,
%
%
\begin{equation}
\label{eq10.7.16}
D_{R}(N^{3}) \ge D_{R}(\tau) \ge|x|_{R} - |X(\tau)|_{R} \ge\tfrac{1}{2}
\kappa^{2}_{N}N^{17} > 2\kappa_{N}N^{4}
\end{equation}
for $\omega\in B^{c}_{3}$.

The case $B_{3}$ requires some work. We first note that, by the first
part of
(\ref{eq10.8.3}),
%
%
\begin{eqnarray}
\label{eq10.7.17}
D'_{R}(t) &\ge& \frac{C_{8}M_{1}}{|X(t)|}
\sum_{r}\kappa_{N,r} \biggl(\int^{\infty}_{N_{H_{r}}} \biggl(\frac{s}{N}
\vee1 \biggr)Z^{*}_{r}(t,s)\,ds \biggr) \nonumber\\[-8pt]\\[-8pt]
&\ge& \frac{C_{8}M_{1}}{2N^{6}}
\sum_{r}\kappa_{N,r} \biggl(\int^{\infty}_{N_{H_{r}}} \biggl(\frac{s}{N}
\vee1
\biggr)Z^{*}_{r}(t,s)\,ds \biggr),\nonumber
\end{eqnarray}
when $\omega\in B_{1}$.

We will truncate the second integral in (\ref{eq10.7.17}) in order to
be able to introduce an additional factor $s$ into the integrand. We
first note that, since $\Phi(0)=0$, if a document with residual service
time at least $s$ is present at time $t$ on some route $r$, then, for
large $N$,
%
%
\begin{equation}
\label{eq10.5.1}
|X(t)|_{R} \ge M_{1}\kappa_{N,r}s^{2}/N \ge M_{1}s^{2}/N.
\end{equation}
Hence, there are no documents with residual service time
%
%
\begin{equation}
\label{eq10.5.2}
s > s_{1} \stackrel{\mathrm{def}}{=} \bigl((N/M_{1})|X(t)|_{R} \bigr)^{1/2}.
\end{equation}
%
It follows that, for appropriate $C_{9}> 0$, (\ref{eq10.7.17}) is at least
%
%
\begin{eqnarray}
\label{eq10.5.3}
&&\frac{C_{8}M_{1}}{2N^{6}}\sum_{r}
\kappa_{N,r} \biggl(\int^{s_{1}+1}_{N_{H_{r}}} \biggl(\frac{s}{N}\vee1 \biggr)
Z^{*}_{r}(t,s)\,ds \biggr) \nonumber\\
&&\qquad\ge\frac{C_{8}M_{1}^{3/2}}{4N^{13/2}
|X(t)|^{1/2}_{R}} \sum_{r}\kappa_{N,r} \biggl(\int^{s_{1}+1}_{N_{H_{r}}}
N_{r}(s)Z^{*}_{r}(t,s)\,ds \biggr) \nonumber\\[-8pt]\\[-8pt]
&&\qquad\ge\frac{2C_{9}M_{1}^{3/2}}{N^{13/2}|X(t)|^{1/2}_{R}}\sum
_{r}\kappa_{N,r}
\biggl(\int^{\infty}_{{N_{H_{r}}}}N_{r}(s)Z^{*}_{r}(t,s)\,ds \biggr) \nonumber\\
&&\qquad= \frac{2C_{9}M_{1}^{1/2}}{N^{13/2}}|X(t)|^{1/2}_{R} \ge C_{9}M_{1}^{1/2}
\kappa_{N}N^{2}\nonumber
\end{eqnarray}
for all $t \in[0,N^{3}]$. The exponential tail of $\Phi(\cdot)$ is
used in the
last inequality; the equality relies on $\omega\in B_{3}$.

Employing the bound on $D'_{R}(t)$ obtained from (\ref{eq10.7.17}) and
(\ref{eq10.5.3}), and integrating over $t \in[0,N^{3}]$, it follows that,
for large $N$,
\[
D_{R}(N^{3}) \ge C_{9}M_{1}^{1/2}\kappa_{N}N^{5} > 2\kappa_{N}N^{4}
\]
on $B_{1} \cap B_{3}$. Together with (\ref{eq10.7.16}), this implies that
$D_{R}(N^{3}) > 2\kappa_{N}N^{4}$ on $B_{1}$. Inequality (\ref{eq10.7.15})
follows from this and $P(B_{1}) \ge1/2$.
\end{pf}

Proposition \ref{prop10.7.1} follows immediately from (\ref {eq10.7.3})
and Propositions \ref{prop10.7.4}, \ref{prop10.7.10} and~\ref{prop10.7.14}.

\section{Upper bounds on $E_{x}[|X(N^{3})|_{L}]$: Basic layout and
bounds on exceptional sets}\label{sec6}

In this section, we begin our investigation of upper bounds on
$E_{x}[|X( N^{3})|_{L}] - |x|_{L}$. Since these bounds will require us
to examine a number of subcases in Sections \ref{sec6}--\ref{sec9}, we
will only arrive at the desired bounds in Section \ref{sec10}. In the
current section, we first state certain elementary inequalities, mostly
involving $|\cdot|_{r,s}$, that will be useful later on. We then define
the ``good'' sets $\mathcal{A}(\cdot)$ of realizations of $X(\cdot)$ to
which our bounds in Sections \ref{sec7}--\ref{sec9} will apply. The
remainder of the section is spent demonstrating Proposition
\ref{prop21.1.1}, which gives an upper bound on
$E_{x}[|X(t)|_{L}-|x|_{L}$; $\mathcal{A}(t)^{c}]$, where
$\mathcal{A}(t)^{c}$ is the small exceptional set.

\subsection*{Elementary inequalities}
Here we state a number of elementary inequalities that will be useful later
on. Let $z_{i}(\cdot)$, $i=1,2,3$, denote configurations of particles on
$\mathbb{R}^{+}$, with $z_{i}(B)$ denoting the number of particles (or
documents) in $B \subseteq\mathbb{R}^{+}$. If one assumes
%
%
\begin{equation}
\label{eq39.4.1}
z_{3}(B) = z_{1}(B) + z_{2}(B) \qquad\mbox{for all } B \subseteq\mathbb{R}^{+},
\end{equation}
it follows that
%
%
\begin{equation}
\label{eq39.4.2}
z^{*}_{3}(B) = z^{*}_{1}(B) + z^{*}_{2}(B) \qquad\mbox{for all } B
\subseteq \mathbb{R}^{+},
\end{equation}
where $z_{i}^{*}(B)$ is defined analogously to $z^{*}_{r}(B)$ below
(\ref{eq5.1.6}), with convolution being with respect to $\phi(\cdot )$.
Several elementary equalities follow from (\ref{eq39.4.2}), including
%
%
\begin{equation}
\label{eq39.4.15}
|x_{3}|_{r,s} = |x_{1}|_{r,s} + |x_{2}|_{r,s} \qquad\mbox{for all } r \in
\mathcal{R} \mbox{ and } s > 0,
\end{equation}
where $x_{i}$ are states in 
the metric space $S$ introduced in Section \ref{sec2} for which the
analog of
(\ref{eq39.4.1}) is satisfied for each $r$ and $|\cdot|_{r,s}$ is
given by
(\ref{eq5.3.1}).
%

Recall that $\tilde{X}(\cdot)$ and $X^{A}(\cdot)$ are the processes
constructed from $X(\cdot)$ that were introduced in Section
\ref{intro}, where service of each document is pathwise identical to
$X(\cdot)$, but where, for $\tilde{X}(\cdot)$, the arrival of documents
is suppressed and, for $X^{A}(\cdot)$, only new documents are included.
One has
\[
Z(t,B) = \tilde{Z}(t,B) + Z^{A}(t,B) \qquad\mbox{for } t \ge0 \mbox{ and }
B \subseteq\mathbb{R}^{+},
\]
where the processes $Z(\cdot), \tilde{Z}(\cdot)$, and $Z^{A}(\cdot )$
correspond to $X(\cdot), \tilde{X}(\cdot)$ and $X^{A}(\cdot)$. From
(\ref{eq39.4.2}),
%
%
\begin{equation}
\label{eq39.4.16}
Z^{*}(t,B) = \tilde{Z}^{*}(t,B) + Z^{A,*}(t,B)
\qquad\mbox{for } t \ge0 \mbox{ and } B \subseteq\mathbb{R}^{+},
\end{equation}
and from (\ref{eq39.4.15}),
%
%
\begin{equation}
\label{eq39.1.3}
|X(t)|_{r,s} = |\tilde{X}(t)|_{r,s} + |X^{A}(t)|_{r,s}\qquad
\mbox{for } t \ge0, r \in\mathcal{R}, s > 0.
\end{equation}

Another elementary equality involving $X(\cdot)$ is given by
%
%
\begin{equation}
\label{eq39.1.4}
\tilde{Z}_{r}(t,B) = z_{r}\bigl(B + \Delta_{r}(t)\bigr)
\qquad\mbox{for } t \ge0, r \in \mathcal{R}, B \subseteq\mathbb{R}^{+},
\end{equation}
where, we recall, $\Delta_{r}(t)$ is the translation that gives the
amount of service an original document that has not yet completed
service has received by time $t$. The equality relies on all documents
on a given route $r$ receiving equal service at each time. [If
$\tilde{Z}_{r}(t,\mathbb{R}^{+})=0$, set $\Delta_{r}(t)=\infty$ and
$z_{r}(\mathbb{R}^{+} + \infty) = 0$.] From (\ref{eq39.1.4}), one
obtains
%
%
\begin{equation}
\label{eq39.1.5}
\tilde{Z}^{*}_{r}(t,B) \le z^{*}_{r}\bigl(B +
\Delta_{r}(t)\bigr) \qquad\mbox{for } t \ge 0, r \in\mathcal{R}, B
\subseteq\mathbb{R}^{+};
\end{equation}
the inequality arises from the possibility that original documents have
completed service by time $t$.

A consequence of (\ref{eq5.3.1}) and (\ref{eq39.1.5}) is that
%
%
\begin{equation}
\label{eq39.1.6}
|\tilde{X}(t)|_{r,s} \le|x|_{r,s+\Delta_{r}(t)} \qquad
\mbox{for }t\ge0, r \in \mathcal{R}, s > 0.
\end{equation}
Combining (\ref{eq39.1.3}) and (\ref{eq39.1.6}) produces
%
%
\begin{equation}
\label{eq39.1.7}
|X(t)|_{r,s} \le|x|_{r,s+\Delta_{r}(t)} +
|X^{A}(t)|_{r,s} \qquad\mbox{for } t \ge0, r \in\mathcal{R}, s > 0;
\end{equation}
taking the supremum over all $r$ and $s$ therefore gives
%
%
\begin{equation}
\label{eq39.2.1}
|X(t)|_{L} \le|x|_{L} + |X^{A}(t)|_{L} \qquad\mbox{for all } t \ge0.
\end{equation}
Application of (\ref{eq39.1.5}) also implies
%
%
\begin{equation}
\label{eq39.1.8}
\tilde{Z}^{*}_{r}(t,s) \le z^{*}_{r}\bigl(s+\Delta_{r}(t)\bigr)\qquad
\mbox{for } t \ge0, r \in\mathcal{R}, B \subseteq\mathbb{R}^{+},
\end{equation}
and application of (\ref{eq39.1.5}), together with (\ref{eq39.4.16}), implies
that
%
%
\begin{equation}
\label{eq39.1.9}
|X(t)|_{2} \le|x|_{2} + |X^{A}(t)|_{2} \qquad\mbox{for } t \ge0,
\end{equation}
where $|\cdot|_{2}$ is given in (\ref{eq5.3.11}). The term on the LHS
of (\ref{eq10.8.5}) can also be derived using (\ref{eq39.1.8}).

\subsection*{The sets $\mathcal{A}(t)$}

In this subsection, we define the random set $\mathcal{A}(t)$, which is
a function of $X(t')$, for $t' \in[0,t]$. In Sections \ref{sec7}--\ref
{sec10}, we will establish upper bounds on $|X(N^{3})|_{r,s}$ for
$\omega\in\mathcal {A}(N^{3})$; the exceptional small set
$\mathcal{A}(N^{3})^{c}$ will be treated in the next subsection. The
set $\mathcal{A}(t)$ will be a ``good'' set in the sense that the
number of arrivals over $[0,t]$, for given $t$, is restricted by upper
bounds, which will enable us to show that $|X(\cdot)|_{L}$ decreases in
an appropriate manner.

The set $\mathcal{A}(t)$ is given by $\mathcal{A}(t) = \mathcal
{A}_{1}(t) \cap \mathcal{A}_{2}(t)$, with
%
%
\begin{equation}
\label{eq20.1.5}
\mathcal{A}_{i}(t) = \bigcap_{r,j} \mathcal{A}_{i,r,j}(t)\qquad
\mbox{for } i=1,2,
\end{equation}
where $\mathcal{A}_{i,r,j}(t)$ specify upper bounds on the numbers of weighted
arrivals of documents with different service times. To define
$\mathcal{A}_{i,r,j}(t)$, we denote by $v_{0}, v_{1}, \ldots, v_{J}$ the
increasing sequence with
%
%
\begin{equation}
\label{eq20.1.6}
v_{j+1} = v_{j} + 1/b^{3} \qquad\mbox{for } j = 0, \ldots, J-1,
\end{equation}
with $v_{0}=0$ and $v_{J}=N+1$, and where $b$ is as in (\ref
{eq5.1.4}). Note that it follows from the second half of
(\ref{eq5.1.10}) that, for $b \ge2$,
%
%
\begin{equation}
\label{eq20.1.7}
\bar{H}^{*}_{r}(v_{j+1})/\bar{H}^{*}_{r}(v_{j}) \ge1/2
\qquad\mbox{for all } r \mbox{ and } j.
\end{equation}
We also denote by $S^{1}_{r}(k)$, $k = 1, \ldots, A_{r}(t)$,
the service time of the $k$th arrival at route $r$, where $A_{r}(t)$
is the cumulative number of arrivals at $r$ by time $t$.

We set, for $r \in\mathcal{R}$ and $j = 0, \ldots, J$,
%
%
\begin{equation}
\label{eq20.1.2}
\mathcal{A}_{1,r,j}(t) = \Biggl\{
\omega\dvtx\sum^{A_{r}(t)}_{k=1} \bar{\Phi }\bigl(v_{j} -
S^{1}_{r}(k)\bigr) \le2\nu_{r}\bigl(\bar{H}^{*}_{r}(v_{j})t \vee t^{\eta}\bigr)
\Biggr\}.
\end{equation}
Here, we assume $\eta\in(0, 1/12]$, and, as elsewhere, we set $\bar{H}_{r}
(\cdot) = 1 - H_{r}(\cdot)$ and $\bar{\Phi}(\cdot) = 1 - \Phi
(\cdot)$.
One has, as a special case of (\ref{eq20.1.2}), that
%
%
\begin{equation}
\label{eq20.2.1}
A_{r}(t) \le2\nu_{r}t \qquad\mbox{on } \mathcal{A}_{1,r,0}(t).
\end{equation}

Since
%
%
\begin{equation}
\label{eq20.1.3}
E\bigl[\bar{\Phi}\bigl(v_{j} - S^{1}_{r}(k)\bigr)\bigr] = \int^{\infty}_{0}
\bar{\Phi}(v_{j}-s)\,dH_{r}(s) = \bar{H}^{*}_{r}(v_{j})
\end{equation}
and $A_{r}(t) \sim\nu_{r}t$ for large $t$, the probability of the
complement $\mathcal{A}_{1,r,j}(t)^{c}$ can be bounded above by using
standard large derivation estimates. The term $t^{\eta}$ is included on
the RHS of (\ref{eq20.1.2}) so that, when $\bar{H}^{*}_{r}(v_{j})$ is
small, the probability
of the event remains small. 

We also set, for $r \in\mathcal{R}$ and $j = 0, \ldots, J$,
%
%
\begin{equation}
\label{eq20.2.2}\quad
\mathcal{A}_{2,r,j}(t) = \Biggl\{
\omega\dvtx\sum^{A_{r}(t)}_{k=1} \phi\bigl( v_{j} - S^{1}_{r}(k)\bigr)
\le(1+\varepsilon_{5})\nu_{r}\bigl(h^{*}_{r}(v_{j})t\vee t^{\eta}\bigr) \Biggr\},
\end{equation}
where $\varepsilon_{5}> 0$. Analogous to (\ref{eq20.1.3}), one has
%
%
\begin{equation}
\label{eq20.2.4}
E\bigl[\phi\bigl(v_{j} - S^{1}_{j}(k)\bigr)\bigr] =
\int^{\infty}_{0}\phi(v_{j} - s) \, dH_{r}(s) = h^{*}_{r}(v_{j}).
\end{equation}
The probabilities $P_{x}(\mathcal{A}_{2,r,j}(t)^{c})$ will satisfy
large deviation bounds as well. The constant $\varepsilon_{5}$ here
will later be required to satisfy
$\varepsilon_{5}\le\varepsilon_{7}/4$, where $\varepsilon_{7}$ is
specified in (\ref {eq17.1.5}) and measures how subcritical the network
is. In (\ref{eq20.1.2}), we only need to employ the constant 2, rather
than $1+\varepsilon_{5}$ as in (\ref {eq20.2.2}), because
(\ref{eq20.1.2}) will be applied to the right tail of $\bar
{H}^{*}_{r}(\cdot)$, rather than the ``main body'' of
$H^{*}_{r}(\cdot)$, as will (\ref{eq20.2.2}).

\subsection*{Upper bounds on $\mathcal{A}(t)^{c}$}

The main result in this last subsection is the following proposition.
\begin{prop}
\label{prop21.1.1}
For large enough $t$,
%
%
\begin{equation}
\label{eq21.1.2}
E_{x}[|X(t)|_{L} - |x|_{L}; \mathcal{A}(t)^{c}] \le
N^{3} e^{-C_{10}t^{\eta}}
\end{equation}
for all $N, x$ and appropriate $C_{10}> 0$.
\end{prop}

Proposition \ref{prop21.1.1} gives strong bounds on the growth of
$|X(t)|_{L}$ on $\mathcal{A}(t)^{c}$. This behavior is primarily due to
the small probability $P_{x}(\mathcal{A}(t)^{c})$, which is given in
the next proposition.
\begin{prop}
\label{prop21.1.3}
For large enough $t$,
%
%
\begin{equation}
\label{eq21.1.4}
P_{x}(\mathcal{A}(t)^{c}) \le N e^{-C_{11}t^{\eta}}
\end{equation}
for all $N, x$ and appropriate $C_{11}> 0$.
\end{prop}

The interarrival times are assumed to be independent, and large initial
residual interarrival times only delay future arrivals. The initial
state $x$ will therefore not affect the bounds in (\ref{eq21.1.2}) and
(\ref{eq21.1.4}). Note that only the arrival process $A(\cdot)$ is
relevant for the bounds in (\ref{eq21.1.4}).

Proposition \ref{prop21.1.3} will serve as the main step in
demonstrating Proposition \ref{prop21.1.1}; it will also be used along
with Proposition \ref{prop21.1.1} in Section \ref{sec10}. When we apply
(\ref{eq21.1.2}) and (\ref{eq21.1.4}) there, we will set $t = N^{3}$
and so the factors $N^{3}$ and $N$ can be absorbed into the
corresponding exponentials. We note that $C_{10}$ and $C_{11}$ in
(\ref{eq21.1.2}) and (\ref{eq21.1.4}), and the bound on $t$ depend on
our choices of $\varepsilon_{5}$ and $b$, and on $\nu_{r}$ and $w_{r}$.

In order to show Proposition \ref{prop21.1.3}, we will employ
elementary large deviation estimates, which are given in the following
two lemmas.
\begin{lem}
\label{lem21.1.5}
Let $W(1), W(2), \ldots$ denote nonnegative i.i.d.
random variables with mean $\beta< \infty$. Then, for each
$\varepsilon> 0$, there exists $C_{12}> 0$, so that
%
%
\begin{equation}
\label{eq21.2.1}
P \Biggl( \sum^{n}_{k=1} W(k) \le(1 - \varepsilon) \beta n \Biggr) \le
e^{-C_{12}n}.
\end{equation}
When the support of $W(1)$ is contained in $[0,1]$ and
$\varepsilon\in(0,1]$,
%
%
\begin{equation}
\label{eq21.2.6}
P \Biggl( \sum^{n}_{k=1} W(k) \ge(1 + \varepsilon) \beta n \Biggr)
\le e^{{-C_{13}}{\varepsilon^{2} \beta n}},
\end{equation}
where $C_{13}> 0$ does not depend on the distribution of $W(1)$ or on
$\varepsilon$.
\end{lem}
\begin{pf}
Both (\ref{eq21.2.1}) and (\ref{eq21.2.6}) are elementary large
deviation bounds. We summarize the argument for (\ref{eq21.2.6}); (\ref
{eq21.2.1}) can be shown directly or by applying (\ref{eq21.2.6}) after
truncating $W(k)$.

As usual, one employs the moment generating function
%
%
\begin{equation}
\label{eq21.2.2}
\psi_{\theta}(n) = E
\bigl[e^{\theta\sum^{n}_{k=1}(W(k)-\beta)} \bigr]\qquad \mbox{for } \theta> 0.
\end{equation}
By expanding the exponential for $n=1$, it follows that for appropriate
$C_{14} \ge1$ and for $\theta\in(0,1]$,
%
%
\begin{equation}
\label{eq21.2.3}
\psi_{\theta}(1) \le1 + C_{14}\beta\theta^{2},
\end{equation}
and hence
%
%
\begin{equation}
\label{eq21.2.4}
\psi_{\theta}(n) \le(1 + C_{14}\beta\theta^{2})^{n}
\le e^{C_{14} \beta\theta^{2} n}.
\end{equation}
By applying Markov's inequality and setting $\theta=
\varepsilon/2C_{14} $, it follows that the LHS of (\ref{eq21.2.6}) is
at most
%
%
\begin{equation}
\label{eq21.2.5}
e^{-\varepsilon\beta\theta n}\psi_{\theta}(n) \le
e^{-\varepsilon ^{2}\beta n/4C_{14}} \le e^{-C_{13}\varepsilon^{2}\beta n}
\end{equation}
for $C_{13}= 1/4 C_{14}$, as desired.
\end{pf}

Let $W(1), W(2), \ldots$ denote the successive interarrival times for a
renewal process (with delay), with $A(t)= \max\{n\dvtx\sum^{n}_{k=1}
W(k)\le t\}$ denoting the number of renewals by time $t$. Here, $W(2),
W(3), \ldots $ are i.i.d., with $W(1)$ being the residual interarrival
time. We also introduce i.i.d. random variables $Y(1), Y(2), \ldots,$
with $Y(1) \in[0,1]$ that are defined on the same space as $W(k)$. Set
$E[W(2)]=\beta> 0$ and $E[Y(1)] = m$.
\begin{lem}
\label{lem21.3.1}
Let $W(1), W(2), \ldots$ and $Y(1), Y(2), \ldots$ be
as above. Then, for given $\varepsilon\in(0,1]$ and large $t$,
%
%
\begin{equation}
\label{eq21.3.2}
P \Biggl(\sum^{A(t)}_{k=1}Y(k) >
(1+\varepsilon)\beta^{-1}mt \Biggr) \le e^{-C_{15}mt},
\end{equation}
where $C_{15}> 0$ does not depend on the distribution of $Y(1)$.
\end{lem}
\begin{pf}
$\{A(t) \ge n \}$ is contained in the event $ \{\sum^{n}_{k=1}W(k)\le t
\}$. Consequently, by (\ref{eq21.2.1}) of Lemma \ref{lem21.1.5}, substitution of
$\varepsilon/3$ for $\varepsilon$ there implies that, for $n(t) =
\lceil (1 - \varepsilon/3)^{-1}\beta^{-1}t\rceil$,
%
%
\begin{equation}
\label{eq21.3.3}
P\bigl(A(t) > n(t)\bigr) \le P
\Biggl(\sum^{n(t)+1}_{k=2}W(k) \le t \Biggr) \le e^{-C_{16}t}
\end{equation}
for appropriate $C_{16}> 0$ and large $t$ (which may depend on
$\varepsilon $ and the distribution of $W$).

We next consider the set where $A(t) \le n(t)$. It follows from (\ref
{eq21.2.6}) of Lem\-ma~\ref{lem21.1.5} that
%
%
\begin{eqnarray}
\label{eq21.3.4}
&&
P \Biggl(\sum^{A(t)}_{k=1}Y(k) > (1+\varepsilon)\beta^{-1}mt; A(t)\le
n(t) \Biggr)
\nonumber\\[-8pt]\\[-8pt]
&&\qquad
\le P \Biggl(\sum^{n(t)}_{k=1}Y(k) > (1+\varepsilon)\beta^{-1}mt \Biggr) \le
e^{-C_{13}\varepsilon^{2}\beta^{-1}mt/9}.\nonumber
\end{eqnarray}
Inequality (\ref{eq21.3.2}) follows from (\ref{eq21.3.3}) and
(\ref{eq21.3.4}).
\end{pf}

We now employ Lemma \ref{lem21.3.1} to prove Proposition \ref{prop21.1.3}.
\begin{pf*}{Proof of Proposition \protect\ref{prop21.1.3}}
We first note that since $\mathcal{A}(t) = \mathcal{A}_{1}(t) \cap
\mathcal{A}_{2}(t)$, with $\mathcal{A}_{i}(t) = \bigcap_{r \in
\mathcal{R}} \bigcap_{j=0}^{J} \mathcal{A}_{i,r,j}(t)$, where $J =
b^{3}(N+1)+1 \le2b^{3}N$, it suffices to show that for each
$\mathcal{A}_{i,r,j}(t)$,
%
%
\begin{equation}
\label{eq21.4.1}
P_{x}(\mathcal{A}_{i,r,j}(t)^{c}) \le e^{-C_{17}t^{\eta}}
\end{equation}
for $t \ge t_{0}$, for some fixed $t_{0}$ and appropriate $C_{17}> 0$.

We consider the case where $i=1$. Denote by $W(1), W(2), \ldots$ the
interarrival times of documents on route $r$ and set $Y(k) = \bar{\Phi}(v_{j}
- S^{1}_{r}(k))$. Then, $Y(k)$ are i.i.d. random variables and, except for
$W(1)$, so are $W(k)$. One has
%
%
\begin{equation}
\label{eq21.4.2}
\beta\stackrel{\mathrm{def}}{=} E[W(2)]=\nu^{-1}_{r} \quad\mbox{and}\quad m
\stackrel{\mathrm{def}}{=}E[Y(1)]=\bar{H}^{*}_{r}(v_{j})
\end{equation}
with the last equality following from (\ref{eq20.1.3}).

We break the problem into two cases, depending on whether or not
$\bar{H}^{*}_{r}(v_{j}) \ge t^{\eta- 1}$, in each case applying Lemma
\ref{lem21.3.1}, with $\varepsilon= 1$. Under $\bar{H}^{*}_{r} (v_{j})
\ge t^{\eta- 1}$, one has
%
%
\begin{equation}
\label{eq21.4.3}\qquad
P_{x}(\mathcal{A}_{1,r,j}(t)^{c}) = P_{x} \Biggl(\sum^{A_{r}(t)}_{k=1}\bar
{\Phi}
\bigl(v_{j}-S^{1}_{r}(k)\bigr) > 2\nu_{r}\bar{H}^{*}_{r} (v_{j})t \Biggr) \le
e^{-C_{15}t^{\eta}}
\end{equation}
for large $t$ and $C_{15}> 0$ as in the lemma, where neither depends on the
particular value of $\bar{H}^{*}_{r} (v_{j})$.

For $\bar{H}^{*}_{r} (v_{j}) < t^{\eta-1}$, we replace the random
variables defined above (\ref{eq21.4.2}) by i.i.d. random variables
$Y'(k) \in(0,1]$, with $Y'(k)\ge Y(k)$ and $E[Y'(k)] = t^{\eta- 1}$.
Then, again applying Lemma \ref{lem21.3.1}, but this time to $Y'(k)$,
$k=1,2,\ldots,$
%
%
\begin{equation}
\label{eq21.5.5}
P_{x}(\mathcal{A}_{1,r,j}(t)^{c}) \le P_{x} \Biggl(\sum
^{A_{r}(t)}_{k=1}Y'(k) > 2
\nu_{r} t^{\eta} \Biggr) \le e^{-C_{15}t^{\eta}}
\end{equation}
as before. Together with (\ref{eq21.4.3}), this implies (\ref
{eq21.4.1}) for $i=1$, with $C_{17}= C_{15}$.

The reasoning for (\ref{eq21.4.1}) when $i=2$ is the same, except that
one now sets $Y(k) = \phi(v_{j} - S^{1}_{r}(k))$, from which one
obtains
%
%
\begin{equation}
\label{eq21.5.3}
m\stackrel{\mathrm{def}}{=}E[Y(1)] = h^{*}_{r}(v_{j}).
\end{equation}
Also, the coefficient 2 on the RHS of (\ref{eq20.1.2}) is replaced by
the coefficient $1+\varepsilon_{5}$ in~(\ref{eq20.2.2}). Setting
$\varepsilon= \varepsilon_{5} $ in Lemma \ref{lem21.3.1}, one obtains
%
%
\begin{equation}
\label{eq21.5.4}
P_{x}\Biggl(\sum^{A_{r}(t)}_{k=1}\phi\bigl(v_{j}-S^{1}_{r}(k)\bigr) > (1+\varepsilon
_{5})\nu_{r} \bigl(h^{*}_{r}(v_{j})\vee t^{\eta}\bigr) \Biggr) \le
e^{-C_{15}t^{\eta}}
\end{equation}
for large $t$ and appropriate $C_{15}> 0$, chosen as in the lemma.
Setting $C_{17}= C_{15}$, one obtains (\ref{eq21.4.1}) for $i=2$ as
well.
\end{pf*}

Setting $|A(t)| = \sum_{r}A_{r}(t)$, where $A_{r}(t)$ is the number of
arrivals at each route by time $t$, it follows from elementary renewal
theory that for appropriate $C_{18}$ and $t \ge1$,
%
%
\begin{equation}
\label{eq21.6.2}
E [|A(t)|^{2} ] \le C_{18}t^{2}
\end{equation}
(see, e.g., \cite{r3}, page 136). Inequality (\ref{eq21.6.2}) is not
difficult to show by applying a standard truncation argument.

Here and later on, we will also use the two inequalities
%
%
\begin{equation}
\label{eq6.88.1}
z^{*}_{r}(s) \le bz^{*}_{r}((s,\infty)) \qquad\mbox{for all } s > 0,
\end{equation}
and
%
%
\begin{equation}
\label{eq6.88.2}
\Gamma\bigl(\bar{H}^{*}_{r}(N_{H_{r}}+1)\bigr)/\Gamma(\bar
{H}^{*}_{r}(N_{H_{r}})) \ge e^{-b},
\end{equation}
which follow from the definition of $\phi(\cdot)$ and the second inequality
in (\ref{eq5.1.10}). Employing Proposition \ref{prop21.1.3} and these
inequalities, we now demonstrate Proposition~\ref{prop21.1.1}.
\begin{pf*}{Proof of Proposition \protect\ref{prop21.1.1}}
By H\"{o}lder's inequality,
%
%
\begin{equation}
\label{eq21.7.1}\qquad
E_{x}[|X(t)|_{L} - |x|_{L}; \mathcal{A}(t)^{c}]
\le\sqrt{P_{x} (\mathcal{A}(t)^{c})} \sqrt{E_{x} \bigl[\bigl(|X(t)|_{L}
- |x|_{L}\bigr)^{2} \bigr]}.
\end{equation}
Also, by Proposition \ref{prop21.1.3}, one has
%
%
\begin{equation}
\label{eq21.7.2}
\sqrt{P_{x}(\mathcal{A}(t)^{c})} \le\sqrt{N} e^{-C_{11}t^{\eta/2}}
\end{equation}
for all $N,x$ and appropriate $C_{11}> 0$. So it remains to bound the
expectation on the RHS of (\ref{eq21.7.1}).

It follows from the definitions of $|\cdot|_{L}$, 
$\phi(\cdot)$ and $\Gamma(\cdot)$, 
and from (\ref{eq90.3.1}), (\ref{eq6.88.1}) and (\ref{eq6.88.2}), that
%
%
\begin{equation}
\label{eq21.7.3}
|x'|_{L} \le\biggl(\sup_{r}\frac{w_{r}}{\nu_{r}} \biggr)
\frac{2be^{b}(1+aN)|x'|}{\Gamma(1/N^{4})} \le C_{19}N^{2}|x'|
\end{equation}
for all $x' \in S$ and appropriate $C_{19}$. So application of (\ref
{eq39.2.1}), together with (\ref{eq21.7.3}) for $x'=X^{A}(t)$, implies
that
\begin{eqnarray*}
E_{x}\bigl[\bigl(|X(t)|_{L}-|x|_{L}\bigr)^{2}\bigr] &\le&
E_{x}[|X^{A}(t)|^{2}_{L}]
\le C_{19}^{2} N^{4}E_{x}[|X^{A}(t)|^{2}] \\
&\le& C_{19}^{2} N^{4}E_{x}[|A(t)|^{2}].
\end{eqnarray*}
Together with (\ref{eq21.6.2}), this implies
%
%
\begin{equation}
\label{eq21.7.5}
\sqrt{E_{x} \bigl[\bigl(|X(t)|_{L} - |x|_{L}\bigr)^{2}
\bigr]} \le C_{20}N^{2}t
\end{equation}
for appropriate $C_{20}$ and large $t$.

Substitution of (\ref{eq21.7.2}) and (\ref{eq21.7.5}) into
(\ref{eq21.7.1}) implies that for large enough $t$, (\ref{eq21.1.2})
holds, as desired.
\end{pf*}

\section{Upper bounds on $|X(N^{3})|_{r,s}$ for $s > N_{H_{r}}$}
\label{sec7}

In Section \ref{sec6}, we obtained upper bounds on
$E_{x}[|X(N^{3})|_{L} - |x|_{L}; \mathcal{A}(N^{3})^{c}]$; we still
need to analyze the behavior of $|X(N^{3})|_{L}-|x|_{L}$ on
$\mathcal{A}(N^{3})$. For this, we analyze $|X(N^{3})|_{r,s}$ for
several cases that depend on whether or not $|x|>N^{6}$ and
$s>N_{H_{r}}$.

In this section, we consider the case where $|x|>N^{6}$ and
$s>N_{H_{r}}$, which
is the simplest case. 
The main result here is the following proposition. Recall that $|x|_2$
is defined in (\ref{eq5.3.11}).
\begin{prop}
\label{prop40.2.3} For given $\varepsilon_{3}>0$, large enough $N$,
and $|x|>N^{6}$ and $|x|_{2}/|x|\le1/2$,
%
%
\begin{eqnarray}
\label{eq40.2.4}
&&E_{x} \Bigl[{\sup_{r, s> N_{H_{r}}}}{|X(N^{3})|}_{r,
s}-|x|_{L};
G \Bigr] \nonumber\\[-8pt]\\[-8pt]
&&\qquad\le C_{3}(|x|_{K}/|x|)N^{3} + \varepsilon_{3}
N^{2} \biggl(\frac{1}{2}-P(G) \biggr)\nonumber
\end{eqnarray}
for all measurable sets $G$, with $C_{3}> 0$ not depending on $N$, $G$
or $x$.
\end{prop}

In the proof of Proposition \ref{prop48.3.2}, we will employ Proposition
\ref{prop40.2.3} by setting
%
%
\begin{equation}
\label{eq7.2.1}
G =\mathcal{A}(N^{3})\cap\Bigl\{\omega\dvtx|X(N^{3})|_{L} = {\sup_{r, s
> N_{H_{r}}}} |X(N^{3})|_{r, s} \Bigr\}.
\end{equation}
%

Much of the work needed to demonstrate Proposition \ref{prop40.2.3} is
done in the following proposition. We recall that $i_{r}(s) = s +
\Delta_{r}$, where $\Delta_{r} = \Delta_{r}(N^{3})$.
%
\begin{prop}
\label{prop39.3.2}
For given $\varepsilon> 0$, large enough $N$ and all $x$,
%
%
\begin{equation}
\label{eq39.3.3}
E_{x} \Bigl[\sup_{r, s > N_{H_{r}}} \bigl\{|X(N^{3})|_{r,s} -
|x|_{r,i_{r}(s)} \bigr\}; G \Bigr] \le C_{21}\varepsilon N^{2}
\end{equation}
for all measurable sets $G$, with $C_{21}$ not depending on
$\varepsilon$, $N$, $G$ or $x$.
\end{prop}
\begin{pf}
We will instead show that
%
%
\begin{equation}
\label{eq39.3.1}
E_{x} \Bigl[{\sup_{r, s > N_{H_{r}}}}|X^{A}(N^{3})|_{r,s} \Bigr]
\le C_{21}\varepsilon N^{2}.
\end{equation}
Inequality (\ref{eq39.3.3}) follows immediately from this and inequality
(\ref{eq39.1.7}) since\break $|X^{A}(N^{3})|_{r,s} \ge0$.

To show (\ref{eq39.3.1}), we first note that for all $r$ and $s$,
%
%
\begin{equation}
\label{eq39.2.3}
|X^{A}(N^{3})|_{r,s} \le C_{22}\kappa_{N,r}
N_{H_{r}}Z^{A,*}_{r}(N^{3}, s)
\end{equation}
for appropriate $C_{22}$, where $Z^{A,*}_{r}(N^{3}, s)\stackrel
{\mathrm{def}}{=} (Z^{A})^{*}_{r}(N^{3}, s)$. The inequality uses
(\ref{eq5.3.1}) and (\ref{eq6.88.2}).
%
%
On $s>N_{H_{r}}$, the RHS of (\ref{eq39.2.3}) is at most
%
%
\begin{eqnarray}
\label{eq39.2.5}\qquad\quad
C_{22}b \kappa_{N,r} N_{H_{r}} Z^{A,*}_{r}(N^{3},
(s,\infty)) \le C_{22}b \kappa_{N,r} \int^{\infty
}_{N_{H_{r}}}N_{r}(s')Z^{A,*}_{r}(N^{3}, s')\,ds'
\end{eqnarray}
%
%
%
on account of (\ref{eq6.88.1}) and $N_{r}(s)\ge s$.

On the other hand, by (\ref{eq5.3.2}), the RHS of (\ref{eq39.2.5}) is
at most
%
%
\begin{equation}
\label{eq39.2.7}
(C_{22}b/M_{1})|X^{A}(N^{3})|_{R} \le C_{21}I_{R}(N^{3}),
\end{equation}
where $C_{21}\stackrel{\mathrm{def}}{=}C_{22}b/M_{1}$ and $I_{R}(\cdot
)$ is as in Section \ref{sec5}. Putting
(\ref{eq39.2.3})--(\ref{eq39.2.7}) together, it follows that, for large
$N$,
%
%
\begin{equation}
\label{eq39.2.8}
{\sup_{r, s > N_{H_{r}}}}|X^{A}(N^{3})|_{r,s} \le C_{21}
I_{R}(N^{3}).
\end{equation}
Also, by Proposition \ref{prop10.7.4}, 
for given $\varepsilon$, one has that for large enough $N$,
%
%
\begin{equation}
\label{eq39.2.9}
E_{x} [I_{R}(N^{3}) ] \le\varepsilon N^{2} \qquad\mbox{for all } x.
\end{equation}
Taking expectations in (\ref{eq39.2.8}) and applying (\ref{eq39.2.9}) implies
(\ref{eq39.3.1}).
\end{pf}

In order to demonstrate Proposition \ref{prop40.2.3}, we need Lemma
\ref{lem40.1.1}, which bounds $|x|_{K}$ from below in terms of $|x|$ when
$({\sup_{r, s \ge N_{H_{r}}}}|x|_{r,s})/|x|_{L}$ is not small. For
the lemma, we require the inequality
%
%
\begin{equation}
\label{eq40.1.5}
z^{*}_{r}((0, N_{H_{r}}]) \le C_{23}|x|_{L} \qquad\mbox{for } r \in
\mathcal{R},
\end{equation}
for appropriate $C_{23}$. This is a weaker version of (\ref{eq17.1.3}),
which we prove in Lem\-ma~\ref{lem17.1.6}. [Equation (\ref{eq40.1.5})
does not require any additional assumptions on $a$ or $b$, unlike
(\ref{eq17.1.3}).]

If one supposes that $|x|_{2} \le|x|/2$, it then follows easily by summing
(\ref{eq40.1.5}) over $r$ that
%
%
\begin{equation}
\label{eq7.5.1}
|x|_{L} \ge C_{24}|x|
\end{equation}
for $C_{24}= 1/2 C_{23}|\mathcal{R}|$. This inequality will be used in
Proposition \ref{prop40.2.3} and will also be used in Sections
\ref{sec8} and \ref{sec9}.
\begin{lem}
\label{lem40.1.1}
Suppose that, for some $r_{0}$ and $s_{0} \ge N_{H_{r_{0}}}$,
%
%
\begin{equation}
\label{eq40.1.2}
|x|_{r_{0},s_{0}} \ge|x|_{L}/2.
\end{equation}
Then, for appropriate $\varepsilon_{6}> 0$ not depending on $N$,
%
%
\begin{equation}
\label{eq40.1.3}
|x|_{K} \ge\varepsilon_{6}|x|/N.
\end{equation}
\end{lem}
\begin{pf}
Applying (\ref{eq40.1.5}), and then substituting (\ref{eq40.1.2}) into
(\ref{eq5.3.1}), one obtains for given $r$ that
%
%
\begin{eqnarray}
\label{eq40.1.4}
z^{*}_{r}((0, N_{H_{r}}]) &\le& C_{25}Nz^{*}_{r_{0}}(s_{0})/\Gamma
\bigl(\bar{H}^{*}_{r_{0}}(N_{H_{r_{0}}}+1)\bigr) \nonumber\\[-8pt]\\[-8pt]
&\le&
C_{25}be^{b}N\kappa_{N,r_{0}}z^{*}_{r_{0}}((N_{H_{r_{0}}},\infty))\nonumber
\end{eqnarray}
for appropriate $C_{25}> 0$, where the second inequality employs the
assumption $s_{0}\ge N_{H_{r_{0}}}$, together with (\ref{eq6.88.1}) and
(\ref{eq6.88.2}).
Addition of $z^{*}_{r}((N_{H_{r}},\infty))$ to both sides of (\ref{eq40.1.4})
gives
\begin{eqnarray*}
z^{*}_{r}(\mathbb{R}^{+}) &\le& z^{*}_{r}((N_{H_{r}},\infty)) + C_{25}
be^{b} N
\kappa_{N,r_{0}} z^{*}_{r_{0}}((N_{H_{r_{0}}},\infty)) \\
&\le& (1+C_{25}be^{b})N \sum_{r'} \kappa_{N,r'}z^{*}_{r'}
((N_{H_{r'}},\infty)).
\end{eqnarray*}
Summing over $r$ then implies
%
%
\begin{equation}
\label{eq40.2.2}
|x| \le\varepsilon_{6}^{-1}N|x|_{K}
\end{equation}
with $\varepsilon_{6}= [|\mathcal{R}|(1+C_{25}be^{b}) ]^{-1}$.
\end{pf}

We now apply Proposition \ref{prop39.3.2}, together with Lemma
\ref{lem40.1.1} and (\ref{eq7.5.1}), to demonstrate Proposition
\ref{prop40.2.3}.
\begin{pf*}{Proof of Proposition \protect\ref{prop40.2.3}}
Suppose first that $|x|_{r_{0},s_{0}}>|x|_{L}/2$ for some $r_{0}$ and
$s_{0} > N_{H_{r_{0}}}$. Choosing $\varepsilon> 0$ and $C_{21}$ as in
Proposition \ref{prop39.3.2}, with $\varepsilon$ small enough so
$\varepsilon< \varepsilon_{3} /C_{21}$ for given $\varepsilon_{3}> 0$,
it follows from the proposition and Lemma \ref{lem40.1.1} that for
large $N$ and any $G$, the LHS of (\ref{eq40.2.4}) is at most
%
%
\begin{eqnarray}
\label{eq40.3.1}
C_{21}\varepsilon N^{2} &\le& \varepsilon_{3}N^{2}
\le2\varepsilon _{3}\varepsilon_{6} ^{-1}(|x|_{K}/|x|)N^{3} -
\varepsilon_{3}N^{2} \nonumber\\[-8pt]\\[-8pt]
&\le& C_{3}(|x|_{K}/|x|)N^{3} - \varepsilon_{3}N^{2},\nonumber
\end{eqnarray}
if $C_{3}$ is chosen to be at least
$2\varepsilon_{3}\varepsilon_{6}^{-1}$, where $\varepsilon_{6} $ is as
in the lemma. This is at most the RHS of (\ref{eq40.2.4}).

Suppose, on the other hand, that $|x|_{r,s} \le|x|_{L}/2$ for all
$s>N_{H_{r}}$ and $r$. Under $|x|>N^{6}$ and $|x|_{2} \le|x|/2$, it
follows from (\ref{eq7.5.1}) that $|x|_{L} \ge C_{24}N^{6}$. Hence,
%
%
\begin{equation}
\label{eq40.3.2}
{\sup_{r, s>N_{H_{r}}}}|x|_{r,s} - |x|_{L} \le-\frac{1}{2}
C_{24}N^{6}.
\end{equation}
Since $i_{r}(s) \ge s > N_{H_{r}}$, it follows from Proposition \ref
{prop39.3.2}
and (\ref{eq40.3.2}) that the LHS of (\ref{eq40.2.4}) is at most
%
%
\begin{eqnarray}
\label{eq40.3.3}
C_{21}\varepsilon N^{2} + {\sup_{r, s>N_{H_{r}}}}|x|_{r,s} -
|x|_{L} &\le& C_{21}\varepsilon
N^{2}-\frac{1}{2}C_{24}N^{6}P(G)\nonumber\\[-8pt]\\[-8pt]
%
&\le& \varepsilon_{3}N^{2} \biggl(\frac{1}{2} - P(G) \biggr)\nonumber
\end{eqnarray}
for large $N$, if we choose $\varepsilon\le\varepsilon_{3}/2C_{21}$.
This is at most the RHS of (\ref{eq40.2.4}), which completes the proof.
\end{pf*}

\section{Pathwise upper bounds on $|X(N^{3})|_{r,s}$ for $s \le
N_{H_{r}}$ and $\Delta_{r}>1/b^{3}$}\label{sec8}

In the previous section, we analyzed the behavior of $|X(N^{3})|_{r,s}
- |x|_{L}$ for $s>N_{H_{r}}$. When $s \le N_{H_{r}}$, we analyze the
cases where $\Delta_{r} \le1/b^{3}$ and $\Delta_{r}>1/b^{3}$
separately. The latter case is quicker and we do it in this section,
postponing the case $\Delta_{r} \le1/b^{3}$ until Section \ref{sec9}.
For both cases, we will require certain pathwise upper bounds on
$|X^{A}(N^{3})|_{r,s}$ that hold on $\mathcal{A}_{1}(N^{3})$, which are
given in Proposition \ref{prop24.4.1}. We begin the section with these
bounds.

\subsection*{Upper bounds on $|X^{A}(N^{3})|_{r,s}$ on $\mathcal{A}_{1}
(N^{3})$}

In order to derive bounds on $|X^{A}(N^{3})|_{r,s}$, we first
require bounds on $Z^{A,*}_{r}(\cdot,\cdot)$ that measure how quickly
documents with the corresponding service times enter a route $r$ up to
a given time. In Lemma \ref{lem24.2.5}, we provide uniform bounds on
$Z^{A,*}_{r}(t,s)$ for $t \in [0,N^{3}]$ and
$\omega\in\mathcal{A}_{1}(N^{3})$. As in previous sections,
$S^{1}_{r}(k)$, $k=1, \ldots, A_{r}(t)$, denotes the positions of the
arrivals of documents up to time $t$. We also denote here by
$S^{2}_{r}(t,k)$ the amount of service such a document has received by
time $t$; $S^{1}_{r}(k) - S^{2}_{r}(t,k)$ is therefore the residual
service time of the $k$th document at time $t$.
\begin{lem}
\label{lem24.2.5}
Suppose $\omega\in\mathcal{A}_{1}(N^{3})$ for some
$N$. Then, for all $r$, $s \in[0, N+1]$ and $t \in[0, N^{3}]$,
%
%
\begin{equation}
\label{eq24.2.6}
Z^{A,*}_{r}(t,s)
\le4b\nu_{r}\bigl(\bar{H}^{*}_{r}(s)N^{3} \vee N^{3\eta}\bigr).
\end{equation}
If instead $s>N+1$, then
%
%
\begin{equation}
\label{eq24.3.1}
Z^{A,*}_{r}(t,s) \le2b\nu_{r}\bigl(\bar{H}^{*}_{r}(N+1)N^{3} \vee
N^{3\eta}\bigr).
\end{equation}
\end{lem}
\begin{pf}
For all $r$, $s \in[0,N+1]$ and $t \in[0,N^{3}]$,
%
%
\begin{eqnarray}
\label{eq24.3.2}\qquad
Z^{A,*}_{r}(t,s) &=& \sum^{A_{r}(t)}_{k=1} \phi\bigl(s - S^{1}_{r}(k) +
S^{2}_{r}(t,k)\bigr)\nonumber\\[-8pt]\\[-8pt]
&\le& \sum^{A_{r}(N^{3})}_{k=1} \sup_{s^{\prime}\in[0,\infty)}\phi
\bigl(s - S^{1}_{r}(k)+s'\bigr)\le b\sum^{A_{r}(N^{3})}_{k=1} \bar{\Phi}\bigl(s
- S^{1}_{r}(k)\bigr)\nonumber
\end{eqnarray}
with the latter inequality employing $\phi(s) \le b\bar{\Phi}(s)$
and the
monotonicity of $\bar{\Phi}(\cdot)$. Letting $j_{0}$ denote the
largest $j$ with
$v_{j} \le s$, the last term in (\ref{eq24.3.2}) is at most
%
%
\begin{equation}
\label{eq24.3.3}
b\sum^{A_{r}(N^{3})}_{k=1} \bar{\Phi}\bigl(v_{j_{0}} - S^{1}_{r}(k)\bigr)
\le2b\nu_{r}
\bigl(\bar{H}^{*}_{r}(v_{j_{0}})N^{3} \vee N^{3\eta}\bigr)
\end{equation}
on $\mathcal{A}_{1}(N^{3})$. The inequality in (\ref{eq24.3.1})
follows from
this, with $j_{0}=J$. The inequality in (\ref{eq24.2.6}) follows by applying
(\ref{eq20.1.7}) to the RHS of (\ref{eq24.3.3}).
\end{pf}

We now derive uniform upper bounds on $|X^{A}(t)|_{r,s}$ for $t \in[0,N^{3}]$
and $\omega\in\mathcal{A}_{1}(N^{3})$. In applications, we will be primarily
interested in the behavior at $t = N^{3}$.
\begin{prop}
\label{prop24.4.1}
Suppose $\omega\in\mathcal{A}_{1}(N^{3})$ for some
$N$. Then, for all $r$ and $s$,
%
%
\begin{equation}
\label{eq24.4.7}
|X^{A}(t)|_{r,s} \le C_{26}N^{3} \qquad\mbox{for } t \in[0, N^{3}]
\mbox{ and all } x,
\end{equation}
for appropriate $C_{26}$ not depending on $x, N, \omega, r$ or $s$. In
particular,
%
%
\begin{equation}
\label{eq24.4.2}
|X(t)|_{L} - |x|_{L} \le C_{26}N^{3} \qquad\mbox{for } t \in[0, N^{3}]
\mbox{ and all } x.
\end{equation}
\end{prop}
\begin{pf}
By (\ref{eq39.2.1}), 
%
%
\begin{equation}
\label{eq24.4.4}
|X(t)|_{L} - |x|_{L} \le|X^{A}(t)|_{L} \qquad\mbox{for all } t,
\end{equation}
and so (\ref{eq24.4.2}) follows immediately from (\ref{eq24.4.7}).

We now investigate $|X^{A}(t)|_{r,s}$. From (\ref{eq5.3.1}) and Lemma
\ref{lem24.2.5}, it follows that, for $t \in[0, N^{3}]$,
%
%
\begin{eqnarray}
\label{eq24.4.5}
|X^{A}(t)|_{r,s} &=& \frac{w_{r}(1+as_{N})Z^{A,*}_{r}(t,s)}{\nu
_{r}\Gamma
(\bar{H}^{*}_{r}(s_{N}))} \nonumber\\
&\le& 4bw_{r}N^{3}(1+as_{N})\bar{H}^{*}_{r}(s_{N})/\Gamma
(\bar{H}^{*}_{r}(s_{N})) \\
&&{} + 4bw_{r}N^{3\eta}(1+as_{N})/\Gamma(\bar{H}^{*}_{r}(s_{N})).
\nonumber
\end{eqnarray}
We proceed to analyze the two terms on the RHS of (\ref{eq24.4.5}).

It follows from the definition of $\Gamma(\cdot)$ in (\ref{eq5.1.3}) that,
for all $s$,
%
%
\begin{equation}
\label{eq24.4.6}
(1+as_{N})\bar{H}^{*}_{r}(s_{N})/\Gamma(\bar{H}^{*}_{r}(s_{N})) \le
(1+as_{N}) (\bar{H}^{*}_{r}(s_{N}))^{1-\gamma}/a C_{2}.
\end{equation}
Since by assumption, $\bar{H}^{*}_{r}(\cdot)$ has more than two moments
and $\gamma\le1/2$, the RHS of (\ref{eq24.4.6}) goes to 0 as $s_{N}\to
\infty$. Hence, it is bounded for all $s_{N}$, which implies that the
first term on the RHS of (\ref{eq24.4.5}) is bounded above by
$C_{27}N^{3}$, for some $C_{27}$ not depending on $t, r$ or $s$.

On the other hand, for all $s$,
%
%
\begin{eqnarray}
\label{eq24.5.1}
(1+as_{N})/\Gamma(\bar{H}^{*}_{r}(s_{N}))
&\le& \bigl(1+a(N+1)\bigr)\bar {H}^{*}_{r}(N_{H_{r}}
+1)^{-\gamma}/a C_{2}\nonumber\\[-8pt]\\[-8pt]
&\le& \bigl(1+a(N+1)\bigr)(e^{b}N^{4})^{\gamma}/C_{2}a.\nonumber
\end{eqnarray}
Since $\gamma\le1/4$, $\eta\le1/3$ and $aN \ge1$, the latter term on
the RHS of (\ref{eq24.4.5}) is bounded above by $C_{28}N^{2}$, for some
$C_{28}$ not depending on $t$, $r$ or $s$.

The above bounds for the two terms on the RHS of (\ref{eq24.4.5}) sum
to $(C_{27} + C_{28})N^{3}$. Setting $C_{26}= C_{27}+ C_{28}$, this
implies (\ref {eq24.4.7}).
\end{pf}

\subsection*{Upper bounds on $|X(N^{3})|_{r,s}$ for $s \le N_{H_{r}}$ and
$\Delta_{r}> 1/b^{3}$}
Proposition \ref{prop41.2.1} gives an upper
bound on $|X(N^{3})|_{r,s}-|x|_{L}$ when $s \le N_{H_{r}}$ and
$\Delta_{r}>1/b^{3}$. The proof, which employs Proposition
\ref{prop24.4.1}, is quick.
\begin{prop}
\label{prop41.2.1}
Suppose that $|x|>N^{6}$, with $|x|_{2}/|x| \le1/2$.
Then, for large enough $N$,
%
%
\begin{equation}
\label{eq41.2.2}\quad
{\sup_{\Delta_{r}>1/b^{3}}\, \sup_{s\le N_{H_{r}}}}
|X(N^{3})|_{r,s} - |x|_{L} \le-N^{4} \qquad\mbox{for all }\omega\in
\mathcal{A}_{1}(N^{3}),
\end{equation}
where $N$ does not depend on $x$ or $\omega$.
\end{prop}
\begin{pf}
For each $r$ and $s$,
%
%
\begin{eqnarray}
\label{eq41.2.4}
|X(N^{3})|_{r,s}-|x|_{L} &=& |X^{A}(N^{3})|_{r,s}-
\bigl(|x|_{L}-|\tilde{X}(N^{3})
|_{r,s} \bigr) \nonumber\\[-8pt]\\[-8pt]
&\le& C_{26}N^{3} - \bigl(|x|_{L} - |\tilde{X}(N^{3})|_{r,s}
\bigr)\nonumber
\end{eqnarray}
with the last line following from Proposition \ref{prop24.4.1}. We
consider two cases, depending on whether $|x|_{r,i_{r}(s)}>|x|_{L}/2$
for given $r$ and $s$.

Suppose first that $|x|_{r,i_{r}(s)}>|x|_{L}/2$, with $s \le N_{H_{r}}$
and $|x|_{2}\le|x|/2$. One has
%
%
\begin{equation}
\label{eq41.2.6}
|x|_{r,i_{r}(s)} - |\tilde{X}(N^{3})|_{r,s}
\ge\frac{w_{r}}{\nu _{r}} \cdot
\frac{a}{b^{3}}\cdot\frac{z^{*}_{r}(i_{r}(s))}{\Gamma(\bar
{H}^{*}_{r}(i_{r} (s)_{N}) )}.
\end{equation}
To see this, one applies (\ref{eq39.1.5}) to the definition of
$|x|_{r,s}$ in (\ref{eq5.3.1}), noting that since $s \le N_{H_{r}}$,
%
%
\begin{equation}
\label{eq8.9.1}
i_{r}(s)_{N} - s_{N} = i_{r}(s) \wedge(N_{H_{r}}+1) - s
\ge\Delta _{r} \wedge1
> 1/b^{3},
\end{equation}
and that $\Gamma(\bar{H}^{*}_{r}(i_{r}(s)_{N})) \le\Gamma(\bar
{H}^{*}_{r}(s))$. On account of (\ref{eq5.3.1}) and
$|x|_{r,i_{r}(s)}>|x|_{L}/2$,
one obtains, from the RHS of (\ref{eq41.2.6}),
\begin{eqnarray*}
&& \frac{a}{b^{3}} \cdot\biggl(\frac{w_{r}z^{*}_{r}(i_{r}(s))}{\nu
_{r}\Gamma
(\bar{H}^{*}_{r}(i_{r}(s)_{N}))} \Big/|x|_{r,i_{r}(s)} \biggr) \cdot
\frac{|x|_{r,i_{r}(s)}}{|x|_{L}} \cdot|x|_{L} \\
&&\qquad \ge\frac{a}{b^{3}} \cdot\bigl(1+ai_{r}(s)_{N}\bigr)^{-1} \cdot\frac
{1}{2} \cdot
|x|_{L}. 
\end{eqnarray*}
Because of $|x|_{2}\le|x|/2$, (\ref{eq7.5.1}), $i_{r}(s)_{N}\le N$,
$|x|>N^{6}$ and $aN \ge1$, this is at most $C_{29}N^{5}$, where
$C_{29}> 0$ does not depend on $N$, $x$ or $\omega$. It follows from
(\ref{eq41.2.6}) and the succeeding inequalities that
%
%
\begin{equation}
\label{eq41.3.2}
|x|_{L}-|\tilde{X}(N^{3})|_{r, s} \ge|x|_{r,i_{r}(s)} - |\tilde
{X}(N^{3})|_{r,s}
\ge C_{29}N^{5}.
\end{equation}
Together with (\ref{eq41.2.4}), this gives the RHS of (\ref{eq41.2.2}).

Suppose, on the other hand, that $|x|_{r,i_{r}(s)} \le|x|_{L}/2$, with $|x|_{2}
\le|x|/2$. Then, by (\ref{eq7.5.1}) and (\ref{eq39.1.6}), the RHS of
(\ref{eq41.2.4}) is at most
%
%
\begin{eqnarray}
\label{eq41.2.5}
&&
C_{26}N^{3} - \tfrac{1}{2} |x|_{L} - \bigl(|x|_{r,i_{r}(s)} - |\tilde{X}
(N^{3})|_{r,s} \bigr) \nonumber\\[-8pt]\\[-8pt]
&&\qquad \le C_{26}N^{3} - \tfrac{1}{2} C_{24}N^{6} \le- N^{5}\nonumber
\end{eqnarray}
for large $N$. This implies (\ref{eq41.2.2}) for $|x|_{r,i_{r}(s)} \le
|x|_{L}/2$, and hence completes the proof.
\end{pf}

\section{Pathwise upper bounds on $|X(N^{3})|_{r,s}$ for $s \le
N_{H_{r}}$ and $\Delta_{r} \le1/b^{3}$}\label{sec9}

In Sections \ref{sec7} and \ref{sec8}, we analyzed the behavior of
$|X(N^{3})|_{r,s} - |x|_{L}$ for $s > N_{H_{r}}$, and for $s \le
N_{H_{r}}$ with $\Delta_{r} > 1/b^{3}$. There remains the case $s \le
N_{H_{r}}$ with $\Delta_{r} \le1/b^{3}$, which is the subject of this
section. This is, in essence, the ``main case'' one needs to show in
order to establish the stability of the network since the other cases
dealt with less sensitive behavior and did not employ the
subcriticality of the system that was given in (\ref{eq5.3.8}). The
same was also true for the computations of the $|\cdot|_{A}$ and
$|\cdot|_{R}$ norms in Sections \ref{sec4} and \ref{sec5}.

Section \ref{sec9} consists of three subsections. First, in Proposition
\ref{prop46.3.2}, we give lower bounds on the minimal service rates
$\lambda^{w} (\cdot)$ of documents in terms of the norm $|\cdot|_{L}$.
In the next subsection, we begin our analysis of $|X(N^{3})|_{r,s}$ for
$s\le N_{H_{r}}$ and $\Delta_{r}\le1/b^{3}$. We decompose
$|X(N^{3})|_{r,s} - |x|_{r,i_{r}(s)}$ into several parts that are
easier to analyze. In Proposition \ref {prop44.1.2}, we then obtain
upper bounds on the factor $Z^{*}_{r}(N^{3},s) - z^{*}_{r}(i_{r}(s))$
of one of the parts. In the third subsection, we do a detailed analysis
of the decomposition from the previous subsection, which also employs
the bounds on $\lambda^{w}(\cdot)$ from the first subsection. From
this, we obtain in Proposition \ref{prop48.1.1} the desired bound on
$|X(N^{3})|_{r,s} - |x|_{L}$. We note that, whereas in Section
\ref{sec8}, our results pertained to $\omega \in \mathcal{A}_{1}
(N^{3})$, starting from the second subsection here, we require
$\omega\in\mathcal{A}_{2}(t)$. Our final results on
$|X(N^{3})|_{r,s}-|x|_{L}$, for $s \le N_{H_{r}}$, will therefore be
valid on $\mathcal{A}(N^{3}) = \mathcal{A}_{1}(N^{3})
\cap\mathcal{A}_{2}(N^{3})$.

\subsection*{Lower bounds on $\lambda^{w}(\cdot)$}
In order to demonstrate the stability of the network, its
subcriticality needs to be employed at some point. With this in mind,
we choose $\varepsilon _{7}\in(0,1]$ small enough so that
%
%
\begin{equation}
\label{eq17.1.5}
(1+\varepsilon_{7})^{2} \sum_{r
\in\mathcal{R}}A_{l,r}\rho_{r} \le c_{l} \qquad\mbox{for all } l,
\end{equation}
which is possible because of (\ref{eq5.3.8}). We henceforth assume
$\varepsilon_{5}\le \varepsilon_{7}/4$, where $\varepsilon_{5}$ was
employed in (\ref {eq20.2.2}) in the definition of
$\mathcal{A}_{2}(\cdot)$.

The main results in this subsection are Propositions \ref{prop17.2.1}
and \ref{prop46.3.2}. Proposition \ref{prop17.2.1} gives a lower bound
on $\lambda^{w}(t)$ in terms of $|X(t)|_{S}$; Proposition
\ref{prop46.3.2}, under additional assumptions, gives the bound in
terms of $|x|_{L}$.
\begin{prop}
\label{prop17.2.1}
Assume (\ref{eq17.1.5}) holds for some
$\varepsilon_{7}> 0$. Then, for large enough $b$ and small enough $a$,
%
%
\begin{equation}
\label{eq17.2.1}
\lambda^{w}(t) \ge(1+\varepsilon_{7})/|X(t)|_{S}
\end{equation}
for almost all $t$.
\end{prop}

In this and the previous subsection, we need to employ certain
properties of $\Gamma(\bar{H}^{*}_{r}(\cdot))$, which appears in the
denominator in (\ref{eq5.3.1}). In Lemma \ref{lem5.2.6}, we state two
such properties; the first is employed for Lemma \ref{lem46.8.2} and
the second is employed for Lemma \ref{lem17.1.6}. Recall that $m_r$ is
the mean of $H_{r}(\cdot)$.
\begin{lem}
\label{lem5.2.6}
For $\Gamma(\cdot)$ as defined in (\ref{eq5.1.3}),
%
%
\begin{equation}
\label{eq5.2.2}
\Gamma'(\bar{H}^{*}_{r}(s)) \ge1 + as \qquad\mbox{for all $r$ and $s$}.
\end{equation}
Moreover, for large enough $b$ and small enough $a$,
%
%
\begin{equation}
\label{eq5.2.5}
\int^{\infty}_{0}\frac{\Gamma(\bar{H}^{*}_{r}(s))}{1+as}\,ds \le
(1+\varepsilon_{7})
m_{r}
\end{equation}
for $\varepsilon_{7}> 0$ satisfying (\ref{eq17.1.5}).
\end{lem}
\begin{pf}
By (\ref{eq5.1.3}) and then (\ref{eq5.2.1}), one has, for all $r$ and
$s$,
%
%
\begin{eqnarray}
\label{eq5.2.7}
\Gamma'(\bar{H}^{*}_{r}(s)) &=& 1+C_{2}\gamma a(\bar
{H}^{*}_{r}(s))^{\gamma- 1}
\nonumber\\[-8pt]\\[-8pt]
&\ge& 1+C_{2}C_{1}^{\gamma- 1}\gamma a(1+s)^{(1-\gamma)(2+\delta_{1})}
\ge1 + as,\nonumber
\end{eqnarray}
where the last inequality uses $\gamma\le1/2$ and $C_{2}\ge
{C_{1}}^{(1-\gamma)/\gamma}$. This implies (\ref{eq5.2.2}).

For (\ref{eq5.2.5}), we note from (\ref{eq5.1.3}) and (\ref
{eq5.2.1}) that
%
%
\begin{equation}
\label{eq5.2.3}\hspace*{28pt}
\int^{\infty}_{0}\frac{\Gamma(\bar{H}^{*}_{r}(s))}{1+as}\,ds \le
\int^{\infty}_{0}\bar{H}^{*}_{r}(s)\,ds + C_{1}^{2} a \int^{\infty}_{0}
(1+s)^{-2\gamma}(1+as)^{-1}\,ds.
\end{equation}
The constant $b$ can be chosen large enough so the first term on the
RHS of (\ref{eq5.2.3}) is at most $(1+\varepsilon_{7}/2)m_{r}$. Also,
by choosing $a>0$ small enough, since the second term can be chosen as
close to 0 as desired, by monotone convergence,
%
%
\begin{equation}
\label{eq5.2.4}\hspace*{28pt}
a \int^{\infty}_{0}(1+s)^{-2\gamma}(1+as)^{-1}\,ds = \int^{\infty}_{0}
(1+s)^{-(1+2\gamma)}\frac{1+s}{1/a+s}\,ds \to0
\end{equation}
as $a \searrow0$. So, for large enough $b$ and small enough $a$,
(\ref{eq5.2.5}) holds.
\end{pf}

By employing (\ref{eq5.2.5}), we obtain upper bounds for $z^{*}_{r}
((0,N_{H_{r}}])$ and $z^{*}_{r}(\mathbb{R}^{+})$ in terms of $|x|_{L}$ and
$|x|_{S}$. Inequality (\ref{eq17.1.4}) will be crucial for Proposition
\ref{prop17.2.1}.
\begin{lem}
\label{lem17.1.6}
For large enough $b$ and small enough $a$,
%
%
\begin{equation}
\label{eq17.1.3}
z^{*}_{r}((0, N_{H_{r}}]) \le(1 + \varepsilon_{7})w^{-1}_{r}\rho_{r}|x|_{L}
\end{equation}
and
%
%
\begin{equation}
\label{eq17.1.4}
z^{*}_{r}(\mathbb{R}^{+}) \le(1 + \varepsilon_{7})w^{-1}_{r}\rho_{r}|x|_{S}
\end{equation}
for all $N$ and $r$, where $\varepsilon_{7}> 0$ is as in (\ref{eq17.1.5}).
\end{lem}
\begin{pf}
We note that by (\ref{eq5.3.1}),
%
%
\begin{equation}
\label{eq17.1.1}\quad
z^{*}_{r}((0,N_{H_{r}}]) = \int^{N_{H_{r}}}_{0}z^{*}_{r}(s)\,ds \le w^{-1}_{r}
\nu_{r}|x|_{L} \int^{N_{H_{r}}}_{0} \frac{\Gamma(\bar{H}^{*}_{r}
(s))}{1+as}\,ds.
\end{equation}
By (\ref{eq5.2.5}), for large enough $b$ and small enough $a$, the last
term in (\ref{eq17.1.1}) is at most
%
%
\begin{equation}
\label{eq17.1.2}
(1+\varepsilon_{7})w^{-1}_{r}\nu_{r}m_{r}|x|_{L} = (1+\varepsilon
_{7})w^{-1}_{r}\rho_{r}|x|_{L}
\end{equation}
for all $N$ and $r$, which implies (\ref{eq17.1.3}). It follows from
(\ref{eq17.1.3}) and the definition of $|\cdot|_{S}$ in (\ref
{eq5.3.12}) that
\begin{eqnarray*}
z^{*}_{r}(\mathbb{R}^{+}) &\le& (1+\varepsilon_{7})w^{-1}_{r}\rho_{r} \biggl[
|x|_{L} +
\frac{w_{r}}{\rho_{r}}z^{*}_{r}((N_{H_{r}}, \infty)) \biggr] \\
&\le& (1+\varepsilon_{7})w^{-1}_{r}\rho_{r}|x|_{S},
\end{eqnarray*}
which implies (\ref{eq17.1.4}).
\end{pf}

A weaker version of the bound (\ref{eq17.1.3}) was used in
(\ref{eq40.1.5}), where the RHS of (\ref{eq17.1.3}) was replaced by
$C_{23}|x|_L$, and no additional assumptions on $b$ and $a$ were
required. This follows by noting that the second term on the RHS of
(\ref{eq5.2.3}) does not depend on $a$ (since $a\le1$).

We now demonstrate Proposition \ref{prop17.2.1}.
\begin{pf*}{Proof of Proposition \protect\ref{prop17.2.1}}
On account of (\ref{eq17.1.5}), a feasible protocol is given by
assigning service to each nonempty route $r$ at rate
$\Lambda_{r,F}\stackrel {\mathrm{def}}
{=}(1+\varepsilon_{7})^{2}\rho_{r}$. By (\ref{eq17.1.4}), the rate at
which each document is served is
%
%
\begin{equation}
\label{eq17.2.3}
\lambda_{r,F} = \frac{(1+\varepsilon_{7})^{2}\rho
_{r}}{Z_{r}(t,\mathbb
{R}^{+})} = \frac{
(1+\varepsilon_{7})^{2}\rho_{r}}{Z^{*}_{r}(t,\mathbb{R}^{+})} \ge
\frac
{(1+\varepsilon_{7})w_{r}}{
|X(t)|_{S}}
\end{equation}
at almost all times $t$. It follows from this and the definition of the weighted
max--min fair protocol that
\[
\lambda^{w}(t) = \min_{r \in\mathcal{R}'}
\frac{\lambda_{r}(t)}{w_{r}} \ge\min_{r \in\mathcal{R}'}
\frac{\lambda_{r,F}}{w_{r}} \ge\frac{(1+\varepsilon_{7})}{|X(t)|_{S}}
\qquad\mbox{for almost all } t,
\]
%
which implies (\ref{eq17.2.1}).
\end{pf*}

We apply Proposition \ref{prop17.2.1} to derive the following lower
bound of
$\lambda_{r}(t)$ on $[0,N^{3}]$. We note that, by (\ref{eq24.4.2}) of
Proposition \ref{prop24.4.1} and (\ref{eq7.5.1}), for $\omega\in
\mathcal{A}_{1}(N^{3})$, $|x|>N^{6}$ and $|x|_{2}\le|x|/2$,
%
%
\begin{equation}
\label{eq46.3.1}
|X(t)|_{L} \le|x|_{L} + C_{26}N^{3} \le(1+\varepsilon)|x|_{L}
\end{equation}
holds for given $\varepsilon> 0$ and large enough $N$. In the
proposition, we
will use
\[
\varepsilon_{8}\stackrel{\mathrm{def}}{=} \biggl[\frac{C_{24}}{8}
\biggl(\max_{r}\frac{w_{r}}{\rho_{r}} \biggr)^{-1}\varepsilon_{7}\biggr]
\wedge\frac{1}{2}.
\]
\begin{prop}
\label{prop46.3.2}
Suppose that (\ref{eq17.1.5}) holds for some $\varepsilon_{7}\in(0,1]$, and
that $|x|
> N^{6}$, with $|x|_{2} \le\varepsilon_{8}|x|$. Then, for large enough
$N$ and
$b$, and
small enough $a$,
%
%
\begin{equation}
\label{eq46.3.8}
\lambda^{w}(t) \ge(1+\varepsilon_{7}/2)/|x|_{L}
\end{equation}
for almost all $t \in[0,N^{3}]$ on $\omega\in\mathcal{A}_{1}(N^{3})$.
\end{prop}
\begin{pf}
It follows from Proposition \ref{prop17.2.1} that
%
%
\begin{equation}
\label{eq46.4.1}
\lambda^{w}(t) \ge(1+\varepsilon_{7})/|X(t)|_{S} \qquad\mbox{almost
everywhere},
\end{equation}
for large enough $b$ and small enough $a$. On the other hand, it
follows from
(\ref{eq5.3.12}), (\ref{eq46.3.1}), (\ref{eq39.1.9}) and (\ref{eq20.2.1})
that, since $|x|>N^{6}$ and $|x|_{2} \le\varepsilon_{8}|x|$,
%
%
\begin{eqnarray}
\label{eq46.4.2}
|X(t)|_{S} &\le& |X(t)|_{L} + \biggl(\max_{r} \frac{w_{r}}
{\rho_{r}} \biggr)Z^{*}_{r}(t,(N_{H_{r}}, \infty))
\nonumber\\[-8pt]\\[-8pt]
&\le& (1+\varepsilon)|x|_{L} + \biggl(\max_{r} \frac{w_{r}}
{\rho_{r}} \biggr) \Bigl[|x|_{2} + 2 \Bigl(\max_{r} \nu_{r}
\Bigr)N^{3} \Bigr]\nonumber
\end{eqnarray}
%
%
holds for given $\varepsilon> 0$ and large enough $N$, for all $\omega
\in \mathcal{A}_{1}(N^{3})$ and $t \in[0,N^{3}]$. Applying $|x|>N^{6}$,
$|x|_2 \le\varepsilon_{8}|x|$ and (\ref{eq7.5.1}) to the RHS of (\ref
{eq46.4.2}) implies that it is at most
\[
\biggl(1+\varepsilon+\frac{\varepsilon_{7}}{8} \biggr)|x|_{L} + 2 \biggl(\max_{r}
\frac{w_{r}}{\rho_{r}} \biggr) \Bigl(\max_{r} \nu_{r} \Bigr)|x|_{L}/C_{24}N^{2}.
\]
Consequently, for small enough
$\varepsilon> 0$,
%
%
\begin{equation}
\label{eq46.4.4}
|X(t)|_{S} \le(1+\varepsilon_{7}/4)|x|_{L} \qquad\mbox{for all } t \in[0,N^{3}].
\end{equation}
Together with (\ref{eq46.4.1}), this implies (\ref{eq46.3.8}).
\end{pf}

\subsection*{Decomposition of $|X(N^{3})|_{r,s}-|x|_{r,i_{r}(s)}$}
In this short subsection, we decompose
$|X(N^{3})|_{r,s}-|x|_{r,i_{r}(s)}$ into several parts, one of which
contains the factor $Z^{*}_{r}(N^{3},s) - z^{*}_{r}(i_{r}(s))$. In
Proposition \ref{prop44.1.2}, we then obtain upper bounds on this
factor. In this and the remaining subsection, the estimates need to be
more precise than in previous sections in order to make use of the
subcriticality of $X(\cdot)$.

The decomposition that was referred to above is given by
%
%
\begin{eqnarray}
\label{eq42.1.2}
&&|X(N^{3})|_{r,s} - |x|_{r,i_{r}(s)} \nonumber\\
&&\qquad= \frac
{w_{r}(1+as)(Z^{*}_{r}(N^{3},s) -
z^{*}_{r}(i_{r}(s)))}{\nu_{r}\Gamma(\sigma_{r})} \\
&&\qquad\quad{} - |x|_{r,i_{r}(s)}\frac{1+as}{1+ai_{r}(s)} \frac{\Gamma(\sigma
_{r}) - \Gamma
({\sigma_{r}'})}{\Gamma(\sigma_{r})} - \frac{aw_{r}\Delta_{r}z^{*}_{r}
(i_{r}(s))}{\nu_{r}\Gamma({\sigma_{r}'})}, \nonumber
\end{eqnarray}
and holds for $s \le N_{H_{r}}$ and $\Delta_{r} \le1/b^{3}$. It will
be employed in
Corollary \ref{col46.2.1}. Here and later on,
we abbreviate, setting $\sigma_{r} = \bar{H}^{*}_{r}(s)$ and ${\sigma
_{r}'} =
\bar{H}^{*}_{r}(i_{r}(s))$. [One can check that (\ref{eq42.1.2})
holds as given,
without employing either $s_{N}$ or $i_{r}(s)_{N}$, as in (\ref
{eq5.3.1}), since
$i_{r}(s) = s + \Delta_{r} \le N_{H_{r}}+1$, and hence $s_{N} = s$ and
$i_{r}(s)_{N} = i_{r}(s)$.]

To apply the bound (\ref{eq20.2.2}) on $\omega\in\mathcal
{A}_{2}(N^{3})$ and
derive an upper bound on $Z^{*}_{r}(N^{3},s) - z^{*}_{r}(i_{r}(s))$, we
need to
select a $v_{j}$ from among $v_{0}, \ldots, v_{J}$, as given by
(\ref{eq20.1.6}). For this, we denote by $v(s)$ the value $v_{j}$ with
%
%
\begin{equation}
\label{eq44.2.1}
v_{j} \in\bigl[i_{r}(s), i_{r}(s) + 1/b^{3}\bigr).
\end{equation}
Under $s \le N_{H_{r}}$ and $\Delta_{r} \le1/b^{3}$, such a $v(s)$ exists.
\begin{prop}
\label{prop44.1.2} Suppose $\omega\in\mathcal{A}_{2}(N^{3})$, for some
$N$ and $b$, with $b$ as in (\ref{eq5.1.4}). Then,
%
%
\begin{equation}
\label{eq44.1.3}\qquad
Z^{*}_{r}(N^{3},s) - z^{*}_{r}(i_{r}(s)) \le(1 +
\varepsilon_{5})(1 + 4/b^{2}) \nu_{r} [ h^{*}_{r}(v(s))N^{3} \vee
N^{3\eta} ]
\end{equation}
for all $r$ and $s$ with $\Delta_{r} \le1/b^{3}$ and $s \le
N_{H_{r}}$, where
$\varepsilon_{5}>0$ is as in (\ref{eq20.2.2}) and $v(s)$ is given by
(\ref{eq44.2.1}).
\end{prop}
\begin{pf}
By (\ref{eq39.1.5}), the LHS of (\ref{eq44.1.3}) is at most $Z^{A,*}_{r}
(N^{3},s)$. For $s \le N_{H_{r}}$, this equals
%
%
\begin{eqnarray}
\label{eq44.1.5}\qquad\quad
\sum^{A_{r}(N^{3})}_{k=1}\phi\bigl(s - S^{1}_{r}(k) + S^{2}_{r}(N^{3},k)\bigr)
&\le& e^{2/b^{2}} \sum^{A_{r}(N^{3})}_{k=1}\phi\bigl(v(s) -
S^{1}_{r}(k)\bigr) \nonumber\\[-8pt]\\[-8pt]
&\le&(1 + 4/b^{2}) \sum^{A_{r}(N^{3})}_{k=1}\phi\bigl(v(s) -
S^{1}_{r}(k)\bigr).\nonumber
\end{eqnarray}
To see (\ref{eq44.1.5}), we note that since $S^{2}_{r}(N^{3}, k) \le
\Delta_{r}
\le1/b^{3}$,
%
%
\begin{equation}
\label{eq44.2.5}\hspace*{28pt}
v_{j} - S^{1}_{r}(k) \in[s - S^{1}_{r}(k) +
S^{2}_{r}(N^{3},k), s - S^{1}_{r}(k) + S^{2}_{r}(N^{3},k) + 2/b^{3}].
\end{equation}
Together with the second half of (\ref{eq5.1.10}), this implies the first
inequality. The second inequality follows by expanding $e^{2/b^{2}}$. Since
$\omega\in\mathcal{A}_{2}(N^{3})$, the RHS of (\ref{eq44.1.3}) then
follows by
applying (\ref{eq20.2.2}).
\end{pf}

In the next subsection, we will also employ the following bound on $h^{*}_{r}
(s_{2}) - h^{*}_{r}(s_{1})$ for $s_{1} \le s_{2}$.
\begin{prop}
\label{prop46.5.1}
For any $r$, $s_{1} \le s_{2}$ and $b$,
%
%
\begin{equation}
\label{eq46.5.2}
h^{*}_{r}(s_{2}) - h^{*}_{r}(s_{1}) \le eb^{2}(s_{2} - s_{1})\bar
{H}^{*}_{r}(s_{1}).
\end{equation}
\end{prop}
\begin{pf}
Since $h^{*}_{r}(s) = \int^{\infty}_{0}\phi(s-s')\,dH_{r}(s')$ for each
$s$, the LHS of (\ref{eq46.5.2}) equals
%
%
\begin{equation}
\label{eq46.5.3}
\int^{\infty}_{0}\bigl(\phi(s_{2}-s') - \phi(s_{1}-s')\bigr)\,dH_{r}(s').
\end{equation}
By the first part of (\ref{eq5.1.10}) and the definition of $\phi
(\cdot)$,
$\phi'(s) \le b^{2}$ for all $s$ and $\phi(\cdot)$ is decreasing on $[1/b,
\infty)$. So, (\ref{eq46.5.3}) is at most
%
%
\begin{eqnarray}
\label{eq46.5.4}\hspace*{28pt}
\int^{\infty}_{0}b^{2}(s_{2}-s_{1})1\{s'>s_{1}-1/b\}\,dH(s') &\le& b^{2}
(s_{2}-s_{1})\bar{H}_{r}(s_{1}-1/b) \nonumber\\
&\le& b^{2}(s_{2}-s_{1})\bar{H}^{*}_{r}(s_{1}-1/b)\\
&\le& eb^{2}(s_{2}-s_{1})
\bar{H}^{*}_{r}(s_{1}).\nonumber
\end{eqnarray}
\upqed\end{pf}

\subsection*{Upper bounds on $|X(N^{3})|_{r,s}$}

In this subsection, we employ the previous two subsections to obtain
upper bounds on $|X(N^{3})|_{r,s}-|x|_{L}$ for $\omega\in\mathcal
{A}(N^{3})$, when $s \le N_{H_{r}}$ and $\Delta_{r} \le1/b^{3}$. Our
main result is the following proposition. As elsewhere in this paper,
we are assuming that $aN\ge1$.
\begin{prop}
\label{prop48.1.1} Suppose that (\ref{eq17.1.5}) holds for some
$\varepsilon_{7}\in[0,1]$ and that $|x| > N^{6}$, with $|x|_{2}
\le\varepsilon_{8}|x|$, where $\varepsilon_{8}$ is specified below
(\ref{eq46.3.1}). Then, for large enough $N$ and $b$, and small enough
$a$,
%
%
\begin{equation}
\label{eq48.1.2}
|X(N^{3})|_{r,s} - |x|_{L} \le-\tfrac{1}{2}w_{r}N^{2}
\end{equation}
for $\omega\in\mathcal{A}(N^{3})$, and all $r$ and $s$ with $\Delta
_{r}\le 1/b^{3}$ and $s \le N_{H_{r}}$.
\end{prop}

Our main step in demonstrating Proposition \ref{prop48.1.1} will be to
demonstrate the following proposition.
\begin{prop}
\label{prop46.1.1}
Under the same assumptions as in Proposition \ref{prop48.1.1},
%
%
\begin{eqnarray}
\label{eq46.1.2}\hspace*{32pt}
\frac{w_{r}(1+as)(Z^{*}_{r}(N^{3},s) -
z^{*}_{r}(i_{r}(s)))}{\nu_{r} \Gamma (\sigma_{r})} &\le& |x|_{L}
\cdot\frac{1+as}{1+ai_{r}(s)} \cdot \frac{\Gamma (\sigma_{r}) -
\Gamma(\sigma_{r}')}{\Gamma(\sigma_{r})}\nonumber\\[-8pt]\\[-8pt]
&&{} + \frac{C_{30}w_{r}N^{3}}{ab(1+as)} + C_{31}w_{r}N^{3/2}\nonumber
\end{eqnarray}
for appropriate $C_{30}$ and $C_{31}$ not depending on $w, N, a, b, r$
or $s$.
\end{prop}

In order to demonstrate Proposition \ref{prop46.1.1}, we note that, on
account of Proposition~\ref{prop44.1.2}, the LHS of (\ref{eq46.1.2})
is, under the assumptions for the latter proposition, at most
%
%
\begin{eqnarray}
\label{eq46.2.3}
&&d_{r}(s) \bigl(h^{*}_{r}(v(s))N^{3}\vee N^{3\eta} \bigr) \nonumber\\
&&\qquad\le
d_{r}(s) \Bigl(\inf_{s^{\prime} \in[s,i_{r}(s)]}h^{*}_{r}(s') \Bigr)N^{3}
\\
&&\qquad\quad{} + d_{r}(s) \Bigl(h^{*}_{r}(v(s)) - \inf_{s^{\prime} \in[s,i_{r}(s)]}
h^{*}_{r}(s') \Bigr)N^{3} + d_{r}(s)N^{3\eta},\nonumber
\end{eqnarray}
where
\[
d_{r}(s) = (1+\varepsilon_{5})(1+4/b^{2})w_{r}(1+as)/\Gamma(\sigma_{r}).
\]
We will show in Lemmas \ref{lem46.8.2}, \ref{lem46.5.5} and \ref
{lem46.6.3} that each of the three terms on the RHS of (\ref{eq46.2.3})
is bounded above by the corresponding term on the RHS of
(\ref{eq46.1.2}). Proposition \ref {prop46.1.1} then follows.

We first show Lemma \ref{lem46.8.2}, which applies to the first term on
the RHS of (\ref{eq46.2.3}), and should be thought of as the ``main
term'' there.
\begin{lem}
\label{lem46.8.2}
Under the same assumptions as in Proposition \ref{prop48.1.1},
%
%
\begin{equation}
\label{eq46.8.3}
d_{r}(s) \Bigl(\inf_{s^{\prime}
\in[s,i_{r}(s)]}h^{*}_{r}(s') \Bigr) N^{3}
\le|x|_{L}\frac{1+as}{1+ai_{r}(s)} \cdot\frac{\Gamma
(\sigma_{r})-\Gamma (\sigma_{r}')}{\Gamma(\sigma_{r})}.
\end{equation}
\end{lem}
\begin{pf}
It follows from Proposition \ref{prop46.3.2} that
%
%
\begin{equation}
\label{eq46.8.4}
\lambda^{w}(t) \ge(1+\varepsilon_{7}/2)/|x|_{L} \qquad\mbox{for almost all }
t \in [0,N^{3}],
\end{equation}
for large enough $N$ and $b$, and small enough $a$, and therefore
%
%
\begin{equation}
\label{eq46.8.5}
\Delta_{r} \ge(1+\varepsilon_{7}/2)w_{r}N^{3}/|x|_{L} \qquad\mbox{for all } r.
\end{equation}
Consequently, the LHS of (\ref{eq46.8.3}) is at most
%
%
\begin{eqnarray}
\label{eq46.8.6}
&&d_{r}(s)|x|_{L}
\Bigl(\inf_{s^{\prime} \in[s,i_{r}(s)]}h^{*}_{r}(s') \Bigr)
\Delta_{r}/w_{r}(1+\varepsilon_{7}/2) \nonumber\\[-8pt]\\[-8pt]
&&\qquad\le d_{r}(s)|x|_{L} \bigl(\bar{H}^{*}_{r}(s) -
\bar{H}^{*}_{r}(i_{r}(s)) \bigr)/ w_{r}(1+\varepsilon_{7}/2).\nonumber
\end{eqnarray}
This last quantity can be rewritten as
%
%
\begin{eqnarray}
\label{eq46.8.7}
&& \frac{(1+\varepsilon_{5})(1+4/b^{2})}{1+\varepsilon_{7}/2} \cdot|x|_{L}
\cdot\frac{1+as}
{1+ai_{r}(s)} \cdot\frac{\Gamma(\sigma_{r})-\Gamma(\sigma_{r}')}
{\Gamma(\sigma_{r})} \nonumber\\[-8pt]\\[-8pt]
&&\qquad{}\times \frac{1+ai_{r}(s)}{ ({\Gamma
(\sigma_{r})
- \Gamma(\sigma_{r}')})/({\sigma_{r} - \sigma_{r}'} )}.\nonumber
\end{eqnarray}

We proceed to bound the components of (\ref{eq46.8.7}). Since
$\varepsilon_{5}\le\varepsilon_{7}/4$, one has for large enough $b$,
depending on $\varepsilon_{7}$, that
%
%
\begin{equation}
\label{eq46.8.8}
\frac{(1+\varepsilon_{5})(1+4/b^{2})}{1+\varepsilon_{7}/2} \le(1+1/b^{2})^{-1}.
\end{equation}
Since $\Gamma(\cdot)$ is concave and $\sigma_{r}>\sigma_{r}'$,
%
%
\begin{equation}
\label{eq46.8.9}
\frac{\Gamma(\sigma_{r}) - \Gamma(\sigma_{r}')}{\sigma_{r} -
\sigma_{r}'} \ge
\Gamma'(\sigma_{r}) \ge1 + as,
\end{equation}
with the second inequality holding on account of (\ref{eq5.2.2}). So
the last term in (\ref{eq46.8.7}) is at most
%
%
\begin{equation}
\label{eq46.8.10}
\frac{1+ai_{r}(s)}{1+as} = 1 + \frac{a\Delta_{r}}{1+as} \le1 + 1/b^{3},
\end{equation}
where the inequality uses $\Delta_{r} \le1/b^{3}$. Consequently,
(\ref{eq46.8.7}) is, for large $b$, at most
\[
(1+1/b^{2})^{-1}(1 + 1/b^{3})|x|_{L}\frac{1+as}{1+ai_{r}(s)} \cdot
\frac{\Gamma(\sigma_{r}) - \Gamma(\sigma_{r}')}{\Gamma(\sigma_{r})},
\]
which is at most as large as the RHS of (\ref{eq46.8.3}). This implies
the lemma.
\end{pf}

We next demonstrate Lemma \ref{lem46.5.5}, which applies to the second
term on the RHS of (\ref{eq46.2.3}).
\begin{lem}
\label{lem46.5.5}
For all $r$ and $s$ with $\Delta_{r} \le1/b^{3}$ and $s \le N_{H_{r}}$,
%
%
\begin{equation}
\label{eq46.5.6}
d_{r}(s) \Bigl(h^{*}_{r}(v(s)) - \inf_{s^{\prime} \in[s,i_{r}(s)]}
h^{*}_{r}(s') \Bigr)N^{3} \le\frac{C_{30}w_{r}N^{3}}{ab(1+as)}
\end{equation}
for appropriate $C_{30}$ not depending on $w, N, a, b, r$ or $s$.
\end{lem}
\begin{pf}
Since $v(s)-s \le2/b^{3}$, it follows from Proposition \ref
{prop46.5.1} that the LHS of (\ref{eq46.5.6}) is at most
%
%
\begin{equation}
\label{eq46.6.1}
(1+\varepsilon_{5})(1+4/b^{2})\frac{2e b^{2}}{b^{3}}w_{r}N^{3}\bar
{H}^{*}_{r}(s)
\frac{(1+as)}{\Gamma(\sigma_{r})}.
\end{equation}
On account of (\ref{eq5.1.3}), since $\gamma\le\delta_{1}/4$, $b \ge2$
and $\varepsilon_{5}\le1$, this is at most
%
%
\begin{eqnarray}
\label{eq46.6.2}\qquad
\frac{24w_{r}}{C_{2}ab}N^{3}(\bar{H}^{*}_{r}(s))^{1-\gamma}(1+as)
&\le&
\frac{24C_{1}w_{r}}{C_{2}ab}N^{3}(1+as)^{1-(1-\gamma)(2+\delta_{1})}
\nonumber\\[-8pt]\\[-8pt]
&\le&\frac{24 C_{1}w_{r}N^{3}}{C_{2}ab(1+as)} .\nonumber
\end{eqnarray}
Recall that $C_{1}$ and $C_{2}$ do not depend on $w, N, a, b, r$ or
$s$. The RHS of (\ref{eq46.5.6}) follows from this last term by setting
$C_{30}= 24 C_{1}/ C_{2}$.
\end{pf}

We now demonstrate Lemma \ref{lem46.6.3}, which applies to the third
term on the RHS of (\ref{eq46.2.3}).
\begin{lem}
\label{lem46.6.3}
For all $s \le N_{H_{r}}$, 
%
%
\begin{equation}
\label{eq46.6.4}
d_{r}(s)N^{3\eta} \le C_{31}w_{r}N^{3/2}
\end{equation}
for appropriate $C_{31}$ not depending on $w, N, a, b,
r$ or $s$.
\end{lem}
\begin{pf}
Since $s \le N_{H_{r}} \le N$, $\gamma\le1/24$, $\eta\le1/12$, $b \ge2$
and $\varepsilon_{5}\le1$, it follows from (\ref{eq5.1.3}) and
(\ref{eq90.3.1}) that the LHS of (\ref{eq46.6.4}) is at most
%
%
\begin{equation}
\label{eq46.7.1} \frac{4w_{r}N^{3\eta}(1+aN)}{C_{2}a \bar
{H}^{*}_{r}(N_{H_{r}})^{\gamma}} \le
\frac{4}{C_{2}a}w_{r}N^{1/2}(1+aN).
\end{equation}
Since $aN \ge1$, this is at most $8w_{r}N^{3/2}/C_{2}$, which gives the
RHS of (\ref{eq46.6.4}) for $C_{31}= 8/C_{2}$.
\end{pf}

Proposition \ref{prop46.1.1} follows by applying Lemmas
\ref{lem46.8.2}, \ref{lem46.5.5} and \ref{lem46.6.3} to
(\ref{eq46.2.3}).

We will apply the following corollary of the proposition to Proposition
\ref {prop48.1.1}. The corollary combines the inequality
(\ref{eq46.1.2}) with (\ref{eq42.1.2}).
\begin{col}
\label{col46.2.1}
Under the same assumptions as in Propositions \ref{prop48.1.1} and~\ref{prop46.1.1},
%
%
\begin{equation}
\label{eq46.2.2}\qquad |X(N^{3})|_{r,s}-|x|_{L}
\le\frac{C_{30}w_{r}N^{3}}{ab(1+as)} + C_{31} w_{r}N^{3/2} -
\frac{aw_{r}\Delta_{r}z^{*}_{r}(i_{r}(s))}{\nu_{r}\Gamma(\sigma_{r}')}
\end{equation}
for appropriate $C_{30}$ and $C_{31}$ not depending on $w$, $N$, $a$,
$b$, $r$ or $s$.
\end{col}
\begin{pf}
The first term on the RHS of (\ref{eq46.1.2}) of Proposition
\ref{prop46.1.1} is at most
%
%
\begin{equation}
\label{eq46.1.3} |x|_{r,i_{r}(s)}\frac{1+as}{1+ai_{r}(s)}
\cdot\frac{\Gamma(\sigma _{r})-\Gamma
(\sigma_{r}')}{\Gamma(\sigma_{r})} + |x|_{L} - |x|_{r,i_{r}(s)}
\end{equation}
since the coefficients of $|x|_{r,i_{r}(s)}$ in the first term in
(\ref{eq46.1.3}) are at most 1. Substituting (\ref{eq46.1.3}) into
(\ref{eq46.1.2}) and then applying the resulting inequality to the RHS
of (\ref{eq42.1.2}), we note that the term on the LHS of (\ref
{eq46.1.2}) is the first term on the RHS of (\ref{eq42.1.2}) and the
first term in (\ref {eq46.1.3}) is the negative of the second term on
the RHS of (\ref{eq42.1.2}). After the resulting cancellation, the last
two terms on the RHS of (\ref{eq46.1.2}), together with the last term
on the RHS of (\ref{eq42.1.2}), give the RHS of (\ref{eq46.2.2}).
\end{pf}

In order to show Proposition \ref{prop48.1.1}, we will need a lower
bound on the last term on the RHS of (\ref{eq46.2.2}) and an upper
bound on each of the first two terms. In the following lemma, we obtain
the former. Note that the assumptions in the lemma are those of
Proposition \ref{prop46.3.2}, with the additional assumption that
%
%
\begin{equation}
\label{eq48.1.8} |x|_{r,i_{r}(s)} \ge|x|_{L}/(1+\varepsilon_{7}/2)
\qquad\mbox{for some } s \le N_{H_{r}}.
\end{equation}
\begin{lem}
\label{lem48.1.3} Suppose that (\ref{eq17.1.5}) holds for some
$\varepsilon_{7}\in(0,1]$, that $|x|>N^{6}$ with $|x|_{2}
\le\varepsilon_{8}|x|$, and that (\ref{eq48.1.8}) is satisfied for a
given $s$. Then, for large enough $N$ and $b$, and small enough $a$,
%
%
\begin{equation}
\label{eq48.1.4}
\frac{\Delta_{r}z^{*}_{r}(i_{r}(s))}{\nu_{r}\Gamma(\sigma_{r}')}
\ge\frac{N^{3}} {1+ai_{r}(s)} \qquad\mbox{on }
\omega\in\mathcal{A}_{1}(N^{3}).
\end{equation}
\end{lem}
\begin{pf}
By Proposition \ref{prop46.3.2},
%
%
\begin{equation}
\lambda^{w}(t) \ge(1+\varepsilon_{7}/2)/|x|_{L} \qquad\mbox{for almost
all } t \in [0,N^{3}].
\end{equation}
Consequently,
%
%
\begin{equation}
\label{eq48.2.1}
\Delta_{r} \ge(1+\varepsilon_{7}/2)w_{r}N^{3}/|x|_{L} \qquad\mbox{for all } r.
\end{equation}
It follows from (\ref{eq48.2.1}), (\ref{eq5.3.1}) and (\ref {eq48.1.8})
that
the LHS of (\ref{eq48.1.4}) is at least
%
%
\begin{equation}
\label{eq48.2.2}\hspace*{28pt}
\frac{(1+\varepsilon_{7}/2)w_{r}N^{3}z^{*}_{r}(i_{r}(s))}{\nu
_{r}\Gamma (\sigma_{r}')|x|_{L}} = \frac{(1+\varepsilon
_{7}/2)N^{3}|x|_{r,i_{r}(s)}}{(1+ai_{r}(s))|x|_{L}} \ge
\frac{N^{3}}{1+ai_{r}(s)}.
\end{equation}
\upqed\end{pf}

We now apply Corollary \ref{col46.2.1} and Lemma \ref{lem48.1.3} to
demonstrate Proposition \ref{prop48.1.1}.
\begin{pf*}{Proof of Proposition \protect\ref{prop48.1.1}}
We will consider two cases for a given $s \le N_{H_{r}}$, depending on
whether (\ref{eq48.1.8}) holds. Suppose it does. Then, by Lemma
\ref{lem48.1.3},
%
%
\begin{equation}
\label{eq48.2.3}
\frac{aw_{r}\Delta_{r}z^{*}_{r}(i_{r}(s))}{\nu_{r}\Gamma(\sigma _{r}')}
\ge \frac{aw_{r}N^{3}}{1+ai_{r}(s)},
\end{equation}
which is a lower bound for the third term on the RHS of (\ref{eq46.2.2}).

On the other hand, if one chooses $b \ge8C_{30}/a^{2}$, then, since $a
\le1$ and $\Delta_{r}\le1/b^{3} \le1$, the first term on the RHS of
(\ref {eq46.2.2}) satisfies
%
%
\begin{equation}
\label{eq48.2.6}
\frac{C_{30}w_{r}N^{3}}{ab(1+as)} \le\frac{aw_{r}N^{3}}{4(1+ai_{r}(s))},
\end{equation}
which is 1/4 of the RHS of (\ref{eq48.2.3}). Since $s \le N_{H_{r}}$,
$i_{r}(s) \le N+1$. So, the sum of the first and third terms on the RHS
of (\ref{eq46.2.2}) is, for large $N$, at most
%
%
\begin{equation}
\label{eq48.2.7}
-\frac{3aw_{r}N^{3}}{4(1+ai_{r}(s))} \le-\frac{5}{8}w_{r}N^{2}.
\end{equation}
The second term on the RHS of (\ref{eq46.2.2}) satisfies
\[
C_{31}w_{r}N^{3/2} \le\tfrac{1}{8}w_{r}N^{2}
\]
for large $N$. Combining this with (\ref{eq48.2.7}), one obtains from
Corollary \ref{col46.2.1} that
\[
|X(N^{3})|_{r,s} - |x|_{L} \le- \tfrac{1}{2}w_{r}N^{2},
\]
which implies (\ref{eq48.1.2}) under (\ref{eq48.1.8}).

When (\ref{eq48.1.8}) fails for $s$, one has, for large $N$,
%
%
\begin{eqnarray}
\label{eq48.3.1}
|X(N^{3})|_{r,s} - |x|_{L} &=& \bigl(|X(N^{3})|_{r,s} - |x|_{r,i_{r}(s)}\bigr) - \bigl(|x|_{L}
- |x|_{r,i_{r}(s)}\bigr) \nonumber\\
&\le& |X(N^{3})|_{r,s} - |x|_{r,i_{r}(s)} - \tfrac{1}{4}\varepsilon_{7}|x|_{L}
\nonumber\\[-8pt]\\[-8pt]
&\le&
|X^{A}(N^{3})|_{r,s} - \tfrac{1}{4} C_{24}\varepsilon_{7}|x| \nonumber\\
&\le& C_{26}N^{3} - \tfrac{1}{4} C_{24}\varepsilon_{7}N^{6} \le-N^{5},
\nonumber
\end{eqnarray}
where, in the second inequality, we applied (\ref{eq39.1.7}) and
(\ref{eq7.5.1}), and in the third inequality, we applied (\ref
{eq24.4.7}) of Proposition \ref{prop24.4.1} and $|x|>N^{6}$. This
implies (\ref {eq48.1.2}) when (\ref{eq48.1.8}) fails.
\end{pf*}

\section{Conclusion: Upper bounds on $E_{x}[|X(N^{3})|_{L}]$}\label{sec10}

In the preceding four sections, we obtained upper bounds on
\[
|X(N^{3})|_{r,s}-|x|_{L}\quad\mbox{and}\quad E_{x} [|X(N^{3})|_{L} - |x|_{L};
\mathcal{A}(N^{3})^{c} ]
\]
under various assumptions. In Propositions \ref{prop21.1.1} and \ref{prop21.1.3}
we showed that\break
$P_{x}(\mathcal {A}(N^{3})^{c})$ and the corresponding expectation
$E_{x}[|X(N^{3})|_{L}$; $\mathcal{A}(N^{3})^{c}]$ are small. In
Proposition \ref{prop40.2.3}, we showed that the expected value of
$|X(N^{3})|_{r,s}-|x|_{L}$ is small for $s>N_{H_{r}}$. In Sections
\ref{sec8} and \ref{sec9}, we obtained pathwise estimates on
$\mathcal{A}(N^{3})$ when $s \le N_{H_{r}}$, depending on whether
$\Delta_{r}> 1/b^{3}$ or $\Delta_{r} \le1/b^{3}$. Proposition
\ref{prop41.2.1} gives an upper bound in the former subcase and
Proposition \ref{prop48.1.1} gives an upper bound in the latter
subcase. Except for Propositions \ref{prop21.1.1} and \ref{prop21.1.3},
we assumed that $|x|>N^{6}$; for the different results, we also
required various side conditions.

We tie these results together in Proposition \ref{prop50.1.1} to obtain
inequality (\ref{eq5.4.2}) that was cited earlier. We do this in
several steps, first combining the results for $s \le N_{H_{r}}$, then
combining these with Proposition \ref{prop40.2.3} for $s>N_{H_{r}}$,
and lastly including the bound from Proposition \ref{prop21.1.1} on
$\mathcal{A}(N^{3})^{c}$. The first two steps are done in Proposition
\ref{prop48.3.2}. As elsewhere in the paper, $aN \ge1$ is assumed.
\begin{prop}
\label{prop48.3.2} Suppose that (\ref{eq17.1.5}) holds for some
$\varepsilon_{7}\in(0,1]$ and that $|x| > N^{6}$, with $|x|_{2}
\le\varepsilon_{8}|x|$, where $\varepsilon_{8}$ is specified below
(\ref{eq46.3.1}). Then, for large enough $N$ and $b$, and small enough
$a$,
%
%
\begin{equation}
\label{eq48.3.3}
|X(N^{3})|_{r,s} - |x|_{L} \le-\tfrac{1}{2}w_{r}N^{2}
\end{equation}
for $\omega\in\mathcal{A}(N^{3})$, all $r$, and $s$ with $s \le N_{H_{r}}$.
Moreover, for large enough $N$ and $b$, and small enough $a$,
%
%
\begin{eqnarray}
\label{eq48.4.2}
&&E_{x} [|X(N^{3})|_{L} - |x|_{L}; \mathcal{A}(N^{3}) ]
\nonumber\\[-8pt]\\[-8pt]
&&\qquad\le C_{3}(|x|_{K}/|x|)N^{3}- \biggl(\frac{1}{4}\min_{r} w_{r}
\biggr)N^{2} \bigl(2P_{x}(\mathcal{A}(N^{3})) - 1\bigr).\nonumber
\end{eqnarray}
\end{prop}

Note that the assumptions for the first half of Proposition \ref{prop48.3.2}
are the same as for Proposition \ref{prop48.1.1}, except that the restriction
that $\Delta_{r} \le1/b^{3}$ has been removed.
\begin{pf*}{Proof of Proposition \protect\ref{prop48.3.2}}
Inequality (\ref{eq48.1.2}) in Proposition \ref{prop48.1.1} covers the
case where $\Delta_{r} \le1/b^{3}$; (\ref{eq41.2.2}) of Proposition
\ref {prop41.2.1} covers the case where $\Delta_{r}>1/b^{3}$. Together,
they imply (\ref{eq48.3.3}).

In order to demonstrate (\ref{eq48.4.2}), we partition $\mathcal
{A}(N^{3})$ into $G \cup H$, with
\[
G = \Bigl\{ \omega\dvtx|X(N^{3})|_{L} = {\sup_{r, s > N_{H_{r}}}}
|X(N^{3})|_{r,s} \Bigr\}.
\]
Applying Proposition \ref{prop40.2.3} to this $G$, with
$\varepsilon_{3}= (\min_{r} w_{r})/2$, and applying (\ref{eq48.3.3}) on
$H$, it follows that the LHS of (\ref{eq48.4.2}) equals
%
%
\begin{eqnarray}
\label{eq48.4.3}\qquad
&& E_{x} \Bigl[{\sup_{r, s > N_{H_{r}}}}|X(N^{3})|_{r,s} -
|x|_{L}; G \Bigr] + E_{x} \Bigl[{\sup_{r, s \le N_{H_{r}}}}
|X(N^{3})|_{r,s} - |x|_{L}; H \Bigr] \nonumber\\
&&\qquad \le C_{3}(|x|_{K}/|x|)N^{3}- \biggl(\frac{1}{4}\min_{r}
w_{r} \biggr)N^{2} \bigl(2P_{x}(G) + 2P_{x}(H) - 1\bigr) \\
&&\qquad = C_{3}(|x|_{K}/|x|)N^{3}- \biggl(\frac{1}{4}\min_{r}
w_{r} \biggr)N^{2} \bigl(2P_{x}(\mathcal{A}(N^{3})) - 1\bigr).\nonumber
\end{eqnarray}
This implies (\ref{eq48.4.2}).
\end{pf*}

We now obtain our desired result, Proposition \ref{prop50.1.1}, which gives
upper bounds on $E_{x}[|X(N^{3})|_{L}] - |x|_{L}$. The first part of the
proposition applies to all $x$; the second part requires that $|x|>N^{6}$.
\begin{prop}
\label{prop50.1.1}
Suppose that (\ref{eq17.1.5}) holds for some $\varepsilon_{7}\in(0,1]$.
\begin{enumerate}[(a)]
\item[(a)] For large enough $N$,
%
%
\begin{equation}
\label{eq50.6.1}
E_{x}[|X(N^{3})|_{L}] - |x|_{L} \le C_{3}N^{3} \qquad\mbox{for all $x$.}
\end{equation}
%
%
\item[(b)] For $|x|>N^{6}$, large enough $N$ and $b$, and small enough $a$,
%
%
\begin{equation}
\label{eq50.1.2}
E_{x}[|X(N^{3})|_{L}] - |x|_{L} \le C_{3}(|x|_{K}/|x|)N^{3} - \biggl(\frac{1}{4}
\min_{r} w_{r} \biggr)N^{2}.
\end{equation}
\end{enumerate}
In both parts, $C_{3}$ is an appropriate constant that does not depend
on $x$ or
$N$.
\end{prop}
\begin{pf}
We first show (a). By (\ref{eq39.2.1}) and (\ref{eq24.4.2}) of Proposition
\ref{prop24.4.1},
%
%
\begin{equation}
\label{eq50.2.1}
|X(N^{3})|_{L} - |x|_{L} \le C_{26}N^{3}
\end{equation}
for all $\omega\in\mathcal{A}_{1}(N^{3})$ and appropriate $C_{26}> 0$
not depending on $x$, $N$, or $\omega$. Together with Proposition \ref
{prop21.1.1}, this implies
%
%
\begin{eqnarray}
\label{eq50.6.2}
&&E_{x}[|X(N^{3})|_{L}] - |x|_{L} \nonumber\\
&&\qquad= E_{x}[|X(N^{3})|_{L}; \mathcal{A}(N^{3})] + E_{x}
[|X(N^{3})|_{L}; \mathcal{A}(N^{3})^{c}] - |x|_{L} \\
&&\qquad\le C_{26}N^{3} + N^{3}e^{-C_{10}N^{3\eta}} \le2C_{26}N^{3}\nonumber
\end{eqnarray}
for large enough $N$. For $C_{3}\ge2C_{26}$, this implies (\ref{eq50.6.1}).

For (b), we suppose first that $|x|_{2} \le\varepsilon_{8}|x|$, where
$\varepsilon_{8}$
is given
below (\ref{eq46.3.1}). Then, (\ref{eq48.4.2}) of Proposition \ref
{prop48.3.2},
together with Propositions \ref{prop21.1.1} and \ref{prop21.1.3},
implies that
the LHS of (\ref{eq50.1.2}) is equal to
%
%
\begin{eqnarray}
\label{eq50.1.4}
&& E_{x}[|X(N^{3})|_{L}; \mathcal{A}(N^{3})]+E_{x}[|X(N^{3})|_{L};
\mathcal{A}
(N^{3})^{c}] - |x|_{L} \nonumber\\
&&\qquad \le C_{3}(|x|_{K}/|x|)N^{3} - \biggl(\frac{1}{4}\min_{r}
w_{r} \biggr)N^{2}\bigl(2P_{x}(\mathcal{A}(N^{3}))-1\bigr)
\nonumber\\[-8pt]\\[-8pt]
&&\qquad\quad{} + N^{3}e^{-C_{10}N^{3\eta}} \nonumber\\
&&\qquad \le C_{3}(|x|_{K}/|x|)N^{3} - \biggl(\frac{1}{4}\min_{r}
w_{r} \biggr)N^{2} \nonumber
\end{eqnarray}
for large $N$ and $b$, and small $a$. 
%
%
This implies (\ref{eq50.1.2}) for $|x|_{2} \le\varepsilon_{8}|x|$.

Assume now that $|x|_{2}>\varepsilon_{8}|x|$. Choosing $C_{3}\ge
(2C_{26}+\frac{1}{4}
\min_{r}w_{r} ) /\varepsilon_{8}$, it follows from
(\ref{eq50.6.2}) that, for large $N$,
\begin{eqnarray*}
E_{x}[|X(N^{3})|_{L}] - |x|_{L} &\le& \biggl(C_{3}\varepsilon_{8}- \frac{1}{4}
\min_{r}w_{r} \biggr)N^{3} \\
&\le& C_{3}(|x|_{2}/|x|)N^{3} - \biggl(\frac{1}{4}\min_{r}
w_{r} \biggr)N^{3} \\
&\le& C_{3}(|x|_{K}/|x|)N^{3} - \biggl(\frac{1}{4}\min_{r}
w_{r} \biggr)N^{3}.
\end{eqnarray*}
This implies (\ref{eq50.1.2}) for $|x|_{2} > \varepsilon_{8}|x|$.
\end{pf}

\section*{Acknowledgments}
The author thanks the referees for a detailed reading of the
paper and for helpful comments.

\printaddresses


\begin{thebibliography}{99}

\bibitem{r1}
\begin{bbook}[mr]
\bauthor{\bsnm{Bramson},~\bfnm{Maury}\binits{M.}}
(\byear{2008}).
\btitle{Stability of Queueing Networks}.
\bseries{Lecture Notes in Math.}
\bvolume{1950}.
\bpublisher{Springer}, \baddress{Berlin}.
\bid{mr={2445100}}
\end{bbook}
\endbibitem

\bibitem{r2}
\begin{bincollection}[vtex]
\bauthor{\bsnm{Bonald},~\bfnm{T.}\binits{T.}} \AND
\bauthor{\bsnm{Massouli\'{e}},~\bfnm{L.}\binits{L.}}
(\byear{2001}).
\btitle{Impact of fairness on
  Internet performance}.
In \bbooktitle{Proceedings of ACM Sigmetrics}
\bpages{82--91}.
\bpublisher{ACM}, \baddress{New York}.
\end{bincollection}
\endbibitem

\bibitem{r3}
\begin{bbook}[mr]
\bauthor{\bsnm{Chung},~\bfnm{Kai~Lai}\binits{K.~L.}}
(\byear{1985}).
\btitle{A Course in Probability Theory}, \bedition{2nd} ed.
\bpublisher{Academic Press}, \baddress{New York}.
\end{bbook}
\endbibitem

\bibitem{r4}
\begin{bbook}[mr]
\bauthor{\bsnm{Davis},~\bfnm{M.~H.~A.}\binits{M.~H.~A.}}
(\byear{1993}).
\btitle{Markov Models and Optimization}.
\bseries{Monographs on Statistics and Applied Probability}
\bvolume{49}.
\bpublisher{Chapman \& Hall}, \baddress{London}.
\bid{mr={1283589}}
\end{bbook}
\endbibitem

\bibitem{r5}
\begin{barticle}[auto:SpringerTagBib|2009-01-14|16:51:27]
\bauthor{\bsnm{De~Veciana},~\bfnm{G.}\binits{G.}},
  \bauthor{\bsnm{Lee},~\bfnm{T.~J.}\binits{T.~J.}} \AND
  \bauthor{\bsnm{Konstantopoulos},~\bfnm{T.}\binits{T.}}
  (\byear{2001}).
\btitle{Stability and performance
  analysis of networks supporting elastic services}.
\bjournal{IEEE/ACM Transactions on Networking}
\bvolume{9}
\bpages{2--14}.
\bid{doi={10.1109/90.909020}}
\end{barticle}
\endbibitem

\bibitem{r6}
\begin{barticle}[mr]
\bauthor{\bsnm{Gromoll},~\bfnm{H.~Christian}\binits{H.~C.}} \AND
  \bauthor{\bsnm{Williams},~\bfnm{Ruth~J.}\binits{R.~J.}}
(\byear{2009}).
\btitle{Fluid limits for networks with bandwidth sharing and general document
  size distributions}.
\bjournal{Ann. Appl. Probab.}
\bvolume{19}
\bpages{243--280}.
\bid{doi={10.1214/08-AAP541}, mr={2498678}}%
\end{barticle}%
\endbibitem%

\bibitem{r7}
\begin{barticle}[auto:SpringerTagBib|2009-01-14|16:51:27]
\bauthor{\bsnm{Kang},~\bfnm{W.~N.}\binits{W.~N.}},
  \bauthor{\bsnm{Kelly},~\bfnm{F.~P.}\binits{F.~P.}},
  \bauthor{\bsnm{Lee},~\bfnm{N.~H.}\binits{N.~H.}} \AND
  \bauthor{\bsnm{Williams},~\bfnm{R.~J.}\binits{R.~J.}}
  (\byear{2009}).
\btitle{State space collapse and
  diffusion approximation for a network operating under a fair bandwidth
  sharing policy}.
\bjournal{Ann. Appl. Probab.}
\bvolume{19}
\bpages{1719--1780}.
\end{barticle}
\endbibitem

\bibitem{r8}
\begin{barticle}[mr]
\bauthor{\bsnm{Massouli{\'e}},~\bfnm{Laurent}\binits{L.}}
(\byear{2007}).
\btitle{Structural properties of proportional fairness: Stability and
  insensitivity}.
\bjournal{Ann. Appl. Probab.}
\bvolume{17}
\bpages{809--839}.
\bid{doi={10.1214/105051606000000907}, mr={2326233}}
\end{barticle}
\endbibitem

\bibitem{r9}
\begin{barticle}[auto:SpringerTagBib|2009-01-14|16:51:27]
\bauthor{\bsnm{Massouli\'{e}},~\bfnm{L.}\binits{L.}} \AND
\bauthor{\bsnm{Roberts},~\bfnm{J.}\binits{J.}}
(\byear{2000}).
\btitle{Bandwidth sharing and admission
  control for elastic traffic}.
\bjournal{Telecommunication Systems}
\bvolume{15}
\bpages{185--201}.
\bid{doi={10.1023/A:1019138827659}}
\end{barticle}
\endbibitem

\bibitem{r10}
\begin{bincollection}[mr]
\bauthor{\bsnm{Meyn},~\bfnm{S.~P.}\binits{S.~P.}} \AND
  \bauthor{\bsnm{Tweedie},~\bfnm{R.~L.}\binits{R.~L.}}
(\byear{1993}).
\btitle{Generalized resolvents and {H}arris recurrence of {M}arkov processes}.
In \bbooktitle{Doeblin and Modern Probability ({B}laubeuren, 1991)}.
\bseries{Contemp. Math.}
\bvolume{149}
\bpages{227--250}.
\bpublisher{Amer. Math. Soc.}, \baddress{Providence, RI}.
\bid{mr={1229967}}
\end{bincollection}
\endbibitem

\bibitem{r11}
\begin{bbook}[mr]
\bauthor{\bsnm{Nummelin},~\bfnm{Esa}\binits{E.}}
(\byear{1984}).
\btitle{General Irreducible {M}arkov Chains and Nonnegative Operators}.
\bseries{Cambridge Tracts in Mathematics}
\bvolume{83}.
\bpublisher{Cambridge Univ. Press}, \baddress{Cambridge}.
\bid{mr={776608}}
\end{bbook}
\endbibitem

\bibitem{r12}
\begin{bbook}[mr]
\bauthor{\bsnm{Orey},~\bfnm{Steven}\binits{S.}}
(\byear{1971}).
\btitle{Lecture Notes on Limit Theorems for {M}arkov Chain Transition
  Probabilities}.
\bpublisher{Van Nostrand-Reinhold}, \baddress{London}.
\bid{mr={0324774}}
\end{bbook}
\endbibitem


\end{thebibliography}
\end{document}